\documentclass[10pt,twoside]{article}
\usepackage{latexsym,amsfonts,amssymb,stmaryrd,psfig}
\usepackage{graphicx}
\begin{document}

\renewcommand{\theequation}{\thesection.\arabic{equation}}

\newcounter{nsubsec}
\newcommand{\subsec}[1]{\noindent\bf \thesection.\arabic{nsubsec} #1 \rm\addtocounter{nsubsec}{1}}

\newcommand{\enunt}[2]{\medskip\noindent\bf #1 \it #2 \rm \hfill $\Box$}
\newcommand{\te}[2]{\medskip\noindent\bf #1 \it #2 \rm}
\newcommand{\defi}[1]{\medskip\noindent\bf Definition \it #1 \rm}
\newcommand{\dem}[1]{\medskip\noindent\it Proof. \rm #1 \hfill $\Box$}
\newcommand{\demlung}[2]{\medskip\noindent\it #1 \rm #2 \hfill $\Box$}
\newcommand{\ex}[2]{\medskip\noindent\small\sc  #1\hspace{2mm}\rm #2 \normalsize\medskip}
\newcommand{\exh}[3]{\medskip\noindent\small\sc  #1\hspace{2mm}\rm #2 \\\noindent\footnotesize\sc Hint.\hspace{2mm}\rm #3 \normalsize\medskip}
\newcommand{\exhl}[7]{\medskip\noindent\small\sc  #1\hspace{2mm}\rm #2 \\\noindent\footnotesize\sc Hint.\hspace{2mm}\rm #3\\\medskip\noindent\sc #4 \hspace{2mm}\rm #5 \\\noindent\sc #6 \hspace{2mm}\rm #7 \normalsize\medskip}
\newcommand{\exs}[3]{\medskip\noindent\small\sc  #1\hspace{2mm}\rm #2 \\\noindent\footnotesize\sc Solution.\hspace{2mm}\rm #3 \normalsize\medskip}

\newcounter{nbib}
\newenvironment{bib}{\noindent\begin{list}{[\arabic{nbib}]}{\usecounter{nbib}\setlength{\parsep}{0mm}\setlength{\itemsep}{0mm}\setlength{\leftmargin}{7mm}\setlength{\rightmargin}{0mm}}}{\end{list}}

\newcounter{nlsn}
\newenvironment{lsn}[2]{\smallskip\noindent\begin{list}{\arabic{nlsn}}{\usecounter{nlsn}\setlength{\topsep}{0mm}\setlength{\itemsep}{#1}\setlength{\leftmargin}{#2}\setlength{\rightmargin}{\leftmargin}}}{\end{list}}

\newcounter{ntqftax}
\newenvironment{tqftax}[2]{\smallskip\noindent\begin{list}{\it(A\arabic{ntqftax})}{\usecounter{ntqftax}\setcounter{ntqftax}{0}\setlength{\topsep}{0mm}\setlength{\itemsep}{#1}\setlength{\leftmargin}{#2}\setlength{\rightmargin}{\leftmargin}}}{\end{list}}

\newcounter{nls}
\newenvironment{ls}[3]{\smallskip\noindent\begin{list}{#1}{\usecounter{nls}\setlength{\topsep}{0mm}\setlength{\itemsep}{#2}\setlength{\leftmargin}{#3}\setlength{\rightmargin}{\leftmargin}}}{\end{list}}



\def\endprop{\hfill $\Box$}
\def\cc{{\mathbb C}}
\def\zz{{\mathbb Z}}
\def\nn{{\mathbb N}}
\def\rr{{\mathbb R}}
\def\qq{{\mathbb Q}}
\def\ff{{\mathbb F}}

\def\kk{\bf K\it}
\def\ll{\bf L\it}
\def\lll{\bf L\it}
\def\ttt{\bf T\it}
\def\oo{\bf O\it}
\def\T{\cal T\it}
\def\dn{{\cal D}^n}
\def\an{{\cal A}^n}
\def\am{{\cal A}^m}
\def\aa{{\cal A}}
\def\a{{\mathfrak a}}
\def\dd{{\cal D}}
\def\fr{{\mathfrak R}}
\def\ahat{\hat{\cal A}}
\def\pfi{\varphi}
\def\gg{{\mathfrak g}}
\def\e{\{[e_i]\}_{i=1,\dots,2g+n-1}}
\def\f{\{[f_j]\}_{j=1,\dots,2g+n-1}}
\def\hs{H_1(S^3-F,\zz)}
\def\hf{H_1(F,\zz)}
\def\ao{\stackrel{\circ}{\aa(\emptyset)}}
\def\aoh{\widehat{\stackrel{\circ}{\aa(\emptyset)}}}
\def\agh{\widehat{{\cal A}(\G)}}
\def\cnh{\widehat{{\cal C}(n)}}
\def\on{{\cal O}_n}

\newcommand{\G}{\Gamma}
\newcommand{\cZ}{\check Z}
\newcommand{\hZ}{{\Hat Z}}
\newcommand{\ve}{\varepsilon}
\newcommand{\sgn}{{\text sgn}}


\newcommand{\tresa}{\raisebox{0pt}{
                 \begin{picture}(30,25)(-10,-5)
                 \put(12,-12){\line(-1,1){24}}\put(2,2){\line(1,1){10}}
                 \put(-12,-12){\line(1,1){10}}\put(-12,12){i}\put(12,12){i+1}
                 \end{picture}}}
\newcommand{\tresaUnu}{\raisebox{0pt}{
                 \begin{picture}(30,20)(-10,-5)
                 \put(12,-12){\line(-1,1){24}}\put(2,2){\line(1,1){10}}
                 \put(-12,-12){\line(1,1){10}}
                 \end{picture}}}
\newcommand{\tresaDoi}{\raisebox{0pt}{
                 \begin{picture}(30,20)(-10,-5)
                 \put(12,-12){\line(-1,1){10}}\put(-2,2){\line(-1,1){10}}
                 \put(-12,-12){\line(1,1){24}}
                 \end{picture}}}
\newcommand{\vertical}{\raisebox{0pt}{
                 \begin{picture}(30,20)(-10,-5)
                 \put(0,-12){\line(0,1){24}}
                 \end{picture}}}
\newcommand{\flip}{\raisebox{0pt}{
                 \begin{picture}(70,20)(-15,-5)
                 \put(-12,-12){\vector(1,1){24}}\put(-2,2){\vector(-1,1){10}}
                 \put(12,-12){\line(-1,1){10}}
                 \put(15,0){$\leftrightarrow$}
                 \put(52,-12){\vector(-1,1){24}}\put(42,2){\vector(1,1){10}}
                 \put(28,-12){\line(1,1){10}}
                 \end{picture}}}

\newcommand{\tangleOne}{\raisebox{0pt}{
                 \begin{picture}(28,17)(-10,-6)
                 \put(12,-10){\line(-1,1){10}}\put(-8,-10){\line(1,1){10}}
                 \put(2,0){\line(0,1){10}}
                 \end{picture}}}
\newcommand{\tangleTwo}{\raisebox{0pt}{
                 \begin{picture}(28,17)(-10,-6)
                 \put(12,10){\line(-1,-1){10}}\put(-8,10){\line(1,-1){10}}
                 \put(2,0){\line(0,-1){10}}
                 \end{picture}}}
\newcommand{\tangleThree}{\raisebox{0pt}{
                 \begin{picture}(28,17)(-8,-6)
                 \put(12,-10){\line(0,1){16}}\put(-8,-10){\line(0,1){16}}
                 \put(10,6){\line(-1,-1){16}}
                 \end{picture}}}
\newcommand{\tangleFour}{\raisebox{0pt}{
                 \begin{picture}(28,17)(-8,-6)
                 \put(12,-10){\line(0,1){16}}\put(-8,-10){\line(0,1){16}}
                 \put(-6,6){\line(1,-1){16}}
                 \end{picture}}}
\newcommand{\tangleFive}{\raisebox{0pt}{
                 \begin{picture}(32,20)(-8,-6)
                 \put(22,8){\line(-1,-1){10}}\put(2,8){\line(1,-1){10}}
                 \put(12,-2){\line(0,-1){10}}
                 \put(-3,8){\oval(10,10)[t]}\put(-8,8){\line(0,-1){20}}
                 \end{picture}}}
\newcommand{\tangleSix}{\raisebox{0pt}{
                 \begin{picture}(32,20)(-8,-6)
                 \put(12,8){\line(-1,-1){10}}\put(-8,8){\line(1,-1){10}}
                 \put(2,-2){\line(0,-1){10}}
                 \put(17,8){\oval(10,10)[t]}\put(22,8){\line(0,-1){20}}
                 \end{picture}}}
\newcommand{\tangleSeven}{\raisebox{0pt}{
                 \begin{picture}(22,17)(-8,-6)
                 \put(-3,8){\oval(10,10)[t]}\put(-8,8){\line(0,-1){20}}
                 \put(2,8){\line(0,-1){20}}\put(7,13){\line(0,-1){25}}
                 \end{picture}}}
\newcommand{\tangleEight}{\raisebox{0pt}{
                 \begin{picture}(75,0)(-6,-4)
                 \thicklines
                 \put(4,4){\oval(24,24)[b]}
                 \put(56,4){\oval(24,24)[b]}
                 \put(-8,4){\line(1,0){76}}
				 \put(25,-2){$\dots$}
				 \put(15,1){\vector(0,1){3}}
				 \put(67,1){\vector(0,1){3}}
                 \end{picture}}}
\newcommand{\tangleNine}{\raisebox{0pt}{
                 \begin{picture}(60,0)(-6,-4)
                 \thicklines
                 \put(-8,4){\line(0,-1){6}}
                 \put(4,4){\line(0,-1){6}}
                 \put(40,4){\line(0,-1){6}}
                 \put(52,4){\line(0,-1){6}}
                 \put(-8,4){\line(1,0){60}}
				 \put(17,-2){$\dots$}
                 \end{picture}}}

\newcommand{\diaOne}{\raisebox{0pt}{
                 \begin{picture}(30,40)(-10,-5)
                 \put(-8,10){\line(1,-1){10}}\put(12,10){\line(-1,-1){10}}\put(2,0){\line(0,-1){10}}
                 \put(-11,19){\line(0,1){10}}\put(15,19){\line(0,1){10}}\put(2,-20){\line(0,-1){10}}
                 \put(-15,11){$a$}\put(14,10){$b$}\put(0,-17){$c$}
                 \put(-12,14){\circle{10}}\put(16,14){\circle{10}}\put(2,-15){\circle{10}}
                 \end{picture}}}
\newcommand{\diaTwo}{\raisebox{0pt}{
                 \begin{picture}(30,25)(-10,0)
                 \put(-4,15){\oval(20,30)[bl]}\put(10,15){\oval(20,30)[br]}
                 \put(-2,-5){$a^2$}\put(3,-1){\circle{14}}
                 \end{picture}}}
\newcommand{\diaThree}{\raisebox{0pt}{
                 \begin{picture}(45,50)(-15,-15)
                 \put(-8,10){\line(1,-1){10}}\put(12,10){\line(-1,-1){10}}\put(2,0){\line(0,-1){10}}
                 \put(-11,25){\line(0,1){10}}\put(15,25){\line(0,1){10}}\put(2,-34){\line(0,-1){10}}
                 \put(-21,14){$\sqrt[4]{\nu}$}\put(7,14){$\sqrt[4]{\nu}$}\put(-11,-27){$\sqrt[4]{\nu^{-1}}$}
                 \put(-12,17){\circle{16}}\put(16,17){\circle{16}}\put(2,-22){\circle{24}}
                 \end{picture}}}
\newcommand{\diaFour}{\raisebox{0pt}{
                 \begin{picture}(45,55)(-20,-10)
                 \put(-8,-10){\line(1,1){10}}\put(12,-10){\line(-1,1){10}}\put(2,0){\line(0,1){10}}
                 \put(-11,-25){\line(0,-1){10}}\put(15,-25){\line(0,-1){10}}\put(2,34){\line(0,1){10}}
                 \put(-21,-20){$\sqrt[4]{\nu}$}\put(7,-20){$\sqrt[4]{\nu}$}\put(-11,17){$\sqrt[4]{\nu^{-1}}$}
                 \put(-12,-17){\circle{16}}\put(16,-17){\circle{16}}\put(2,22){\circle{24}}
                 \end{picture}}}
\newcommand{\diaFive}{\raisebox{0pt}{
                 \begin{picture}(60,60)(-20,-10)
                 \put(2,-25){\line(0,-1){10}}\put(20,-25){\line(0,-1){10}}
                 \put(20,29){\line(0,1){10}}\put(20,5){\line(0,-1){10}}
                 \put(2,18){\line(0,-1){23}}\put(-7,17){\oval(18,18)[t]}\put(-16,18){\line(0,-1){52}}
                 \put(-2,-18){$\Delta(\sqrt{\nu})$}\put(7,12){$\sqrt{\nu^{-1}}$}
                 \put(12,-15){\oval(40,20)}\put(20,17){\circle{24}}                 
                 \end{picture}}}
\newcommand{\diaSix}{\raisebox{0pt}{
                 \begin{picture}(60,80)(-25,-29)
                 \put(14,-46){\line(-1,1){10}}\put(14,-46){\line(1,1){10}}\put(14,-46){\line(0,-1){10}}
                 \put(24,44){\line(0,1){10}}\put(24,20){\line(0,-1){10}}\put(14,-80){\line(0,-1){10}}
                 \put(4,-10){\line(0,-1){10}}\put(24,-10){\line(0,-1){10}} 
                 \put(-25,10){\line(0,-1){100}} 
                 \put(-10,10){\oval(30,20)[t]}
                 \put(11,27){$\sqrt{\nu^{-1}}$}\put(0,-4){$\Delta(\sqrt{\nu})$}
                 \put(-5,-31){$\sqrt[4]{\nu}$}\put(15,-31){$\sqrt[4]{\nu}$}
                 \put(0,-72){$\sqrt[4]{\nu^{-1}}$}
                 \put(24,32){\circle{22}}\put(14,0){\oval(40,20)}
                 \put(4,-28){\circle{16}}\put(24,-28){\circle{16}}
                 \put(14,-68){\circle{24}}
                 \end{picture}}}
\newcommand{\diaSeven}{\raisebox{0pt}{
                 \begin{picture}(60,70)(-25,-28)
                 \put(14,-29){\line(-1,1){10}}\put(14,-29){\line(1,1){10}}\put(14,-29){\line(0,-1){10}}
                 \put(24,35){\line(0,1){10}}\put(24,11){\line(0,-1){10}}\put(14,-64){\line(0,-1){10}} 
                 \put(-25,1){\line(0,-1){25}}\put(-10,1){\oval(30,20)[t]}\put(-25,-40){\line(0,-1){34}}
                 \put(11,18){$\sqrt[4]{\nu^{-1}}$}\put(0,-13){$\Delta(\sqrt{\nu})$}
                 \put(-34,-35){$\sqrt[4]{\nu}$}\put(0,-56){$\sqrt[4]{\nu^{-1}}$}
                 \put(24,23){\circle{22}}\put(14,-9){\oval(40,20)}
                 \put(-25,-32){\circle{16}}\put(14,-52){\circle{24}}
                 \end{picture}}}
\newcommand{\diaEight}{\raisebox{0pt}{
                 \begin{picture}(60,70)(-25,-20)
                 \put(-3,10){\line(1,-1){10}}\put(17,10){\line(-1,-1){10}}\put(7,0){\line(0,-1){10}}
                 \put(23,31){\line(0,1){10}}\put(7,-26){\line(0,-1){10}}\put(7,-60){\line(0,-1){10}}
                 \put(-10,10){\oval(14,14)[t]}\put(-17,10){\line(0,-1){22}}\put(-17,-28){\line(0,-1){42}}
                 \put(-25,-24){$\sqrt[4]{\nu}$}\put(12,14){$\sqrt[4]{\nu^{-1}}$}
                 \put(-7,-53){$\sqrt[4]{\nu^{-1}}$}\put(-1,-22){$\sqrt{\nu}$}
                 \put(-17,-20){\circle{16}}\put(25,19){\circle{24}}
                 \put(7,-48){\circle{24}}\put(7,-18){\circle{16}}
                 \end{picture}}}
\newcommand{\diaNine}{\raisebox{0pt}{
                 \begin{picture}(84,0)(-12,-2)
                 \put(-4,15){\oval(20,30)[bl]}\put(10,15){\oval(20,30)[br]}
                 \put(0,-3){$\nu$}\put(3,-1){\circle{14}}\put(20,10){\vector(0,1){5}}
                 \put(46,15){\oval(20,30)[bl]}\put(60,15){\oval(20,30)[br]}
                 \put(50,-3){$\nu$}\put(53,-1){\circle{14}}\put(70,10){\vector(0,1){5}}
                 \put(22,0){$\dots$}
                 \end{picture}}}

\newcommand{\strutt}{\raisebox{0pt}{
                 \begin{picture}(14,6)(0,0)
                 \put(4,4){\oval(16,8)[b]}
                 \end{picture}}}
\newcommand{\UP}{\raisebox{0pt}{
                 \begin{picture}(14,8)(-5,-4)
                 \thicklines
                 \put(0,2){\oval(16,16)[b]}
                 \end{picture}}}
\newcommand{\LP}{\raisebox{0pt}{
                 \begin{picture}(14,8)(-5,-4)
                 \thicklines
                 \put(0,-4){\oval(16,16)[t]}
                 \end{picture}}}
\newcommand{\U}{\raisebox{0pt}{
                 \begin{picture}(16,8)(-5,-4)
                 \thicklines
                 \put(0,2){\oval(16,16)[b]}    
                 \put(8,3){\vector(0,1){3}}
                 \end{picture}}}
\newcommand{\D}{\raisebox{0pt}{
                 \begin{picture}(16,8)(-5,-4)
                 \thicklines
                 \put(0,-3){\oval(16,16)[t]}    
                 \put(-8,-4){\vector(0,-1){3}}
                 \end{picture}}}
\newcommand{\zerochords}{\raisebox{0pt}{
                 \begin{picture}(30,22)(-12,-7)
                 \thicklines\put(0,0){\circle{20}}
                 \end{picture}}}
\newcommand{\dashedO}{\raisebox{0pt}{
                 \begin{picture}(12,8)(-5,-4)
                 \thinlines\put(0,0){\circle{10}}
                 \end{picture}}}
\newcommand{\gGraph}{\mathop{\underbrace{\raisebox{0pt}{
                 \begin{picture}(78,5)(-40,-5)
                 \thicklines
                 \put(-38,0){\circle{10}}
                 \put(-33,0){\vector(1,0){10}}
                 \put(-18,0){\circle{10}}
                 \put(-13,0){\vector(1,0){10}}
                 \put(-38,5){\vector(-1,0){3}}
                 \put(-18,5){\vector(-1,0){3}}
                 \put(-38,-5){\vector(1,0){3}}
                 \put(-18,-5){\vector(1,0){3}}
                  \dots
                 \put(0,0){\vector(1,0){10}}
                 \put(15,0){\circle{10}}
                 \put(15,5){\vector(-1,0){3}}
                 \put(15,-5){\vector(1,0){3}}
                 \end{picture}}}_{g\; times }}}
\newcommand{\oGraph}{\raisebox{0pt}{
                 \begin{picture}(76,10)(-42,-5)
                 \thicklines
                 \put(-38,0){\circle{10}}
                 \put(-33,0){\vector(1,0){10}}
                 \put(-18,0){\circle{10}}
                 \put(-13,0){\vector(1,0){10}}
                 \put(-38,5){\vector(-1,0){3}}
                 \put(-18,5){\vector(-1,0){3}}
                 \put(-38,-5){\vector(1,0){3}}
                 \put(-18,-5){\vector(1,0){3}}
                 \dots
                 \put(0,0){\vector(1,0){10}}
                 \put(15,0){\circle{10}}
                 \put(15,5){\vector(-1,0){3}}
                 \put(15,-5){\vector(1,0){3}}
                  \end{picture}}}
\newcommand{\oOriented}{\raisebox{0pt}{
                 \begin{picture}(12,8)(-5,-4)
                 \thicklines
                 \put(0,0){\circle{10}}
                 \put(0,5){\vector(-1,0){3}}
                 \put(0,-5){\vector(1,0){3}}
                 \end{picture}}}
\newcommand{\arrowX}{\mathop{\raisebox{0pt}{
                 \begin{picture}(2,10)(0,-5)
                 \thicklines
                 \put(0,-5){\vector(0,1){10}}
                 \end{picture}}}_X}
\newcommand{\arrowY}{\mathop{\raisebox{0pt}{
                 \begin{picture}(2,10)(0,-5)
                 \thicklines
                 \put(0,-5){\vector(0,1){10}}
                 \end{picture}}}_Y}
\newcommand{\arrowZ}{\mathop{\raisebox{0pt}{
                 \begin{picture}(2,10)(0,-5)
                 \thicklines
                 \put(0,-5){\vector(0,1){10}}
                 \end{picture}}}_Z}
\newcommand{\arrowtwog}{\mathop{\raisebox{0pt}{
                 \begin{picture}(6,10)(0,-5)
                 \thicklines
                 \put(4,-5){\vector(0,1){10}}\put(0,5){\vector(0,-1){10}}
                 \end{picture}}}_{[2g]}}
\newcommand{\arrowdown}{\mathop{\raisebox{0pt}{
                 \begin{picture}(2,10)(0,-5)
                 \thicklines
                 \put(0,2){\vector(0,-1){10}}
                 \end{picture}}}}
\newcommand{\arrowg}{\mathop{\raisebox{0pt}{
                 \begin{picture}(2,10)(0,-5)
                 \thicklines
                 \put(0,-5){\vector(0,1){10}}
                 \end{picture}}}_{[g]}}
\newcommand{\arrowgi}{\mathop{\raisebox{0pt}{
                 \begin{picture}(2,10)(0,-5)
                 \thicklines
                 \put(0,-5){\vector(0,1){10}}
                 \end{picture}}}_{[g_i]}}
\newcommand{\arrowgj}{\mathop{\raisebox{0pt}{
                 \begin{picture}(2,10)(0,-5)
                 \thicklines
                 \put(0,-5){\vector(0,1){10}}
                 \end{picture}}}_{[g_j]}}
\newcommand{\arrowgOne}{\mathop{\raisebox{0pt}{
                 \begin{picture}(2,10)(0,-5)
                 \thicklines
                 \put(0,-5){\vector(0,1){10}}
                 \end{picture}}}_{[g_1]}}
\newcommand{\arrowgTwo}{\mathop{\raisebox{0pt}{
                 \begin{picture}(2,10)(0,-5)
                 \thicklines
                 \put(0,-5){\vector(0,1){10}}
                 \end{picture}}}_{[g_2]}}
\newcommand{\OTHERarrowg}{\uparrow_{[g]}}
\newcommand{\gArrows}{\mathop{\underbrace{\raisebox{0pt}{
                 \begin{picture}(30,5)(-10,-5)
                 \thicklines
                 \put(-10,-7){\vector(0,1){12}}
                 \put(-5,-7){\vector(0,1){12}}\dots
                 \put(3,-7){\vector(0,1){12}}
                 \end{picture}}}_{g\; times }}}

%
%
%
%

\title{\bf A TQFT associated to the LMO invariant of three-dimensional manifolds
\thanks{The 
results of this article were obtained when the authors were at the Department of Mathematics, 
State University of New York at Buffalo, Buffalo, NY 14260-2900, USA}
\author{Dorin Cheptea and Thang T Q Le}
}
\date{}
\maketitle
\pagestyle{myheadings}
\markboth{Dorin Cheptea and Thang T Q Le}{A TQFT associated to the LMO invariant of three-dimensional manifolds}


\begin{abstract}
We construct a Topological Quantum Field Theory (in the sense of Atiyah [\ref{at88}]) associated 
to the universal finite-type invariant of 3-dimensional manifolds, as a functor from the category 
of 3-dimensional manifolds with parametrized boundary, satisfying some additional conditions, 
to an algebraic-combinatorial category. 
It is built together with its truncations with respect to a natural grading, and
we prove that these TQFTs are non-degenerate and anomaly-free.
The TQFT(s) induce(s) a (series of) representation(s) of
a subgroup ${\cal L}_g$ of the Mapping Class Group that contains the Torelli group.
The $N=1$ truncation produces a TQFT for the Casson-Walker-Lescop invariant.
\end{abstract}


\setcounter{nsubsec}{1}
\noindent In [\ref{LMO}] the Kontsevich integral $Z(L)$,
i.e. the universal finite-type invariant for links,
has been extended to an invariant $Z^{LMO}(M)$ of 3-dimensional manifolds.
The later is universal for rational homology 3-spheres [\ref{habiro}]. 
The task of putting $Z^{LMO}$ in the structural
framework of TQFT was partially accomplished in [\ref{MO}], 
however that construction uses a twisted gluing of cobordisms
and the resulting anomaly is complicated. The construction is however 
important to establish, since TQFT natually connects a manifold
invariant (in this case LMO) to the Mapping Class Group.
It also aims to shed some light to the question of topological interpretation of
quantum invariants of manifolds.

Our construction allows us to associate to the LMO invariant 
an infinite-dimensional {\it linear} representation of the Torelli group,
in fact of a larger {\it Lagrangian subgroup} of the Mapping Class Group. The new results of this
paper are:
\begin{itemize}
\item proving an isomorphism (Proposition 2.3) redusing the study of the LMO invariant of 3-dimensional
manifols with parametrized boundary to that of finite-type invariants of string-links
\item the construction of the composition of chord diagrams from truncations (Theorem 2.12) and establishing an 
important limit property (Lemma 3.5)
\item proving the non-degeneracy of the TQFT (Theorem 3.2), which means it is posible to calculate the induced
representation in purely combinatorial terms 
\item a combinatiral description of the map $Z(M)\mapsto Z(\widehat{M})$ (Proposition 4.3), where $M$ denotes
the cobordism, and $\widehat{M}$ the closed 3-manifold obtained by "caping" its boundary
\end{itemize}
The natural truncation induces a TQFT for the Walker-Lescop extension
of the Casson invariant, and we can identify Morita's representation as
its first non-trivial part.

The essential difference between this TQFT and the
Reshetikhin-Turaev TQFT for quantum invariants is that
the present is taylored for integer and rational homology spheres, 
because  $Z^{LMO}$ is strong for them, and weaker if the rank of homology is bigger.
Hence we consider {\it connected} cobordisms between {\it connected}
surfaces. When gluing, we discriminate between the domain and the range of a cobordism. 
In particular, while we regard the
standard surface $\Sigma_g$ of genus $g$ in the domain as the boundary of the
standard handlebody $N_g$, we regard $\Sigma_g$ in the range as the boundary
of the complement of $N_g$ in $S^3$. Thus, gluing identically in our TQFT produces
$S^3$ as opposed to $\#_g(S^2\times S^1)$ in the Reshetikhin-Turaev TQFT.
It remains, however, to interpret our TQFT as a "perturbative expansion around $0$" 
of the Reshetikhin-Turaev TQFT [\ref{C}].

This paper is organized as follows.
In Section 1 we recall the topological categories $\mathfrak Q, 
\mathfrak Z$ introduced in [\ref{CL}], and the pertaining results that 
we will need subsequently in this paper. 
The categories $\aa$ and $\aa^{\leq N}$ of chord diagrams are
explained in Section 2. We use a simplier definition of $Z$ on elementary 
pseudo-quasi-tangles\footnote{"pseudo" stands for the presence of
3-valent vertices}, and an even associator, while the one in [\ref{MO}] is based on 
the Knizhnik-Zamolodchikov associator.
Important results here are Proposition 2.3 and Theorem 2.12. 
We also recall for comparison the multiplication $\bullet$ of $\aa(\arrowg)$ 
and its Campbell-Hausdorff-type property. 
In Section 3 we formally construct the anomaly-free (by [\ref{CL}]) truncated and full TQFTs. 
The main result of that section is showing that the completion of the algebraic image of cobordisms with one 
boundary component is precisely the whole space $\aa(\arrowg)$ of chord diagrams on $g$ vertical
lines. That means that the induced representation can in principle be used for combinatorial calculations 
in solving topological questions about three-dimensional manifolds and the Mapping Class Group. 
We also finish the proof of Theorem 2.12 there.
In Section 4 we restrict to the case $N=1$ to get a TQFT for the 
Casson-Walker-Lescop invariant. Also there we identify (for arbitrary $N$)
the algebraic map that sends the invariant of a manifold with parametrized boundary to
the invariant of the closed manifold obtained from the former by the
natural procedure that we call {\it filling} (see Section 1).


\medskip\noindent
{\bf 0.1 Chord diagrams.}
Let us recall some basic definitions. (For details see
[\ref{BN}, \ref{BLT}].)
An {\em open chord diagram} is a vertex-oriented uni-trivalent
graph, i.e. a graph with univalent and trivalent vertices together
with a cyclic order of the edges incident to the trivalent
vertices. Self-loops and multiple edges are allowed. A univalent
vertex is also called a {\em leg}, and a trivalent vertex is also
called an {\em internal vertex}. In planar pictures, the
orientation of the edges incident to a vertex is the
counterclockwise orientation, unless otherwise stated; the
pictures can not be perfect since not every graph is planar,
therefore when reading pictures one should keep in mind that
four-valent vertices do not exist. The {\em degree} of an open
chord diagram is half the number of all vertices.

Suppose $\Gamma$ is a compact oriented 1-manifold (possibly with
boundary) and $X$ a finite set of asterisks. A {\em chord diagram
with support $\Gamma\cup X$} is a vertex-oriented uni-trivalent
graph $D$ together with a decomposition $D=\Gamma\cup E$, where
$E$ is an open chord diagram with some legs labeled by elements
of $X$, such that $D$ is the result of gluing all non-labeled
legs of $E$ to distinct interior points of $\Gamma$ (the 3-valent
vertices resulted from gluing, which will {\sc not} have an associated a
cyclic order of adjacent edges, are called {\em external
vertices}). Repetition of labels is allowed and not all labels
have to be used. The {\em degree} of $D$ is, by definition, the
degree of $E$. $\Gamma$ is also called the {\em skeleton} of $D$,
and in pictures is represented by bold lines. Often the components
of $\Gamma$, as well as different asterisks in $X$, are
distinguished in pictures by labels.

By {\em a graph} $\Gamma$ we will mean a uni-trivalent graph, with 
all edges oriented, and with a cyclic order of edges incident to 
trivalent vertices prescribed. Self-loops and multiple edges are 
allowed. The connected components of the graph have to be always 
ordered (and this order has to be preserved by a homeomorphism). 
Additionally we may label (color) some subgraphs within each 
connected component. One should think of a graph as a generalization 
of the notion of oriented compact 1-manifold. We can repeat the 
definition of the previous paragraph to obtain the notion of a {\em chord 
diagram with support a graph}. The graph is the skeleton of the chord 
diagram. We keep the same definition of the degree: the vertices 
of the graph are {\sc not} counted. As examples of graphs $\Gamma$ let us 
consider:

\begin{itemize}
\item the oriented manifold which is the union of $g\in\nn^{\ast}$
copies of $[0,1]$, each copy labeled (colored) by a distinct element
of a finite abstract ordered set $X$ of asterisks.
This special graph will be denoted $\arrowX$, and in planar
pictures will be represented by vertical lines $\gArrows$
\item the {\it chain graph}\footnote{this terminology we borrow from [\ref{MO}]}
suggestively denoted $\oGraph$.
The order of the edges adjacent to each vertex is everywhere the
counterclockwise with respect to its standard embedding in
$\rr^2\subset\rr^3$, the subgraphs $\oOriented$ are labeled $1$
through $g$ from left to right. Let us denote this special graph
by $\Gamma^g$. Note that $\Gamma^g$ standardly embedded in $\rr^3$
has a preferred (the blackboard) ribbon graph neighbourhood, in the
sense of section 1 below. For $g=1$, set $\Gamma^1=\oOriented$,
one oriented edge, no vertices.
\item it is convenient to set $\Gamma^0 = one\: point$ as a chain 
graph. Chord diagrams on $\Gamma^0$ automatically can have only 
internal vertices.
\end{itemize}

%
%
%
%

\section{The topological category}
\setcounter{equation}{0}
\setcounter{nsubsec}{1}


\subsec{Definitions {\rm (see also [\ref{CL}])}.}Two 
{\it triplets} $(K,G_1,G_2)$ and
$(L,H_1,H_2)$, consisting of a framed oriented link $K$
(respectively $L$) in $S^3$, and two disjoint (and disjoint from the
corresponding link) embedded {\sc framed} chain graphs are {\it equivalent} 
(notation $\cong$) if there is a PL-homeomorphism 
$\phi:S^3\rightarrow S^3$ which preserves the link and the embedded framed graphs, 
i.e. $\phi$ sends $K$ to $L$, the first embedded framed graph $G_1=(\Gamma_1, R_1)$ 
to the first embedded framed graph $H_1=(\Delta_1, S_1)$, and the second 
$G_2=(\Gamma_2, R_2)$ to the second one $H_2=(\Delta_2, S_2)$. Here 
$\emptyset$ is also considered a framed oriented link in $S^3$. Call 
$G_1=(\Gamma_1, R_1)$ the {\it bottom}, and $G_2=(\Gamma_2, R_2)$ 
the {\it top} of the triplet.

Let $M$ be a compact oriented $3$-manifold with boundary
$\partial M = (-S_1)\cup S_2$, suppose that parametrizations
$f_i:\Sigma_{g_i}\rightarrow S_i$ are fixed. Such 
$(M,f_1,f_2)$ is called a (parametrized) {\it (2+1)-cobordism}. 
$S_1$ is {\it the bottom}, $S_2$ -- {\it the top} of 
the cobordism. The cobordisms  $(M,f_1,f_2)$ 
and $(N,h_1,h_2)$ are {\it equivalent (homeomorphic)} if there is a 
PL-homeomorphism $F:M\rightarrow N$ sending bottom to 
bottom and top to top preserving the parametrizations, 
i.e. $F\circ f_i=h_i$, $i=1,2$.

Fix $N_g$ -- a standard neighbourhood of $\Gamma^g$ in $S^3$.
$\Sigma_g=\partial N_g\subset S^3$ is the standard oriented
surface of genus $g$. Let $\overline{N_g}$ be the handlebody
complement of $N_g$ in $S^3$.  We also denote by
$\overline{\Gamma^g}$ the core of $\overline{N_g}$. Clearly
$\partial\overline{N_g}=-\Sigma_g$. When $g=0$ we assume that
$\Gamma^g$ is a point, and $N_g$ is a ball. 
Using the parametrizations, we can glue the standard handlebody
$N_{g_1}$ to the bottom and the standard anti-handlebody
$\overline{N_{g_2}}$ to the top of $M$. The result
$\widehat{M}:=M\cup_{f_1}N_{g_1}\cup_{-f_2}(-\overline{N_{g_2}})$ 
is called {\em the filling} of $(M,f_1,f_2)$.


\medskip
\subsec{Surgery description of gluing cobordisms.}Let 
$\mathfrak G$ denote set of equivalence classes of triplets 
$(L,G_1,G_2)$ in $S^3$. Let $\mathfrak C$ denote the 
set of equivalence classes of  3-cobordisms,
with connected non-empty bottom and top.

\enunt{Proposition.}{1) The map $\kappa:{\mathfrak G}
\rightarrow\mathfrak C$ that associates to every equivalence class of triplets 
$(L,G_1,G_2)$ the equivalence class of cobordisms $(M,f_1,f_2)$, 
obtained by doing surgery on $L\subset S^3$, removing tubular
neighbourhoods $N_1,N_2$ of each $G_1,G_2$, and recording 
the parametrizations of the two obtained boundary components, is well-defined and surjective.
If one glues according to these parametrizations a standard handlebody to 
$-\partial N_1$ and a standard anti-handlebody to $\partial N_2$,
then one obtains $S^3_L$.

2) Let a {\bf first Kirby move} on a triplet be the cancellation / insertion of a $\oo^{\pm 1}$ 
separated by an $S^2$ from anything else, and an {\bf extended (generalized) second Kirby move}
be a slide over a link component of an arc, either from another link component or from a chain graph.
Then if one factors $\mathfrak G$ by the extended Kirby moves and changes of orientations of link components, 
the induced map $\overline{\kappa}$ is a bijection.
}

\medskip
For example, to represent the 3-cobordism $(\Sigma_g\times[0,1], p_1, p_2)$, where 
$p_i:\Sigma_g\rightarrow\Sigma_g\times i\subset S^3$ are two coppies of the standard embedding of
$\Sigma_g$ in $S^3$, we can choose the framed graphs $R_1$, $R_2$ as in Figure 1.
Let $\Gamma$, respectively $\Gamma^{\prime}$ generically denote
the bottom, respectively the top of a triplet.
Call the union of the lower half-circles and the horizontal
segments of $\Gamma$, {\em the horizontal line} of $\Gamma$.
Similarly, call the union of the upper half-circles and the
horizontal segments of $\Gamma^{\prime}$, {\em the horizontal
line} of $\Gamma^{\prime}$. See Figure 2a.

\begin{figure}[thb]
\centerline{\psfig{file=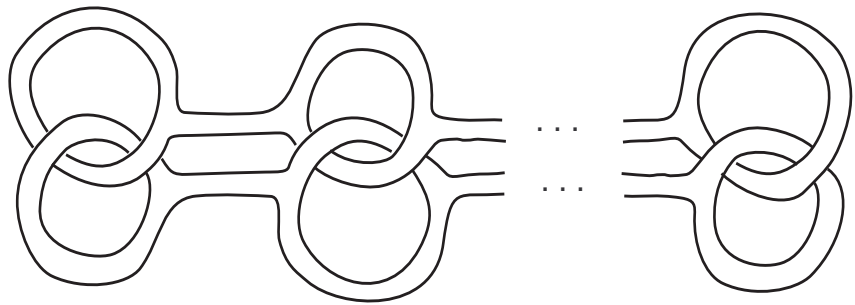,width=12cm,height=5cm,angle=0}}
\caption{\sl The preferred choice of ribbons $R_i$, $i=1,2$ for $(\Sigma_g\times I,p_1,p_2)$.} \label{fig2}
\end{figure}

\begin{figure}[thb]
\centerline{\psfig{file=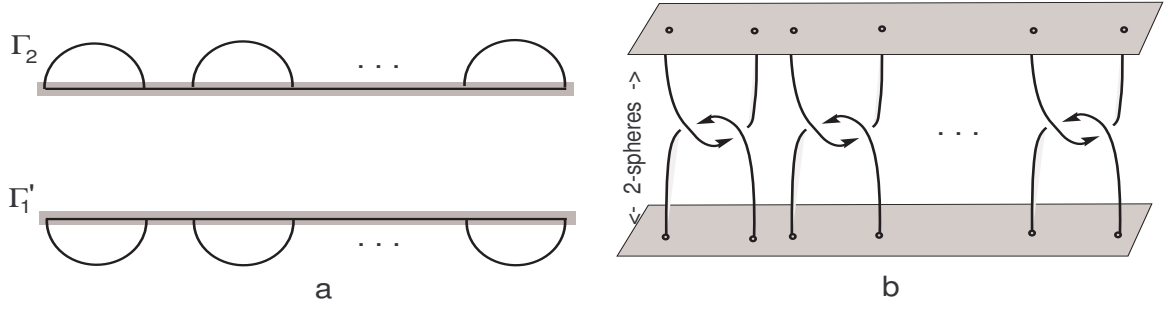,width=16cm,height=5cm,angle=0}}
\caption{\sl {\bf a:} The horizontal segments of $\Gamma_1^{\prime}$ 
and $\Gamma_2$; {\bf b:} Framed tangle
$T_g\subset\overline{B(0,2)-B(0,1)}$.} \label{fig4}
\end{figure}


\medskip
\subsec{Proposition.}{\it Let $(M_1,f_1,f_1^{\prime})$ and
$(M_2,f_2,f_2^{\prime})$ be two
3-cobordisms with connected non-empty bottoms and tops,
represented by triplets $L_1,G_1,G_1^{\prime})$ and
$(L_2,G_2,G_2^{\prime})$.
Remove a 3-ball neighbourhood of the horizontal line of
$G_1^{\prime}\subset S^3$, and identify the remain with 
$\overline{B(0,1)}$. Remove a 3-ball neighbourhood of the
horizontal line of $G_2\subset S^3$, and identify
the remain with $\overline{S^3-B(0,2)}$. Glue the framed tangle
$T_g\subset\overline{B(0,2)-B(0,1)}$ shown in figure \ref{fig4}b to
the ends of the remains of $G_1^{\prime}$ in $\overline{B(0,1)}$
and $G_2$ in $\overline{S^3-B(0,2)}$, strictly
preserving the order of the points, so that the composition of
these framed tangles is a smooth framed oriented link $L_0$ in
$S^3=(\overline{S^3-B(0,2)})\cup(\overline{B(0,2)-B(0,1)})
\cup(\overline{B(0,1)})$. Then 
\begin{equation}
\kappa(L_1\cup L_0\cup L_2, G_1, G_2^{\prime}) 
=  (M_2\cup_{f_2\circ(f_1^{\prime})^{-1}}M_1, f_1, f_2^{\prime})
\label{CompTriplet}
\end{equation}
where the framed graphs $G_1$, $G_2^{\prime}$ 
in this formula are determined in the obvious way by the original $G_1$, 
$G_2^{\prime}$ in the two copies of $S^3$. Hence, any triplet
representing $(M_2\cup_{f_2\circ(f_1^{\prime})^{-1}}M_1, f_1, f_2^{\prime})$
is equivalent to $(L_1\cup L_0\cup L_2, G_1, G_2^{\prime})$ by
extended Kirby moves and changes of orientations of link components.
\endprop}

\medskip
Let $\Gamma^g$ be a chain graph embedded in an arbitrary
way in $S^3$, then $H_1(S^3-\Gamma^g,\zz)\cong\zz^g$, with free
generators the meridians of the circle components.


\medskip
\subsec{Proposition.}{\it Suppose $M$ is a connected compact oriented 
3-manifolds with two distinguished (not necessarily connected) 
boundary components $\partial M = (-S_1)\cup S_2$, let $f_1$, $f_2$ be
parametrizations of these surfaces, let $N_1$ and $\overline{N_2}$ be 
corresponding-to-the-genera disjoint unions of standard handlebodies, 
respectively anti-handlebodies, and let $i=(i_1,i_2): \partial
M\hookrightarrow M$ be the inclusion. The following conditions are
equivalent:

\medskip
\rm(1)\it \quad $ H_1(\widehat{M}, \zz)=0$

\rm(2)\it \quad $ H_1(M, \zz) = i_{\ast}(H_1(\partial M,\zz)/
(f_1,-f_2)_{\ast}H_1(N_1\sqcup -\overline{N_2},\zz))$ 

\medskip\noindent
They imply:

\medskip
\rm(3)\it \quad $ 2\cdot rank\: H_1(M;\zz) = rank\: H_1(\partial
M;\zz)$ \endprop}


\medskip
\subsec{Proposition.}{\it Suppose $M$ is a connected compact oriented 
3-manifolds with two distinguished (not necessarily connected) 
boundary components $\partial M = (-S_1)\cup S_2$, let $f_1$, $f_2$ be
parametrizations of these surfaces, let $N_1$ and $\overline{N_2}$ be 
corresponding-to-the-genera disjoint unions of standard handlebodies, 
respectively anti-handlebodies, and let $i=(i_1,i_2): \partial
M\hookrightarrow M$ be the inclusion. The following conditions are
equivalent:

\medskip
\rm(1)\it \quad $ H_1(\widehat{M}, \qq)=0$

\rm(2)\it \quad $ H_1(M, \qq) = i_{\ast}(H_1(\partial M,\qq)/
(f_1,-f_2)_{\ast}H_1(N_1\sqcup -\overline{N_2},\qq))$

\medskip\noindent
They imply:

\medskip
\rm(3)\it \quad $ 2\cdot rank\: H_1(M;\qq) = rank\: H_1(\partial
M;\qq)$ \endprop}

\medskip
A 3-cobordism satisfying the equivalent conditions (1), (2) of Proposition 1.4
is called an {\it Integer Homology Cobordism} ($\zz$HC). 
A 3-cobordism satisfying the equivalent conditions (1), (2) of Proposition 1.5
is called a {\it Rational Homology Cobordism} ($\qq$HC).
Note that in both definitions of $\qq$HC and $\zz$HC we allow one
or both $S_i$ to be empty, although from the point of this TQFT the case of empty
top and/or bottom is indistinguished from the case when that component is $S^2$.


\medskip
\subsec{Description of the categories $\mathfrak Q\supset\mathfrak Z$.}Objects 
in each of these are {\em natural numbers}. The
morphisms between $g_1$ and $g_2$ are equivalence
(homeomorphism) classes of connected 3-cobordisms
with bottom $S_1$ of genus $g_1$ and top $S_2$ of genus $g_2$,
satisfying the $\ff$-  semi-Lagrangian conditions:
\begin{equation}\label{LagrangeCondition}
f_{1\ast}L^a \supseteq f_{2\ast}L^a\; {\rm and} \; f_{1\ast}L^b \subseteq f_{2\ast}L^b,
\end{equation}
where  $L^a={\rm ker}({ \rm incl}_{\ast}:H_1(\Sigma_g,\ff)\rightarrow H_1(N_g,\ff))$,
and $L^b={\rm ker}({ \rm incl}_{\ast}:H_1(\Sigma_g,\ff)\rightarrow H_1(\overline{N_g},\ff))$, and
$\ff=\zz$ or $\qq$.
The composition-morphism of two cobordisms $(M_1,f_1,f_1^{\prime})$ and $(M_2,f_2,f_2^{\prime})$
is the equivalence class of the 3-cobordism $(M_2\cup_{f_2\circ(f_1^{\prime})^{-1}}M_1,f_1,f_2^{\prime})$.

In general condition (\ref{LagrangeCondition})
over $\zz$ is stronger than (\ref{LagrangeCondition}) over $\qq$. 
It also may hold with strict inclusion. [\ref{CL}]


\medskip
\subsec{Proposition {\rm [\ref{CL}]}.}{\it The composition of two morphisms (say, class of $M$ and class of $N$) in 
category $\mathfrak Q$ (respectively $\mathfrak Z$) is again a 
morphism in the category $\mathfrak Q$ (respectively $\mathfrak Z$). \endprop}

\medskip
Let us restrict to 3-cobordisms
$M$ of the form $(\Sigma_g\times[0,1],f\times 0,f^{\prime}\times 1)$, with $f,f^{\prime}\in Aut(\Sigma_g)$, i.e.
the parametrization of the top differs by that of the bottom by the
automorphism $w=(f^{\prime})^{-1}\circ f$. The equivalence
classes of this cobordism depends only on the isotopy class of $w$ (i.e. we don't need
to specify both $f$, $f^{\prime}$). The equivalence class of $M=(\Sigma_g\times[0,1],f\times 0,f^{\prime}\times 1)$ 
is a $\zz$-semi-Lagrangian cobordism iff it is a $\qq$-semi-Lagrangian cobordism iff it satisfies 
$L^a = w_*(L^a)$ and $L^b = w_*(L^b)$, i.e. is a Lagrangian cobordism.
Then $\widehat{M}$ is always a $\zz$-homology sphere. [\ref{CL}]

The composition of two cobordisms
$(\Sigma_{g}\times I, f_1\times 0, f_1^{\prime}\times
1)\cong(\Sigma_{g}\times I, w_1\times 0, id\times1)$ and
$(\Sigma_{g}\times I, f_2\times 0, f_2^{\prime}\times
1)\cong(\Sigma_{g}\times I, w_2\times 0, id\times1)$ along
$(f_2\times 0)\circ((f_1^{\prime})^{-1}\times 1)$ (respectively
$(w_2\times 0)\circ(id\times 1)^{-1}$) is the 3-cobordism
$(\Sigma_{g}\times I,f_2\circ(f_1^{\prime})^{-1}\circ f_1\times
0,f_2^{\prime}\times 1) \cong (\Sigma_{g}\times
I,(f_2^{\prime})^{-1}\circ f_2\circ(f_1^{\prime})^{-1}\circ
f_1\times 0,id\times 1) \cong (\Sigma_{g}\times I, (w_2 \circ
w_1)\times 0, id\times 1)$. In particular, the composition of two
morphisms of category $\mathfrak L$ is again a morphism in the
same category.


\medskip
\subsec{Definition {\rm [\ref{CL}]}.}{\it Denote by ${\cal L}_g$ the subgroup of the Mapping Class Group,
consisting of isotopy classes of elements $w\in Aut(\Sigma_g)$
such that $w_*(L^a) = L^a$ and $w_*(L^b) = L^b$ (over $\qq$ or
over $\zz$, is equivalent by the above), and call it
the {\bf Lagrangian subgroup of the MCG}.
}

\medskip
The TQFT of the LMO invariant induces a representation of 
${\cal L}_g$. This subgroup of $MCG(g)$ is big enough to be 
interesting, it contains the Torelli group. Its image under
the action on homology is the group (not normal as subgroup of $Sp(2g,\zz)$)
of matrices of the form
$\left(\begin{array}{cc}A&0\\ 0&(A^T)^{-1}\end{array}\right)$,
where $A\in GL(g,\zz)$.

\medskip\noindent
{\it Remark.} Let $\lambda$ denote the Casson invariant of homology 3-spheres. By fixing the standard handlebody 
of genus $g$ in $\rr^3\subset S^3$ we fixed a Heegaard homeomorphism that Morita [\ref{morita}] calls $\iota_g$, 
and by taking the filling $\widehat{(\Sigma_g\times I,\varphi,id)}$ we obtain a manifold denoted by Morita $W_{\pfi}$.

%
%

\section{The algebraic-combinatorial category}
\setcounter{equation}{0}
\setcounter{nsubsec}{1}

In [\ref{MO}], an essential part of constructing a TQFT associated
to $Z^{LMO}$ has been completed. Namely, first the Kontsevich integral
was extended to an invariant $Z(G)$ of oriented framed trivalent
graphs $G$ in $S^3$ (see [\ref{MO}, theorem 1.4]). A framed 
graph $G\subset S^3$ is represented as a plane projection (with 
implicit blackboard framing), then decomposed into elementary
pseudo-quasi-tangles, and $Z$ is defined for each piece (see
[\ref{MO}, figure 2] for the exact definition of $Z$). 
It is easy to observe that in order to verify the independence 
of $Z$ of the decomposition into pseudo-quasi-tangles and the invariance 
under extended (generalized) Reidemeister moves for trivalent 
graphs, one is forced to introduce relations that ``move'' (in the
sense of Proposition 2.2) a box-diagram over a trivalent vertex 
to a box-diagram. This implies the branching relations (figure 
\ref{figbranch} here, figure 1 in [\ref{MO}]). These relations are 
necessary to impose regardless of the definition of $Z$ for the 
neighbourhood of a trivalent vertex.

From the extended $Z$ Murakami and Ohtsuki [\ref{MO}] derived 
an invariant of oriented 3-manifolds with boundary, along the same 
lines the $Z^{LMO}$ is constructed [\ref{LMO}] from the Kontsevich 
integral of framed links. But in [\ref{MO}] a twisted gluing is used
for composing cobordisms, and the resulting TQFT has complicated  anomaly.

\begin{figure}[p]
\centerline{\psfig{file=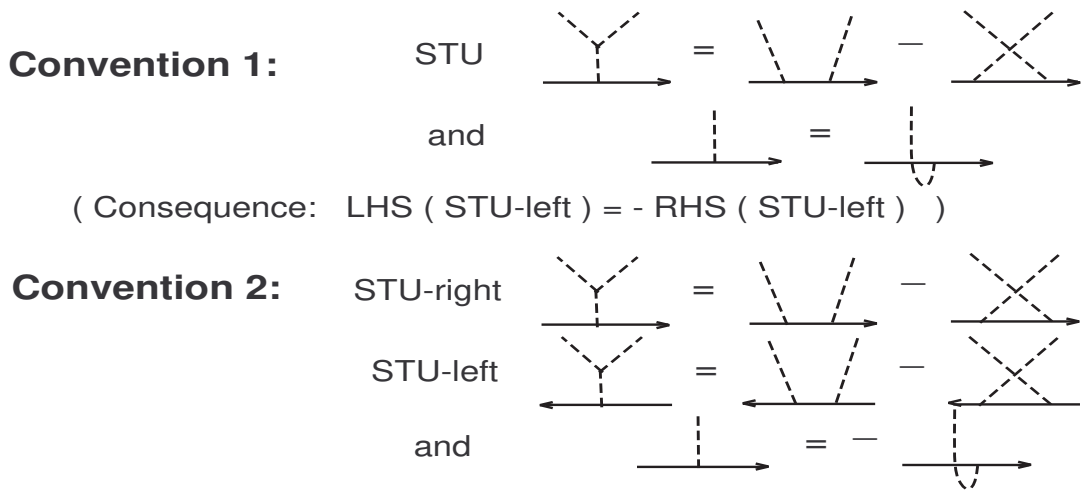,width=15cm,height=7cm,angle=0}}
\caption{\sl The two conventions for chord diagrams}
\label{figconv}
\end{figure}

\medskip
\subsec{The modules of chord diagrams.}Let 
$\Gamma$ be a graph, we will be mainly interested
in the cases $\Gamma=$ a 1-manifold and $\Gamma=$ a chain graph. Let
$\aa(\Gamma)$ be the formal series completion with respect to the
degree of the $\qq$-vector space freely generated by the set of
homeomorphism classes of chord diagrams with support $\Gamma$,
without self-loops and univalent vertices, modulo AS, IHX, STU and
branching relations (which are homogeneous with respect to the
degree).\footnote{There are essentially two conventions in
defining STU and AS relations, and drawing certain elements of
$\aa(\Gamma)$, as shown in figure \ref{figconv}. Note that in
convention 1, which is the one that we use (as well as
[\ref{LMO},\ref{MO}]), AS relations refer only to internal
vertices, and no cyclic order of edges adjacent to external
trivalent vertices is defined. In this convention, as a
consequence, the LHS of STU-left is equal to {\it minus} the RHS
of STU-left. Using the second convention, the definition of a
chord diagram has to be changed as to account for the cyclic order
of edges adjacent to external trivalent vertices. The two
$\aa(\Gamma)$, from the two conventions, are canonically
isomorphic; in fact only the meaning of some diagrams as elements
of $\aa(\Gamma)$ is changed by adding a $-$ sign.}

\begin{figure}[p]
\centerline{\psfig{file=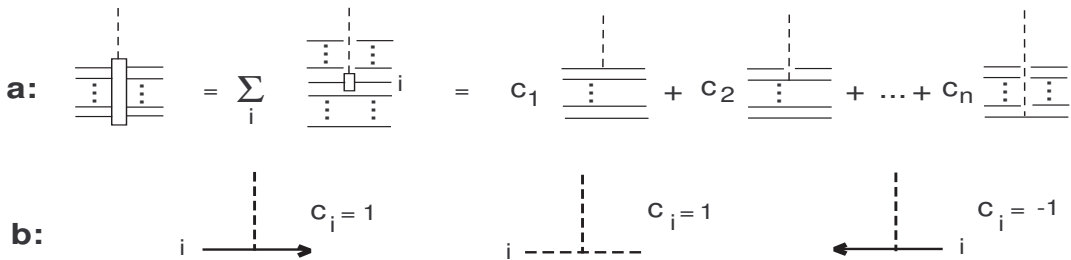,width=15cm,height=4cm,angle=0}}
\caption{\sl The box-diagram} \label{figbox}
\end{figure}

\begin{figure}[p]
\centerline{\psfig{file=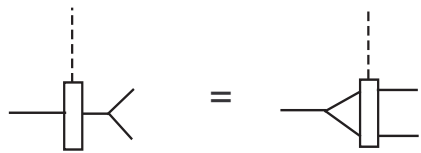,width=6cm,height=25mm,angle=0}}
\caption{\sl The 8 branching relations (all but the vertical edge
are bold): one for each possible orientations of the 3 bold edges}
\label{figbranch}
\end{figure}

We will use the following {\it box-diagram} notation for the
formal sum of chord diagrams, as shown in figure \ref{figbox}a.
There outside the drawn part the diagrams are identical, the
vertical edge is dashed, the horizontal edges are arbitrary. If
the horizontal edge $i$ is dashed, then $c_i=1$, if it is bold,
then $c_i$ is as shown in figure \ref{figbox}b.\footnote{For
convention 2 all coefficients $c_i=1$. Then the box-diagrams in
the two conventions correspond precisely one to the other via the
canonical isomorphism between the conventions.} The branching
relations, introduced in [\ref{MO}, figure 1] are shown in figure
\ref{figbranch} using this box-notation.

Similarly, let $\aa(\emptyset)$ denote the formal series
completion with respect to the degree of the $\qq$-module freely
generated by the set of homeomorphism classes of open chord
diagrams without self-loops and univalent vertices, modulo AS and
IHX relations. For a chord diagram $D$, denote $[D]$ the
corresponding element of $\aa(\Gamma)$. $\aa(\Gamma)$ and
$\aa(\emptyset)$ are co-algebras with respect to the decomposition 
of the dashed part of a diagram in connected components (the elements
represented by diagrams that have non-empty connected dashed part are
defined to be primitive)\footnote{One can check (e.g. by induction 
on the number of internal vertices of chord diagrams) that this 
comultiplication is well-defined (remember the presence of STU 
relations).}. $\aa(\emptyset)$ is an algebra with
respect to disjoint union, and together with (completed) comultiplication
$\Delta$ forms a Hopf algebra. Note that $\aa(\Gamma)$ is 
an $\aa(\emptyset)$-module with respect to the disjoint union.

\begin{figure}[p]
\centerline{\psfig{file=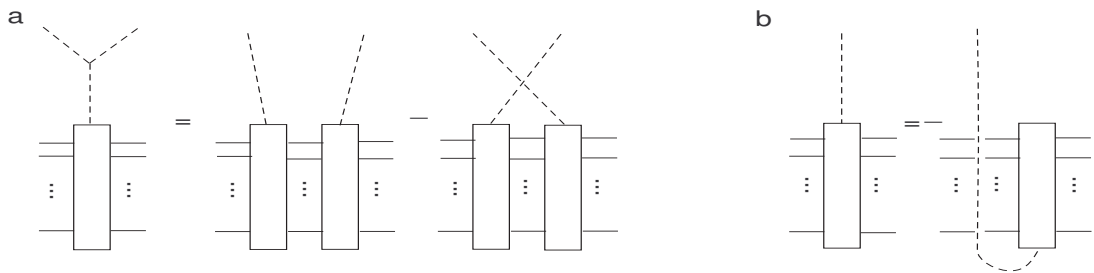,width=15cm,height=4cm,angle=0}}
\caption{\sl {\bf a:} The box-STU relation (each term in the RHS
contains a double sum over the horizontal edges), {\bf b:} The
box-AS relation} \label{figboxrel}
\end{figure}

\begin{figure}[p]
\centerline{\psfig{file=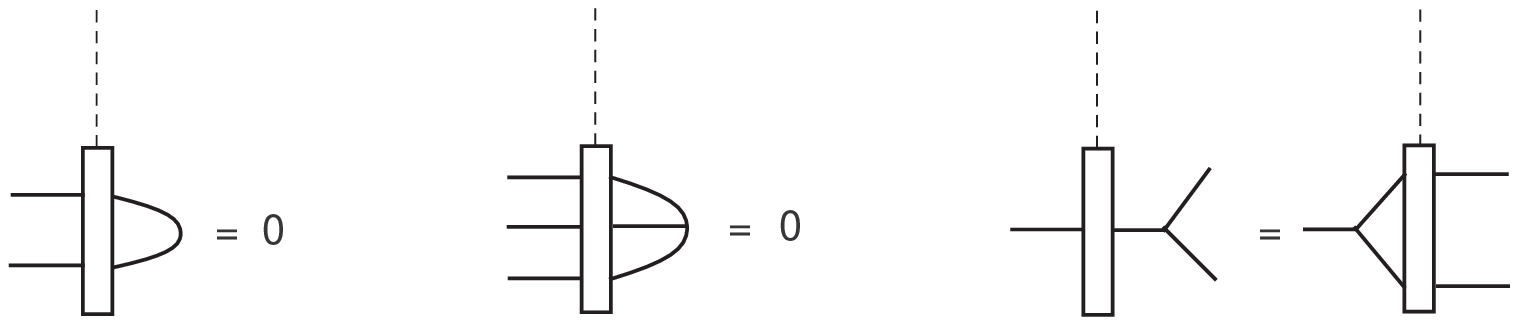,width=16cm,height=35mm,angle=0}}
\caption{\sl Invariance over ``elementary pseudo-tangles'' (the
dashed/bold type of horizontal edges is arbitrary, the vertical
edge is dashed)} \label{fig10}
\end{figure}
\begin{figure}[p]
\centerline{\psfig{file=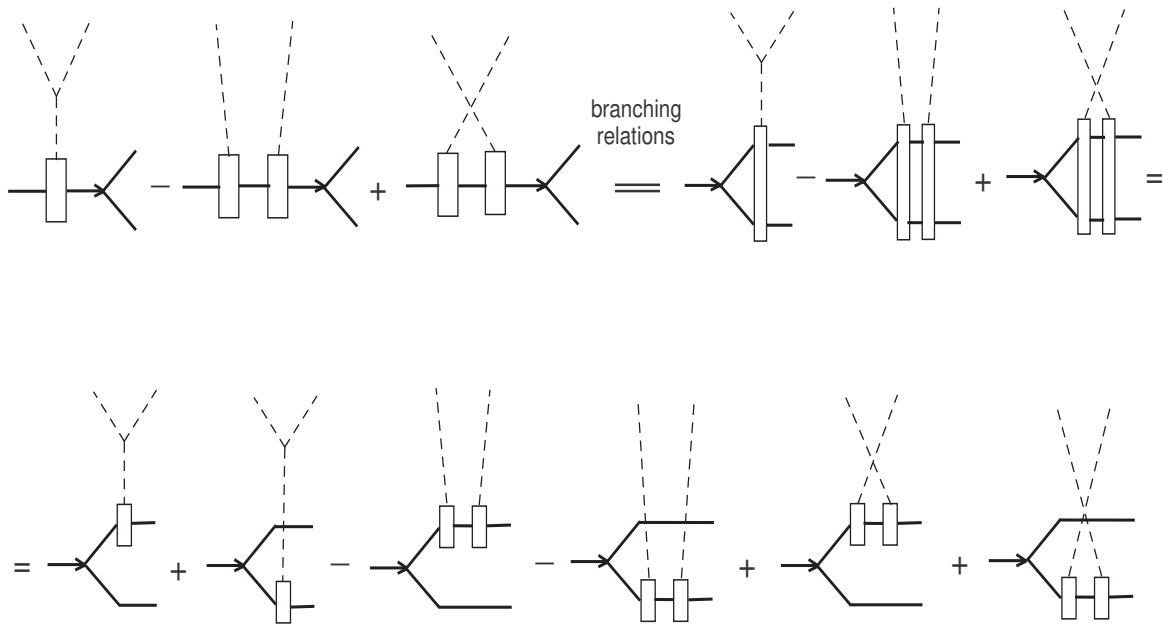,width=16cm,height=9cm,angle=0}}
\caption{\sl ``Moving'' a box-STU relation over a trivalent
vertex} \label{figboxmov}
\end{figure}

\medskip
\subsec{Proposition.}{\it (a) The box-STU and box-AS relations,
schematically shown in figure \ref{figboxrel} hold in $\aa(\Gamma)$.

(b) The three relations in figure \ref{fig10} hold in
$\aa(\Gamma)$.

(c) The box-STU and box-AS relations can be ``moved'' over any
trivalent vertex of $\Gamma$, using only branching relations (see
figure \ref{figboxmov} for an example).\footnote{One can
reformulate this statement: Every IHX (respectively AS) relation
on-the-left-of-the-trivalent-vertex is a consequence of branching
relations and IHX (respectively AS) relations on-the-right-of-the-trivalent-vertex.}
}

\dem{(a) Let $x_i$ denote the horizontal edges. Let $[D^Y]$,
$[D^{II}]$, $[D^X]$ denote the three terms of the box-STU
relation. Note that the brackets are also part of the notation,
$D^Y$ means the box-diagram, which is not a chord diagram. Let
$[D_{x_i}^{Y}]$ denote the element of $\aa(\Gamma)$ corresponding
to the chord diagram obtained from $D^Y$ by replacing the box with
a prolongation of the vertical edge until the the edge $x_i$. With
similar notations $[D^{II}_{x_ix_j}]$ and $[D^{X}_{x_kx_l}]$, note
that for $i\not=j$, $[D^{II}_{x_ix_j}]=[D^{X}_{x_jx_i}]$. Hence:

$$RHS
= \sum_{x_i}\sum_{x_j}c_ic_j[D^{II}_{x_ix_j}]
-\sum_{x_j}\sum_{x_i}c_jc_i[D^{X}_{x_jx_i}] =
\sum\limits_{i=j}c_i^2[D^{II}_{x_ix_j}] +
\sum\limits_{i\not=j}c_ic_j[D^{II}_{x_ix_j}] -
\sum\limits_{i=j}c_i^2[D^{X}_{x_jx_i}] -$$
$$-\sum\limits_{i\not=j}c_jc_i[D^{X}_{x_jx_i}]
= \sum\limits_{i=j}c_i^2[D^{II}_{x_ix_j}] -
\sum\limits_{i=j}c_i^2[D^{X}_{x_jx_i}] = \sum\limits_i
([D^{II}_{x_ix_i}] - [D^{X}_{x_ix_i}]) = \sum\limits_i
c_i[D_{x_i}^{Y}] = LHS$$

\noindent where in the equality before the last we have used an
IHX, STU, or convention-1 form of STU-left for each $x_i$. The
proof of the box-AS relation is elementary, using AS relations and
the definition of coefficients $c_i$.

(b) Consider all dashed/bold possibilities for the edges. The
relations then follow from the AS, IHX, STU and branching
relations.

(c) Every box-diagram is a sum of box-diagrams with small boxes.
For the later follow the calculation shown in figure \ref{figboxmov},
for the box-STU case. The box-AS case is obvious.
}

\medskip
Note that this proposition for the case of $\Gamma$ being a
1-manifold is part of [\ref{vogel}, proposition 1.4].

\begin{figure}[p]
\centerline{\psfig{file=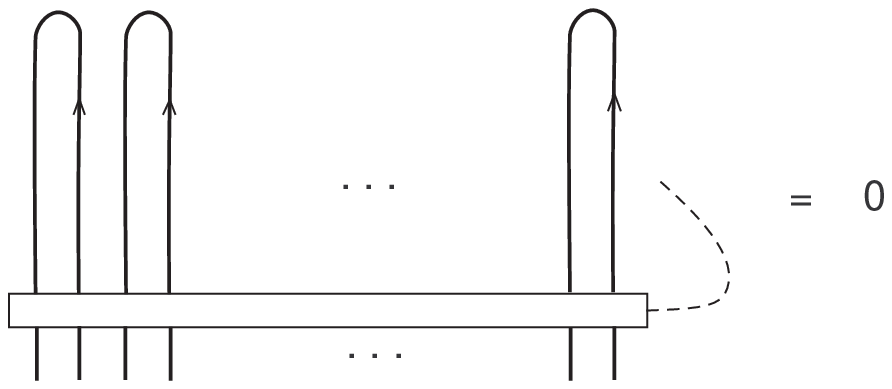,width=12cm,height=4cm,angle=0}}
\caption{\sl Relation 9: the ambiguity of ``moving'' a dashed end
off the horizontal line.} \label{figrel9}
\end{figure}


The ``formal series completion'' (i.e. the topology is given by\footnote{The 
reason we choose $\frac{1}{2^n}$ instead of $\frac{1}{n}$ is Lemma 3.5. See the remark after it.}
{\it distance}$(p,q)\leq\frac{1}{2^n}\Leftrightarrow p-q$ has no terms 
of degree $<n$) is algebraically nothing else but
the direct product over $i\in\nn$ of the vector spaces generated by
diagrams of a fixed order $i$. AS, IXH, STU and branching
relations are homogeneous with respect to the degree. For every
$i$, the degree $i$ part $\aa^i$ is defined as ${\cal D}^i/{\mathfrak
R}^i$, where ${\cal D}^i$ is the $\qq$-module freely generated by
the chord diagrams of degree $i$ (without factoring through
relations), and ${\mathfrak R}^i$ is the $\qq$-module freely
generated by the relations involving only diagrams of order $i$.
By the universal property of the direct product $\aa =
\prod_{i\in\nn}\aa^i \cong \prod_{i\in\nn}{\cal D}^i /
\prod_{i\in\nn}{\mathfrak R}^i = {\cal D}/{\mathfrak R}$, i.e.
factoring and taking completion commute. We will not use anywhere
below the next proposition that $\aa \cong {\cal D}/{\mathfrak
R}$, our object is always $\aa$.

\medskip
\subsec{Proposition.}{\it Denote $[g]\stackrel{\rm not}{=}\{1,\dots,g\}$. 
Let $\phi:\arrowg\rightarrow\oGraph_{[g]}=\Gamma^g$ be the embedding
of $\arrowg$ onto the upper half-circles of $\oGraph$, sending the
arrow labeled $i$ to the $i^{\rm th}$ upper half-circle of
$\Gamma^g$, preserving orientation. Then it extends to an
isomorphism of $\qq$-vector spaces
$\phi_{\ast}:\aa(\arrowg)\rightarrow\aa(\Gamma^g)$.}

\dem{Fix an arbitrary degree $i$ of chord diagrams. Then $\phi$ induces 
a homomorphism of vector spaces
$\phi_{\ast}:\dd^i(\arrowg)\rightarrow\dd^i(\Gamma^g)$, under which
$\fr^i(\arrowg)$ is sent exactly to the set of AS, IHX and STU
relations in $\Gamma^g$ that involve only diagrams with support in
$\phi(\arrowg)$. For simplicity of notation, let us denote
$\phi_{\ast}\dd^i(\arrowg)$ by $\dd^i(\arrowg)$, and
$\phi_{\ast}\fr^i(\arrowg)$ by $\fr^i(\arrowg)$.

Replace each external trivalent vertex in $\Gamma^g-\phi(\arrowg)$
of a chord diagram by a small box (and add a sign to it, the 
coefficient $c_i$), then ``move'', using the branching
relations, one by one all boxes off $\Gamma^g-\phi(\arrowg)$.
This assigns to an arbitrary chord diagram with support in
$\Gamma^g$ a diagram with boxes (with a $\pm$ sign) with support in
$\phi(\arrowg)$. It depends of the choice of the sequence of trivalent 
vertices over which boxes are ``moved'' in $\Gamma^g$.
Observe, however, that different such choices result
in diagrams with boxes, representing elements of $\dd^i(\arrowg)$ that
differ one from the other by a sum (with coefficients $\pm 1$) of
relations depicted in figure \ref{figrel9}. Let us call them {\it
Relations 9} as reference to figure 9 in [\ref{MO}]. By linearity,
this defines a homomorphism of $\qq$-vector spaces
$\alpha:\dd^i(\Gamma^g)\rightarrow\dd^i(\arrowg)/R9$, which when
restricted to $\dd^i(\arrowg)\rightarrow\dd^i(\arrowg)/R9$ is the
canonical quotient map. Here $R9$ is the $\qq$-vector subspace of
$\dd^i(\arrowg)$ generated by the set of Relations 9.

Proposition 2.2(b) implies that Relations 9 are true in
$\aa^i(\arrowg)$, i.e. $R9\subset\fr^i(\arrowg)$. Let
$\beta:\dd^i(\arrowg)/R9\rightarrow\dd^i(\arrowg)/\fr^i(\arrowg)$ be the
canonical projection. Let us observe that for every branching
relation $R$, $\alpha(R)=0$. Therefore $(\beta\circ\alpha)(R)=0$,
so if we denote by ${\mathfrak B}^i$ the $\qq$-vector subspace of
$\dd^i(\Gamma^g)$ generated by the set of branching relations, then
${\mathfrak B}^i\cap\dd^i(\arrowg)\subset\fr^i(\arrowg)$.

On the other hand, any IHX, AS, STU or branching relation on
$\Gamma^g$ is, by Proposition 2.2, a sum of IHX, AS and STU
relations on $\phi(\arrowg)$, plus a sum of branching relations.
Indeed, an IHX or AS relation refers only to a neighbourhood
outside $\Gamma^g-\phi(\arrowg)$, hence the ``moving'' procedure
can be applied simultaneously to all terms of the relation; while
a STU relation is, up to sign, a box-STU relation, therefore using
Proposition 2.2(c) can be ``moved'' to a box-STU relation with
support in $\Gamma^g-\phi(\arrowg)$, the later being a consequence
of $\fr^i(\arrowg)$ by Proposition 2.2(a). The difference between
the start and the end of each step of a ``moving'' procedure is,
of cause, an element of ${\mathfrak B}^i$. Hence
$\fr^i(\Gamma^g)=\fr^i(\arrowg)+{\mathfrak B}^i$.

The two established relations imply
$\fr^i(\Gamma^g)\cap\dd^i(\arrowg)\subset\fr^i(\arrowg)$. Since the
opposite inclusion is obvious,
$\fr^i(\Gamma^g)\cap\dd^i(\arrowg)=\fr^i(\arrowg)$. Then, by the second
isomorphism theorem for vector spaces,
$\dd^i(\arrowg)/\fr^i(\arrowg)\cong\dd^i(\Gamma^g)/\fr^i(\Gamma^g)$.
Composing with $\phi_{\ast}$ from the first paragraph, we obtain
$\phi_{\ast}:\aa^i(\arrowg)\rightarrow\aa^i(\Gamma^g)$, for every $i\geq 0$.
Moreover the induced $\phi_{\ast}:\aa(\arrowg)\rightarrow\aa(\Gamma^g)$
preserves the topology. 
}

\medskip\noindent
{\it Remark.} This Proposition still holds if $\Gamma$ has two or more connected components,
but we can "eliminate" the horizontal line of only {\sc one} component.
If we "eliminate" more than one horizontal line, the corresponding $\phi_{\ast}$ is still
well-defined and surjective.


\medskip
\subsec{The algebra structure of $\aa(\uparrow_{[g]})$ and the set ${\mathfrak C}_{\emptyset}$.}Let 
$\aa_c(\arrowg)$ be the $\qq$-vector subspace of $\aa(\arrowg)$
generated by formal series of diagrams on $\arrowg$ with 
no components of the dashed graph disconnected from the support. 
Viewing each chord diagram as a union of the connected components of
the dashed graph that do not meet the support with the part that meets the support,
we get $\aa(\arrowg)=\aa(\emptyset)\otimes_{\qq}\aa_c(\arrowg)$.
Let $\a(\arrowg)$ be the $\qq$-vector subspace of $\aa_c(\arrowg)$
generated by formal series of diagrams on $\arrowg$ with {\it non-empty and
connected dashed graph} (and connected to the support). 
$\a(\arrowg)$ is precisely the set of primitive elements of 
$\aa_c(\arrowg)$. A similar notation $\a(\Gamma)$ for any
abstract graph $\Gamma$ is self-evident.
$\aa_c(\arrowg)$ is an algebra with respect to justaposition of
the bold vertical arrows. Denote this associative, generally (if $g>1$) non-commutative 
operation $\bullet$. In fact $\aa_c(\arrowg)$ is a co-commutative 
Hopf algebra [\ref{vogel}, Proposition 1.5]. The following is apparently
"common knowledge":

\te{Proposition}{ 1) $\a(\arrowg)$ is a Lie algebra over 
$\qq$ with respect to the operation $(x,y)\mapsto x\bullet y-y\bullet x$.

2) Let $\widehat{I}$ be the topological ideal of $\aa_c(\arrowg)$ generated by 
$\a(\arrowg)$. Then $\exp:\widehat{I}\rightarrow 1+\widehat{I}$ 
and $\log:1+\widehat{I}\rightarrow\widehat{I}$, 
defined by $\exp(x)=\sum_{n=0}^{\infty}\frac{x^n}{n!}$ and 
$\log(1+x)=\sum_{n=1}^{\infty}(-1)^{n+1}\frac{x^n}{n}$, 
where the product is the operation $\bullet$, 
satisfy $\exp\circ\log=id_{1+\widehat{I}}$ and $\log\circ\exp=id_{\widehat{I}}$. 
In particular, $\exp$ and $\log$ are bijections.

3) $\exp$ is a bijection from $\a(\arrowg)\subset\widehat{I}$ to
the set of group-like elements in $1+\widehat{I}$.

4) If $\alpha, \beta\in\a(\arrowg)$, 
then $\exp(\alpha)\bullet\exp(\beta)=\exp(\gamma)$
for some $\gamma\in\a(\arrowg)$. Moreover, $\gamma$ is given 
by the Campbell-Hausdorff formula.

5) $\widehat{I}$ coincides with the set of formal series of chord diagrams of 
degree $\geq 1$.}

\dem{ 1) The statement is sufficient to prove for $x,y=$ diagrams 
with connected dashed graph. Using STU relations, as shown in 
Figure 10, we can interchange two consecutive external vertices, 
one from $x$, the other from $y$, on any bold arrow, up to $\pm$ a 
diagram with connected dashed graph. Therefore iteratively we can 
interchange all external vertices of $x$ with all external vertices 
of $y$, obtaining $x\bullet y-y\bullet x=$ a sum of diagrams with 
connected dashed graph.

2), 3) and 4) are classical statements. The proofs in [\ref{serre}, 
Theorem 7.2, Corollary 7.3 and Theorem 7.4] apply m\^ ot-a-m\^ ot.
For 3) and 4) note that 
if $\gamma$ is primitive, then $\gamma\in\a(\arrowg)$.

5) Since the set of formal series of chord diagrams of 
degree $\geq 1$ is an ideal containing $\a(\arrowg)$, and is closed 
topologically, $\widehat{I}$ certainly belongs to it. Conversely, 
pick an arbitrary connected component $y^{\prime}$ of the dashed 
graph of a chord diagram. Observe that using the "`trick"' in Figure 10, 
up to $\pm$ a sum of diagrams with the number of connected components of the 
dashed graph less by $1$, $y^{\prime}$ can be assumed to have all 
external vertices below all the other external vertices of the 
diagram. Hence an induction on the number of connected components 
of the dashed graph shows that any chord diagram of degree $\geq 1$ 
is a sum of terms of type $\pm z_1\bullet z_2\bullet\dots\bullet z_k$, 
$k\geq 1$, with $z_i$ a diagram in $\a(\arrowg)$. We conclude that 
the set of finite sums of chord diagrams of degree $\geq 1$ 
is contained in $\widehat{I}$. Hence so is its completion.
}

\begin{figure}[thb]
\centerline{\psfig{file=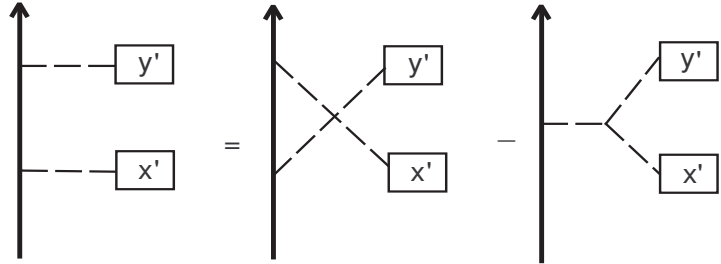,width=10cm,height=4cm,angle=0}}
\caption{\sl Two consecutive external vertices from connected components $x^{\prime}$ and $y^{\prime}$ can be 
interchanged up to $\pm$ a diagram with dashed graph having one component less.} 
\end{figure}

\begin{figure}[thb]
\centerline{\psfig{file=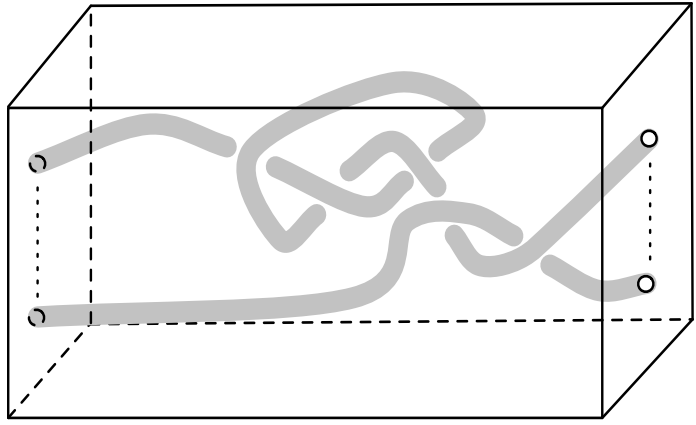,width=10cm,height=7cm,angle=0}}
\caption{\sl $(L,G)$ less a tubular neighbourhood of the horizontal line of $G$.} \label{multiplicationBox}
\end{figure}

This Proposition holds if we replace $\aa_c(\arrowg)$ by $\aa(\arrowg)$ and $\qq$ by $\aa(\emptyset)$.
It suggests us to consider an operation on
${\mathfrak C}_{\emptyset}$, the set
of connected 3-cobordisms with empty bottom
and connected top, to correspond to the multiplication in
$\aa(\arrowg)$. Let $(M_1, \emptyset, f_1),
(M_2, \emptyset, f_2)\in {\mathfrak C}_{\emptyset}$,
and let $(L_1,G_1)\subset S^3$, $(L_2,G_2)\subset S^3$ 
such that $\kappa(L_1,G_1)=(M_1,\emptyset,f_1)$,
$\kappa(L_2,G_2)=(M_2,\emptyset,f_2)$.\footnote{The notion 
of a pair $(L,G)$ is defined in the same way as that of a triplet.}
Remove a tubular neighbourhood of the horizontal line of
$G_1$ ($=$ a ball in $S^3$), and similarly for $G_2$.
Up to isotopy we can assume each looks 
as the box in Figure 11. Glue the two boxes from left to right,
and fill back in the standard way a horizontal line. Denote the result
by $(L_1\cup L_2,G_1\bullet G_2)$, and define:

\begin{equation}\label{multiplication}
(M_1, \emptyset, f_1)\bullet(M_2, \emptyset, f_2)=\kappa(L_1\cup L_2,G_1\bullet G_2)
\end{equation}

\noindent
Observe that the new 3-cobordism does not depend on the 
choice of pairs $(L_i,G_i)$, such that $\kappa(L_i,G_i)=(M_i,\emptyset,f_i)$,
since $(L_1\cup L_2,G_1\bullet G_2)$ needs to be determined only up to extended 
KI, KII relations, and change of orientation of link components.
In the case of $g=0$, $\bullet$ is the connected sum. Hence this 
operation is an alternative (to composition of cobordisms) way of generalizing connected sum.
Note that $\widehat{(M_1,\emptyset,f_1)\bullet(M_2,\emptyset,f_2)}=
\widehat{(M_1,\emptyset,f_1)}\#\widehat{(M_2,\emptyset,f_2)}$.
In particular the sets ${\mathfrak C}_{\emptyset}\cap\{\zz{\rm -cobordism}\}$
and  ${\mathfrak C}_{\emptyset}\cap\{\qq{\rm -cobordism}\}$ are closed under $\bullet$.


\medskip
\subsec{The LMO invariant for closed manifolds and extending the maps $\iota_n$.}In 
[\ref{LMO}] from the Kontsevich integral an invariant of
oriented framed links $L$ was constructed, which does not change under Kirby-1,2
moves and change of orientation of components of $L$. 
We recall it here, together with the maps $\tilde\iota_n$, necessary to extend 
it to an invariant of unions of embedded framed {\it chain} graphs in $S^3$.

Let $\stackrel{\circ}{\aa(\emptyset)}$ be the formal series completion of the
$\qq$-vector space generated by the homeomorphism classes of open chord diagrams without 
univalent vertices (but allowing dashed self-loops - these are set of degree $0$)
modulo AS and IHX relations. Let ${\cal B}(X)$ be the formal power
series completion of the $\qq$-vector space generated by the
homeomorphism classes of open chord diagrams without dashed self-loops,
with the univalent vertices colored by elements of $X$, modulo AS
and IHX relations.

Denote $[m]\stackrel{\rm not}{=}\{1,\dots,m\}$, and let 
$\Gamma=\sqcup_{m}S^1$, where each component is colored by a 
different element of $[m]$. Let $\widetilde{{\cal B}([m])}$ be
the subspace of ${\cal B}([m]\cup\{\ast\})$, generated by the
diagrams with one $\ast$-colored vertex, and
$f_i:\widetilde{{\cal B}([m])} \rightarrow {\cal B}([m])$,
$f_i:=average\; of\; the\; diagrams\; obtained\;$ $by\;
attaching\; the\; \ast-vertex\; near\; all\; i-vertices$. We can
define a map $\varphi:{\cal C}(m) := \frac{{\cal B}([m])}{<im\:
f_i | \forall i>} \rightarrow{\cal A}(\Gamma)$, $\varphi :=
average$ {\it of the diagrams obtained by attaching $i$-coloured
vertices to the $i^{th}$ copy of $S^1$ in $\Gamma$, $\forall i$}.
(One checks that the definition on diagrams extends over relations
to a map between formal series completions.) This map is in fact
an isomorphism of $\qq$-modules (and co-algebras). For details, 
please consult [\ref{LMO},\ref{vogel}]. (If $\Gamma=S^1$, 
${\cal C}(1)$ and ${\cal A}(S^1)$ are algebras, but $\pfi$ is not 
an algebra homomorphism.)

For every $n\geq0$ define a map $\kappa_n:{\cal C}(m)\rightarrow\ao$;
$\kappa_n(K)=0$, if $\exists i$ such that the number of
$i$-colored vertices is not $2n$, $\kappa_n(K)=sum\; of\; all\;
ways\; of\; attaching\; i-coloured\; vertices\; in\;$ $pairs$,
$\forall i$, otherwise. Let ${\cal O}_n$ be the ideal of $\ao$
generated by $\dashedO + 2n$. ($\ao$ is an algebra with respect to
disjoint union.) It can be shown that as modules
(and even as algebras)
$\ao / {\cal O}_n \cong {\cal A}(\emptyset)$. Now, let
$  \iota_n=q_n\circ\kappa_n\circ\pfi^{-1} : \aa(\Gamma) \rightarrow
{\cal C}(m) \rightarrow \ao \rightarrow \ao / \on \cong
\aa(\emptyset)$, where $q_n$ is the quotient map.

Let $\kappa_n^{\ast}:{\cal B}(X\sqcup\{\ast\}) \rightarrow {\cal
B}(X)$ be defined as $\kappa_n$, but only involving
$\ast$-colored vertices (see [\ref{LMO},\ref{vogel}] for details).
Let $P_n = {\rm Im} (\kappa_n^{\ast})$. The map $\kappa_n^{\ast}$
passes to the quotient from the definition on ${\cal C}(m)$, and
hence we get a submodule $P_n$ of ${\cal C}(m)$. The relations 
$P_n$ also commute with $\pfi$. Define the quotient map
$j_n : \ao / \on \rightarrow \ao / \on , P_{n+1}$ of graded modules.
It is isomorphism in degree $\leq n$, and is the main
ingredient in showing that $ \iota_n(\cZ(L))_{\leq n}$ is invariant
under the second Kirby move [\ref{LMO}, \ref{vogel}].

This construction can be extended for $\Gamma=\Gamma^{g_1}\sqcup
(\sqcup_mS^1)\sqcup\Gamma^{g_2}$, disjoint union of two chain
graphs and $m$ copies of $S^1$, i.e.  
$ \iota_n:\aa(\sqcup_mS^1)\rightarrow\ao/\on\cong\aa(\emptyset)$
can be extended (meaning that for $g_1=g_2=0$, $\tilde  \iota_n$ acts 
exactly as $ \iota_n$) to a map:

\begin{equation}\label{in}
\tilde  \iota_n=\tilde q_n\circ\tilde{\kappa_n}\circ\tilde{\pfi}^{-1}:
\aa(\Gamma^{g_1}\sqcup (\sqcup_mS^1)\sqcup\Gamma^{g_2})
\rightarrow \stackrel{\circ}{{\cal
A}(\Gamma^{g_1}\sqcup\Gamma^{g_2})}/\on \cong{\cal
A}(\Gamma^{g_1}\sqcup\Gamma^{g_2})
\end{equation}

\noindent
where the corresponding homomorphism $\tilde{\varphi}^{-1}$ 
refers only to {\it all present circle components} of $\Gamma$. 
Here, to define the preimage of 
$\tilde{\pfi}:{\cal C}(\Gamma^{g_1},[m],\Gamma^{g_2})
\rightarrow\aa(\Gamma^{g_1},\sqcup_mS^1,\Gamma^{g_2})$ we consider
absolutely analogous chord diagrams with support the disjoint 
union of two chain graphs $\Gamma^{g_1}$, $\Gamma^{g_2}$, and 
points indexed by elements of $[m]$ (it is convenient {\sc not} 
to call these points vertices), $\tilde\kappa_n$ is extended in 
the same manner, and $\tilde q_n$ is just the quotient map. $\tilde{\pfi}$ 
is an isomorphism with the proof of Section 2 of [\ref{LMO}]. Moreover, 
the similarly constructed map $\tilde j_n$ is an isomorphism 
in degree $\leq n$. Namely, and this is exactly the statement of 
Lemma 2.3 of [\ref{MO}], if
$\Gamma=\Gamma^{g_1}\sqcup(\sqcup_mS^1)\sqcup\Gamma^{g_2}$, 
then $\tilde{j_n}:({\cal A}(\Gamma^{g_1}\sqcup\Gamma^{g_2})
\cong 
\stackrel{\circ}{{\cal A}
(\Gamma^{g_1}\sqcup\Gamma^{g_2})} / \on)_{\leq n}
\stackrel{\cong}{\longrightarrow} 
(\stackrel{\circ}{{\cal A}(\Gamma^{g_1}\sqcup\Gamma^{g_2})} / 
\on, P_{n+1})_{\leq n}$. To check this fact it is enough to
follow the proof of lemma 3.3 in [\ref{LMO}] or proposition
4.4 in [\ref{vogel}].

Let $Z(L)$ be the usual Kontsevich integral of the (oriented)
framed link $L$, $\nu$ be $Z$ of the zero-framed unknot. Denote
$\cZ(L):=Z(L)\otimes\nu^{|L|}$, meaning we take the ``connected
sum'' of $Z(L)$ on each its component with $\nu$. Like $Z(L)$,
$\cZ(L)$ is also group-like of the form $1+(terms\; of\;
degree\geq 1)$. \footnote{It can be shown by induction that then
for $|L|=1$ the formal graded series $log(element)$ is a primitive 
element of $\aa(\sqcup S^1)$, and has no part of degree $0$, hence it
is a formal power series of chord diagrams with connected dashed
part. More precisely, a statement similar to Proposition 2.4 holds.} 
Let $\sigma_{\pm}$ be the number of positive, resp. negative
eigenvalues of the linking matrix of $L$. Denote $\oo^{+1}$, resp.
$\oo^{-1}$ the unknot with $+1$, resp. $-1$ framing, and $S_L^3$
the $3$-manifold obtained by surgery on the framed link $L$ in
$S^3$. Recall the definition of the LMO invariant for oriented
closed 3-manifolds $M\equiv S^3_L$:

\begin{equation}\label{omega}
\Omega_n(S_L^3) := \left( \frac{ \iota_n(\cZ(L))}{ \iota_n(\cZ(\oo^{+1}))
^{\sigma_+} \cdot  \iota_n(\cZ(\oo^{-1}))^{\sigma_-}}\right)_{\leq n}
\end{equation}

\noindent and:

\vspace{-4mm}
$$Z^{lmo}(M) := \sum\limits_{n \geq 0}\Omega_n(M)_n$$

\noindent 
and for $\qq$-homology spheres also:

$$Z^{LMO}(M) :=  \sum\limits_{n \geq 0}d(M)^{-n}\Omega_n(M)_n$$

\noindent where $d(M)=|det(lk(L))|$, which is $0$ if 
$H_1(S^3_L, \qq) \neq 0$ and $|H_1(M,\zz)|$ otherwise. We use the
convention $|det(lk(\emptyset))|=1$. Then we have
$\Omega_{n+1}(S^3_L)_{\leq n} = d(M)\cdot \Omega_n(S^3_L)$, hence
we can write $Z^{LMO}(M)_{\leq n} = d(M)^{-n}\Omega_n(M)$.
More precisely, the following holds [\ref{vogel}, Proposition 4.5]:

\begin{equation}\label{eqvogel}
[\iota_{n+1}\cZ(L)]_{\leq n} = (-1)^{|L|}det(lk(L))[\iota_n\cZ(L)]_{\leq n}
\end{equation}

\noindent
and therefore we can define:

\begin{eqnarray}\label{cplus} 
c_+ & = & \lim\limits_{n\rightarrow\infty}(-1)^n[\iota_n\cZ(\oo^{+1})]_{\leq n}\\
c_- & = & \lim\limits_{n\rightarrow\infty}[\iota_n\cZ(\oo^{-1})]_{\leq n}\label{cminus}
\end{eqnarray}

\noindent 
These elements of $\aa(\emptyset)$ are canonical constants in the
theory of LMO invariant. (\ref{eqvogel}) implies

\begin{equation}\label{LMOtrunc}
Z^{LMO}_{\leq N}(M)=\frac{(-1)^{N\sigma_+}}{d(M)^{N}}\cdot
\left(\iota_{N}\left(\frac{\cZ(L)}
{c_+^{\sigma_+}c_-^{\sigma_-}}\right)\right)^{[\leq N]}
\end{equation}

\noindent
where the notation $[\leq N]$ means the minimal internal degree in the sense 
of 2.8.

We restrict to the case of $\qq$-homology spheres. Since 
$\cZ(L\sqcup L^{\prime})=\cZ(L)\cZ(L^{\prime})$, we have
$\Omega_n(M\#M^{\prime}) = \Omega_n(M)\Omega_n(M^{\prime})$,
therefore:

\begin{equation}\label{Zconnsum}
Z^{LMO}(M\#M^{\prime}) = Z^{LMO}(M) Z^{LMO}(M^{\prime})
\end{equation}

\noindent
And since $\overline{S^3_L}=S^3_{\overline{L}}$, where $\overline{L}$
is the mirror image of $L$, we have [\ref{LMO}] :

\begin{equation}\label{Zchangeor}
Z^{LMO}(-M)_n = (-1)^{n(b_1+1)} Z^{LMO}(M)_n
\end{equation}

\noindent
where $b_1=rank\: H_1(M,\zz)$.
Let $\omega^{LMO}(M)=log(Z^{LMO}(M))$. Then the two formulas above
can be re-written:

\begin{eqnarray}
\omega^{LMO}(M\#M^{\prime})_n & = & \omega^{LMO}(M)_n+\omega^{LMO}(M^{\prime})_n\label{omega.connsum}\nonumber\\
\omega^{LMO}(-M) & = & \sum\limits_{n=1}^{\infty}(-1)^{n(b_1+1)}\omega^{LMO}(M)_n\label{omega.changeor}\nonumber\\
& = & \left\{\begin{array}{ccc}\overline{\omega^{LMO}(M)}&,&if\:b_1=even\\ \omega^{LMO}(M)&,&if\:b_1=odd\end{array}\right.\nonumber
\end{eqnarray}

\noindent
where the conjugation in $\aa(\emptyset)$ is defined as identity
on chord diagrams of even degree, and multiplication by $-1$ on
chord diagrams of odd degree.

Also define as in [\ref{MO}] $\cZ(L\cup G)=Z(L\cup G)\otimes(\nu^{\otimes |L|})$, 
i.e. add $\nu$ to $Z$ of each component of $L$.


\medskip
\subsec{The definition of $Z$ on elementary pseudo-quasi-tangles.}To 
extend $\Omega_n(S^3_L)\in\aa(\emptyset)$ to invariants $\Omega_n(L,G)\in\aa(\Gamma)$, where
$\Gamma$ is $G$ as abstract graph, we will extend now $Z(L)$ to $Z(L\cup G)$. However we shall do
this differently from Murakami and Ohtsuki [\ref{MO}, see Figure 2 there], who use Knizhnik-Zamolodchikov associator.
We will use the even associator.

Let $G$ be an embedded framed graph in $S^3$. Fix a plane projection
such that $G$ is given the blackboard framing. This 
projection of $G$ can be decomposed into elementary tangles\footnote{The words {\em quasi} and {\em pseudo} are left
out for simplicity of language.}: \tangleOne, \tangleTwo, \tangleThree,
\tangleFour, \UP, \LP, \tresaUnu, \tresaDoi and \vertical. We need only to specify the definition of $Z$ on the first 
two, since on the others we know it from the link case.

Let $\Gamma$ be an abstract (disjoint union of) chain graph(s), and $\epsilon_e\Gamma$
be $\Gamma$ with edge $e$ erased. Suppose $\epsilon_e\Gamma$ is also a chain graph.
A similar notation $\epsilon_e G$ for a framed graph $G$ is self-evident.
Define the map $\epsilon_{(e)}:\aa(\Gamma)\rightarrow\aa(\epsilon_e\Gamma)$,
$\epsilon_{(e)}(D)=0$, if $D$ has an external vertex on the removed edge,
and $\epsilon_{(e)}(D)=D$, otherwise. To verify well-defineness of $\epsilon_{(e)}$
it is enough to check its invariance under branching relations of diagrams on $\Gamma$.
There are 3 diagrams involved in a branching relation. Suppose $v$ is a trivalent
vertex of $\Gamma$, and $e_1, e_2, e_3$ the edges adjacent to $v$. Edge $e$ cannot
be repeated twice among $e_1, e_2, e_3$, since then $\epsilon_e\Gamma$ would not be a 
(union of) chain graph(s). Therefore we can assume $e=e_1$, $e\not=e_2$, $e\not=e_3$.
It is easy to check that then one of the three diagrams in the relation is sent to $0$
by $\epsilon_{(e)}$, while the other two are sent to diagrams that form an AS relation
in $\aa(\epsilon_e\Gamma)$.

If $e$ is an edge of $G$, denote by $S_e G$ the graph obtained from G by reversing the
orientation of the edge $e$ (without changing the framing). If $\Gamma$ is the underlying
abstract graph of $G$, denote by $S_e\Gamma$ the underlying abstract graph of $S_e G$.
Let $S_{(e)}:\aa(\Gamma)\rightarrow\aa(S_e\Gamma)$ be the linear map which sends 
every diagram $D$ in $\aa(\Gamma)$ to the diagram obtained from $D$ by reversing the
orientation of $e$, multiplied by $(-1)^m$, where $m$ is the number of vertices of $D$
on the edge $e$.

We define $Z$ for the elementary
tangles \tangleOne and \tangleTwo to satisfy the following two conditions
(compare with [\ref{MO}, Proposition 1.5]):

\medskip
\rm(1)\it \quad $ Z(S_e G) = S_{(e)}Z(G)$, \quad for any embedded framed graph $G$ and edge $e$

\rm(2)\it \quad $ Z(\epsilon_e G) = \epsilon_{(e)}Z(G)$, \quad for any (disjoint
union of) embedded chain graph(s) $G$ and edge $e$, such that $\epsilon_e G$ is still a
(disjoint union of) embedded chain graph(s)\rm

\noindent
Moreover, we seek to define $Z(\tangleTwo)$ to be of the form $\diaOne$
(for all possible 8 orientation). By condition (2) above we must have $a=b=c^{-1}$, 
and hence also $Z(\UP)=\diaTwo$.
But $Z(\UP)=\nu^{1/2}$, therefore we must require $a=b=c^{-1}=\nu^{1/4}$, i.e.
$Z(\tangleTwo)=\diaThree$.
Similarly $Z(\tangleOne)=\diaFour$.
(These formulas are each for the 8 possible orientations.)


\medskip
\subsec{Theorem.}{\it $Z(G)$ is an isotopy invariant of embedded framed chain graphs.}

\dem{In a big part, this is mostly a repetition of the proofs of the statements in
[\ref{MO}, Section 1], hence we only sketch here the details that are not identical.
First, one shows that $Z(G)$ is invariant under isotopies of the plane. If such isotopies 
fix a neighbourhood of each trivalent vertex, the result is known from the link case.
If isotopies move such neighbourhoods ``as a whole'', the result follows using branching 
relations [\ref{MO}, Lemma 1.2]. Finally, it is sufficient to show $Z(\tangleFive)=Z(\tangleOne)$
and $Z(\tangleSix)=Z(\tangleOne)$.

\smallskip
Secondly, one shows that $Z(G)$ is invariant under extended Reidemeister moves. This is
also easily achieved from results known from the link case and the branching relations
[\ref{MO}, Lemma 1.4].

To prove the two remaining relations, note that in [\ref{LeGrenoble}, page 8]
it is proved (using an even associator) that $Z(\tangleSeven)=\diaFive$. Therefore
$Z(\tangleFive)=\diaSix=\diaSeven\stackrel{branching\: relations}{=}\\ =\diaEight=\diaFour=Z(\tangleOne)$.
Similarly $Z(\tangleSix)=Z(\tangleOne)$.
}

\vspace{2cm}
The above properties (1) and (2) we have now for granted (compare to [\ref{MO},
Proposition 1.5]). It then follows directly from their definitions in section 3.1 that $\tau^{\leq N}$ and $\tau$
also enjoy properties (1) and (2).

It seems that statements in this section (using even associator) 
have been known to different people, but a complete proof was missing from the literature.

\medskip
\te{Conjecture A}{If $G$ is a chain graph, then this definition of $Z(G)$ using even associator coincides
with the definition in [\ref{MO}], which uses KZ associator.}

\medskip
We have been able to obtain only partial results toward the proof of this conjecture. The results of this
paper are equally true for any associator for which Theorem 2.7 holds.
Also, note that branching relations must be introduced (in addition to IHX, AS and STU) 
regardless how one defines $Z(G)$.

\medskip\noindent
{\it Remark.} If we use even associator it is easy to see that $Z(\tangleEight)=\phi_{\ast}(\diaNine)$,
where $\phi_{\ast}$ is the isomorphism from Proposition 2.3.


\medskip
\subsec{The composition of chord diagrams.}Let 
$(\Gamma^{g_1},\Gamma^{g_2})$ be an ordered pair of chain 
graphs. Every time we consider such a pair, $\Gamma^{g_1}$ is the 
union of its horizontal line and the upper half-circles, and 
$\Gamma^{g_2}$ is the union of its horizontal line and lower 
half-circles. Denote $\aa(\Gamma^{g_1},\Gamma^{g_2})=
\aa(\Gamma^{g_1}\sqcup\Gamma^{g_2})$, where the order 
$(\Gamma^{g_1},\Gamma^{g_2})$ has been specified. 
As we have remarked after Proposition 2.3, every element 
of $\aa(\Gamma^{g_1},\Gamma^{g_2})$ has a representative, a formal 
series (with rational coefficients) of chord diagrams, 
whose external vertices don't meet one horizontal line.
We will "remove" a horizontal line when needed using this isomorphism.

Besides (total) degree, which is half the number 
of all vertices, a chord diagram has {\it internal degree}, half the
number of internal vertices. For arbitrary 
$\alpha\in\aa(\Gamma^{g_1},\Gamma^{g_2})$, 
choose a formal sum of chord diagrams (an element of 
${\cal D}(\Gamma^{g_1},\Gamma^{g_2})$) representing $\alpha$, 
express each chord diagram $d$ as a product $a\cdot\beta$, where 
$a$ is an open chord diagram (i.e. an element of 
${\cal D}(\emptyset)$), say of internal degree $n$, and all 
connected components of $\beta$  intersect the support 
$\Gamma$ non-empty. Using STU relations, rewrite $\beta$ (in 
$\aa(\Gamma^{g_1},\Gamma^{g_2})$) as a finite formal sum of chord 
diagrams without internal vertices, hence of internal degree zero. 
Group the terms of $\alpha$ with the same internal degree to 
obtain $\alpha^{[n]}\in\aa(\Gamma^{g_1},\Gamma^{g_2})$, $\forall n$. 
(Note that $\alpha^{[n]}$ is in general
a formal series, not just a finite sum.)
Since STU are the only internal degree non-homogeneous relations,
$\alpha^{[n]}$, $\forall n$ are uniquely determined by $\alpha$. 
Call $\alpha=\sum\limits_{n\geq 0}\alpha^{[n]}$ {\it the internal 
degree decomposition} of $\alpha$, and $\alpha^{[n]}$ its {\it internal degree $n$ part}.

For any non-negative integer $N$, define 
$\aa^{\leq N}(\emptyset)=\aa(\emptyset)/D^{[>N]}$, where $D^{[>N]}$ 
is the subspace spanned by diagrams of degree $>N$. 
The sequence of natural $\qq$-linear projection maps  
$\aa^{\leq N+1}(\emptyset)\rightarrow\aa^{\leq N}(\emptyset)$, 
that forget the degree $N+1$ part, i.e.: 

$$\dots\rightarrow\aa^{\leq N+1}(\emptyset)
\rightarrow\aa^{\leq N}(\emptyset)\rightarrow\dots
\rightarrow\aa^{\leq 1}(\emptyset)\rightarrow\aa^{\leq 0}(\emptyset)
=\qq$$

\noindent
has inverse limit the direct product of homogeneous degree parts, i.e.
$\aa(\emptyset)$. (This has already been used in (\ref{cplus}) and (\ref{cminus}).)
Absolutely similarly, with the only observation that 
the degree is the minimal internal degree of diagrams, we define 
$\aa^{\leq N}(\Gamma^{g_1},\Gamma^{g_2})
=\aa(\Gamma^{g_1},\Gamma^{g_2})/D^{[>N]}$ and present 
$\aa(\Gamma^{g_1},\Gamma^{g_2})$ as the inverse limit of a similar
sequence:

$$\dots\rightarrow\aa^{\leq N+1}(\Gamma^{g_1},\Gamma^{g_2})
\rightarrow\aa^{\leq N}(\Gamma^{g_1},\Gamma^{g_2})\rightarrow\dots
\rightarrow\aa^{\leq 1}(\Gamma^{g_1},\Gamma^{g_2})
\rightarrow\aa^{\leq 0}(\Gamma^{g_1},\Gamma^{g_2})
=\aa_c(\Gamma^{g_1},\Gamma^{g_2})$$

\noindent
Of cause, $\aa(\emptyset)$ and $\aa(\Gamma^{g_1},\Gamma^{g_2})$
are also inverse limits of infinite subsequences of the above, e.g.
if we only consider $N$ even.
Note also the natural isomorphism $\aa^{\leq N}(\Gamma^{g_1},\Gamma^{g_2})
\cong\aa^{\leq N}(\emptyset)\otimes_{\qq}\aa_c(\Gamma^{g_1},\Gamma^{g_2})$.

The $\qq$-vector spaces $\aa^{\leq N}(\Gamma)$ have a 
$\aa^{\leq N}(\emptyset)$-module structure given by
looking at the multiplication by $\aa(\emptyset)$ in $\aa(\Gamma)$ in the quotient: 
the product of $\alpha\in\aa^{\leq N}(\emptyset)$ 
with $\beta\in\aa^{\leq N}(\Gamma)$ 
is $(\alpha\beta)^{[\leq N]}\in\aa^{\leq N}(\Gamma)$,
where the multiplication $\alpha\beta$ is induced by the disjoint union 
of chord diagrams. Similarly the multiplication $\bullet$ in $\aa(\arrowg)$ induces 
one on $\aa^{\leq N}(\arrowg)$, making it an algebra.

However for the comultiplication we don't have the trouble of quotienting through 
$^{[>N]}$. Indeed, $\Delta:\aa(\Gamma)\rightarrow\aa(\Gamma)\widehat{\otimes}\aa(\Gamma)$ 
preserves both the degree and the internal degree parts, as it is easy to observe 
from its definition. Hence it induces 
$\Delta^{\leq N}:\aa^{\leq N}(\Gamma)\rightarrow(\aa^{\leq N}(\Gamma)
\widehat{\otimes}\aa^{\leq N}(\Gamma))^{\leq N}\subset\aa^{\leq N}(\Gamma)
\widehat{\otimes}\aa^{\leq N}(\Gamma)$. We can drop the index in $\Delta^{\leq N}$ since 
the quotient $\aa(\Gamma)\rightarrow\aa^{\leq N}(\Gamma)$ admits a canonical 
section $\aa^{\leq N}(\Gamma)\rightarrow\aa(\Gamma)$, which is to view every 
formal series of diagrams of minimal internal degree $\leq N$ as a formal series 
of diagrams, i.e. $\Delta^{\leq N}$ can be equally viewed as the restriction 
of $\Delta$ to a submodule. An element of $\beta\in\aa^{\leq N}(\Gamma)$ will 
be called {\bf primitive} if $\Delta^{\leq N}\beta=\beta\otimes 1+1\otimes\beta$.
(In fact $\beta$ is primitive in  $\aa^{\leq N}(\Gamma)$ iff it is primitive in $\aa(\Gamma)$.)
$\alpha\in\aa^{\leq N}(\Gamma)$ will be called {\bf group-like} if 
$\Delta^{\leq N}\alpha=(\alpha\widehat{\otimes}\alpha)^{[\leq N]}$. 
Moreover, even though $\Delta^{\leq N}$ are not homomorphisms,
Proposition 2.4 holds true if we replace 
$\aa(\arrowg)$, $\a(\arrowg)$ and $\widehat{I}$ by their internal 
degree $\leq N$ truncations $\aa^{\leq N}(\arrowg)$, $\a^{\leq N}(\arrowg)$ 
and $\widehat{I}^{\leq N}$, provided we use the above definitions of 
primitive and group-like. It follows from the fact that 
$\aa^{\leq N}(\arrowg)=\aa^{\leq N}(\emptyset)\otimes_{\qq}\aa_c(\arrowg)$.
This is a simple, yet important observation.
$\aa^{\leq N}(\arrowg)$, and via Proposition 2.3 also $\aa^{\leq N}(\Gamma^g)$, 
have all the properties of non-commutative co-commutative Hopf algebras,
except $\Delta$ being homomorphisms.\footnote{Alternatively,
we can redefine the tensor products $\aa^{\leq N}(\Gamma)\overline{\otimes}\aa^{\leq N}(\Gamma)
:=\aa^{\leq N}(\Gamma)\widehat{\otimes}\aa^{\leq N}(\Gamma)
/\left(\aa^{\leq N}(\Gamma)\widehat{\otimes}\aa^{\leq N}\right)^{>N}$, i.e.
enlarge the notion of {\em graded} Hopf algebras with many interesting examples.
But apparently this approach is disliked by algebraists.}

Let $(\Gamma^{g_1},\Gamma^{g_2})$, $(\Gamma^{g_2},\Gamma^{g_3})$ be
two  fixed ordered pairs of chain graphs. For 
$\alpha\in\aa(\Gamma^{g_1},\Gamma^{g_2})$, 
$\beta\in\aa(\Gamma^{g_2},\Gamma^{g_3})$, represented by single chord 
diagrams $x$, respectively $y$, let $\alpha\ast\beta$ denote the 
element of $\aa(\Gamma^{g_1},\sqcup_{g_2}S^1,\Gamma^{g_3})$, 
represented by the diagram obtained by attaching $\phi_{\ast}^{-1}y$ (the horizontal line of
$\Gamma^{g_2}$ removed) on top of $\phi_{\ast}^{-1}x$ (the horizontal line of
$\Gamma^{g_2}$ removed). For $g=0$ set $\ast$ to be the formal multiplication. 
Extend $\ast$ by linearity to formal 
power series of chord diagrams. Note that $\ast$ is associative.

Let $z_g\in\aa(\arrowg,\arrowg)$ and 
$z_g^N\in\aa^{\leq N}(\arrowg,\arrowg)$ be defined by:  

\begin{eqnarray}\label{zg}
z_g & = & \frac{Z(T_g)\otimes(\nu^{1/2})^{\otimes 2g}}{c_+^g\cdot c_-^g}\\
z_g^N & = & \left(\frac{Z(T_g)\otimes(\nu^{1/2})^{\otimes 2g}}
{c_+^g\cdot c_-^g}\right)^{[\leq N]} \label{zgN}
\end{eqnarray}

\noindent
where $T_g$ is the q-tangle from figure \ref{fig4}b with the 
non-associative structure
$(\dots(((\bullet\bullet)(\bullet\bullet))(\bullet\bullet))\dots)$,
$\nu=Z(\oo)\in\aa(\oOriented)$ is the Kontsevich integral of the 
zero-framed unknot, and $\otimes$ means taking the connected sum of 
chord diagrams on each of the $2g$ components, $c_+, c_-$ have been 
defined in 2.5. Note that 
$Z(T_g)\otimes(\nu^{1/2})^{\otimes 2g}$ has internal degree $0$.
For $g=0$ define $z_0=z_0^N=1$.


\medskip
\subsec{Proposition.}{\it 1) Let $\ast$ be the gluing operation defined above, let 
$\tilde\iota_N:\aa(\Gamma^{g_1}\sqcup(\sqcup_{2g_2}S^1)\sqcup\Gamma^{g_3})
\rightarrow\aa(\Gamma^{g_1}\sqcup\Gamma^{g_3})$ be the $\aa(\emptyset)$-linear
map defined by (\ref{in}), which refers exactly to all present circle 
components (in this case $2g_2$). Then\footnote{If one of elements 
belongs to the space of chord diagrams
on certain $2g$ arrows that connect $2g$ points on a "bottom line"
with $2g$ points on a "top line", as is $z_g$, we still can define $\ast$.
It extends by linearity and is associative.} 

\begin{equation}\label{ellN}
\ell^{\leq N}(\alpha,\beta)=(-1)^{gN}(\tilde\iota_{N}
(\alpha\ast z_{g_2}^N\ast\beta))^{[\leq N]}
\end{equation}

\noindent
defines a $\aa^{\leq N}(\emptyset)$-bilinear form
$\aa^{\leq N}(\Gamma^{g_1},\Gamma^{g_2})
\otimes\aa^{\leq N}(\Gamma^{g_2},\Gamma^{g_3})
\rightarrow\aa^{\leq N}(\Gamma^{g_1},\Gamma^{g_3})$.

2) For any $N$ and any respective elements 
$\alpha, \beta, \gamma$ we have 
$\ell^{\leq N}(\ell^{\leq N}(\alpha,\beta),\gamma)
=\ell^{\leq N}(\alpha,\ell^{\leq N}(\beta,\gamma))$.
}

\dem{ 1)
Let $\alpha\in\aa^{\leq N}(\Gamma^{g_1},\Gamma^{g_2})$,
$\beta\in\aa^{\leq N}(\Gamma^{g_2},\Gamma^{g_3})$,
$z_{g_2}^N\in\aa^{\leq N}(\Gamma^{g_2},\Gamma^{g_2})$.
$\ell^{\leq N}(\alpha,\beta)$ is well-defined. 
Indeed, it can be calculated in two steps: 
first $\tilde\pfi^{-1}(\alpha\ast z_{g_2}^N\ast\beta)$, then 
apply $\tilde q_{N}\circ\tilde\kappa_{N}$. 
Suppose $\alpha=0\in\aa^{\leq N}(\Gamma^{g_1},\Gamma^{g_2})$. 
Then, since $\tilde{\pfi}$ is an isomorphism,  
$\tilde{\pfi}^{-1}(\alpha\ast z_{g_2}^N\ast\beta)^{[\leq N]}=0$, i.e. if we factor in 
${\cal C}(\Gamma^{g_1},[m],\Gamma^{g_2})$
by the subspace spanned by diagrams with internal degree $>N$, 
hence $\tilde\iota_{N}\left(\alpha\ast z_{g_2}^N\ast\beta\right)^{[\leq N]}
=0\in\aa^{\leq N}(\Gamma^{g_1},\Gamma^{g_2})$. Similarly for $\beta$.
$\aa^{\leq N}(\emptyset)$-bilinearity of $(\alpha,\beta)\mapsto
\tilde\iota_{N}\left(\alpha\ast z_{g_2}^N\ast\beta\right)^{[\leq N]}$
is obvious. 

2) follows from the fact that $\ast$ is associative, and 
$\tilde\iota_{N}(\tilde\iota_{N}(\alpha\ast z_{g_1}^N\ast \beta)
\ast z_{g_2}^N\ast \gamma)=
\tilde\iota_{N}(\alpha\ast z_{g_1}^N\ast\tilde\iota_{N}(\beta
\ast z_{g_2}^N\ast \gamma))=
\tilde\iota_{N}(\alpha\ast z_{g_1}^N\ast\beta\ast z_{g_2}^N\ast \gamma)$.
}

\medskip
Note that when $g_1=0$, $\ell^{\leq N}$ becomes a
$\aa^{\leq N}(\emptyset)$-linear map: 
$$\ell^{\leq N}:\aa^{\leq N}(\Gamma^{g_2})\otimes\aa^{\leq N}(\Gamma^{g_2},\Gamma^{g_3})
\rightarrow\aa^{\leq N}(\Gamma^{g_3})$$

\noindent
Hence every element in $\aa^{\leq N}(\Gamma^{g_2},\Gamma^{g_3})$ defines 
a $\aa^{\leq N}(\emptyset)$-linear map from $\aa^{\leq N}(\Gamma^{g_2})$ to 
$\aa^{\leq N}(\Gamma^{g_3})$, and the induced map $\tilde\ell^{\leq N}:
\aa^{\leq N}(\Gamma^{g_2},\Gamma^{g_3})\rightarrow
\aa^{\leq N}(\Gamma^{g_2})^{\ast}\otimes\aa^{\leq N}(\Gamma^{g_3})$
is $\aa^{\leq N}(\emptyset)$-linear.
Composing with the isomorphism $\phi_{\ast}:\aa^{\leq n}(\arrowgi)\rightarrow\aa(\Gamma^{g_i})$
we obtain an $\aa^{\leq N}(\emptyset)$-linear map also denoted
$\tilde{\ell}^{\leq N}:\aa^{\leq N}(\Gamma^{g_1},\Gamma^{g_2})\rightarrow
\aa^{\leq N}(\arrowgOne)^{\ast}\otimes\aa^{\leq N}(\arrowgTwo)$.

Extend $\Gamma^g\rightarrow(\Gamma^g,\emptyset)$ to 
$\aa^{\leq N}(\Gamma^g)\rightarrow\aa^{\leq N}(\Gamma^g,\emptyset)$, 
and compose with $\tilde\ell^{\leq N}$ to obtain 
a $\aa^{\leq N}(\emptyset)$-linear map
$(\ell^{\leq N})^{\ast}:\aa^{\leq N}(\Gamma^g)\rightarrow
\aa^{\leq N}(\Gamma^g)^{\ast}$. 
Namely $(\ell^{\leq N})^{\ast}(\beta)(\alpha)=
\ell^{\leq N}(\alpha,\beta)$, 
$\forall\alpha,\beta\in\aa^{\leq N}(\Gamma^g)$.
Similarly, there is a map
$(\ell^{\leq N})^{\ast}:\aa^{\leq N}(\arrowg)\rightarrow
\aa^{\leq N}(\arrowg)^{\ast}$

The second part of the above proposition shows that 
$\tilde\ell^{\leq N}(\ell^{\leq N}(\beta,\gamma))
=\tilde\ell^{\leq N}(\gamma)\circ\tilde\ell^{\leq N}(\beta)$
for any corresponding $\beta,\gamma$, i.e. the following diagram is
commutative:

\begin{eqnarray*}
\aa^{\leq N}(\Gamma^{g_1},\Gamma^{g_2})\otimes\aa^{\leq N}(\Gamma^{g_2},\Gamma^{g_3})&
\stackrel{\ell^{\leq N}}{\longrightarrow}&\aa^{\leq N}(\Gamma^{g_1},\Gamma^{g_3})\\
\downarrow\tilde\ell^{\leq N}\otimes\tilde\ell^{\leq N} & &
\downarrow\tilde\ell^{\leq N}\\
\aa^{\leq N}(\Gamma^{g_1})^{\ast}\otimes\aa^{\leq N}(\Gamma^{g_2})\otimes
\aa^{\leq N}(\Gamma^{g_2})^{\ast}\otimes\aa^{\leq N}(\Gamma^{g_3})&
\stackrel{evaluation}{\longrightarrow}&
\aa^{\leq N}(\Gamma^{g_1})^{\ast}\otimes\aa^{\leq N}(\Gamma^{g_3})
\end{eqnarray*}

\medskip\noindent
{\it Remark.} If we were to use Knizhnik-Zamolodchikov or any other associator, in the definition of $\ell$,
between $\alpha$, $z_g$ and $\beta$ we would have to insert an element $A$ in the space of chord diagrams on $2g$
arrows alternatively oriented downward and upward, and its horizontal reflection, such that 
$\left(\phi_{\ast}^{-1}Z(\tangleEight)\right)\ast A=\nu^{\otimes g}$; and similarly for $\beta$.
For the even associator, $A$ can be taken $1$, i.e. it can be omitted.
Conjecture A above claims that $A=1$ for any associator.


\medskip
\subsec{The categories $\aa^{\leq N}$ and $\aa$.}Let 
$\aa^{\leq N}$ be the category with objects $\aa^{\leq N}(\arrowg)\equiv\aa^{\leq N}(\Gamma^0,\arrowg)$, 
$g\geq 0$, and morphisms the set of $\aa^{\leq N}(\emptyset)$-homomorphisms between these modules. 
Similarly define the category $\aa$. Note that via the isomorphism $\phi_{\ast}$ in
Proposition 2.3 we can identify $\aa(\Gamma_g)$ and $\aa(\arrowg)$.

\begin{figure}[htbp]
\centerline{\psfig{file=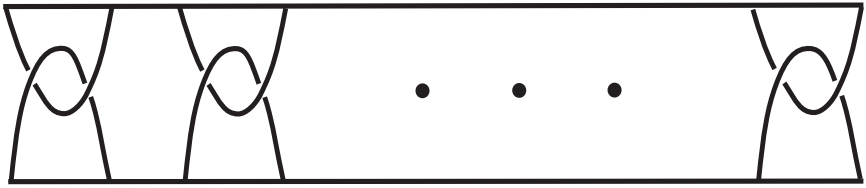,width=12cm,height=3cm,angle=0}}
\caption{\sl The embedded framed graph $W_g$} \label{figWg}
\end{figure}

\medskip
\subsec{Proposition.}{\it Let $\ell^{\leq N}$ be the bilinear form defined
by the previous proposition, and $\tilde\ell^{\leq N}$ be the induced 
map $\aa^{\leq N}(\Gamma^{g_1},\Gamma^{g_2})\rightarrow
\aa^{\leq N}(\arrowgOne)^{\ast}\otimes\aa^{\leq N}(\arrowgTwo)
=Hom(\aa^{\leq N}(\arrowgOne),\aa^{\leq N}(\arrowgTwo))$. 
Denote $w_g=Z(W_g)=Z(W_g)^{[\leq N]}\in\aa^{\leq N}(\Gamma^g,\Gamma^g)$,
where $W_g$ is the embedded framed graph in figure \ref{figWg}, where the first $\Gamma^g$
corresponds to the lower of the two chain graphs in the picture.
(For $g=0$, $W_g=1$.) Then $\tilde\ell^{\leq N}(w_g)$ is the identity operator on 
$\aa^{\leq N}(\arrowg)$.
}

\medskip
Note that $w_g$ has zero as internal degree $\geq 1$ parts.

\medskip
\subsec{Theorem.}{\it 1) There is a (unique) 
$\aa(\emptyset)$-bilinear form
$\ell:\aa(\Gamma^{g_1},\Gamma^{g_2})
\otimes\aa(\Gamma^{g_2},\Gamma^{g_3})
\rightarrow\aa(\Gamma^{g_1},\Gamma^{g_3})$, such that
its $\leq N$ internal degree part coincides with $\ell^{\leq N}$.

2)  Let $\tilde\ell$ be the induced map 
$\aa(\Gamma^{g_1},\Gamma^{g_2})\rightarrow
\aa(\arrowgOne)^{\ast}\otimes\aa(\arrowgTwo)=
Hom_{\aa(\emptyset)}(\aa(\arrowgOne),\aa(\arrowgTwo))$, and 
denote as before $w_g=Z(W_g)\in\aa(\Gamma^g,\Gamma^g)$,
where $W_g$ is shown in figure \ref{figWg}. 
Then $\tilde\ell(w_g)$ is the identity operator 
on $\aa(\arrowg)$.
}

\medskip
We will prove these statements in Section 3 (after 3.6, resp. 3.7).

%
%

\section{The TQFT}
\setcounter{equation}{0}
\setcounter{nsubsec}{1}

Now we can formally construct the truncated and full TQFTs. 
We will show that they are non-degenerate and anomaly-free. 
The truncated TQFTs are with respect to the internal 
degree. Since the map 
$\tilde\iota_N$, which we had to introduce if we want to have invariance 
under Kirby moves for chain graphs, decreases the total degree 
of a diagram by $2gN$, and since $\tilde\iota_N$ must be applied 
every time we glue two cobordisms, one should not expect the theory 
to truncate with respect to the total degree of chord diagrams.

For our purposes a TQFT $(\T,\tau)$ based on the cobordism category 
$\mathfrak Q$ (or a subcategory of it) is 1) a covariant functor $\T$ 
from the category those objects are the objects of ${\mathfrak Q}$ 
(i.e. natural numbers) and morphisms are the 
homeomorphisms of parametrized surfaces to a subcategory
${\mathfrak V}_K$ of the full category of $K$-modules, such that 
$\T(0)=K$, where $K$ is a commutative module; 
and 2) a map $\tau$ that associates
to each 3-cobordism $(M,f_1,f_2)$ a $K$-homomorphism 
$\tau(M):\T(\Sigma_1)\rightarrow\T(\Sigma_2)$,
satisfying the following axioms:

\it
\medskip
\begin{tqftax}{0mm}{10mm}

\item (Naturality) If $(M_1,\Sigma_1,\Sigma_1^{\prime})$, 
$(M_2,\Sigma_2,\Sigma_2^{\prime})$ are two 3-cobordisms, 
and $f:M_1\rightarrow M_2$ is a homeomorphism of 3-cobordisms, 
preserving the parametrizations, then the following diagram is 
commutative:
$$\T(\Sigma_1)\stackrel{\tau(M_1)}{\longrightarrow}
\T(\Sigma_1^{\prime})$$
$$\T(f|_{\Sigma_1})\downarrow\hspace{2cm}
\downarrow\T(f|_{\Sigma_2})$$
$$\T(\Sigma_2)\stackrel{\tau(M_2)}{\longrightarrow}
\T(\Sigma_2^{\prime})$$

\item (Functoriality) If $M_1$, $M_2$ are 3-cobordisms, 
$f=f_2\circ(f_1^{\prime})^{-1}:\partial_{top}(M_1)\rightarrow
\partial_{bottom}(M_2)$ is the gluing homeomorphism, and denote 
$M=M_2\cup_f M_1$, then $\tau(M)=k\cdot\tau(M_2)\circ\tau(M_1)$.
$k\in K$ is called {\sc the anomaly}.

\item (Normalization) Let 
$(\Sigma\times[0,1],(\Sigma\times0,p_1),(\Sigma\times1,p_2))$ 
be the 3-cobordism mentioned in 1.2, then  
$$\tau(\Sigma_g\times[0,1],(\Sigma_g\times0,p_1),
(\Sigma_g\times1,p_2))=id_{\T(\Sigma_g)}$$

\item (pseudo-Hermitian structure) There is a superstructure on each 
element $V$ of ${\mathfrak V}_K$, i.e. it admits an antimorphism 
$\overline{\cdot}:V\rightarrow V$ (a map linear in 0-supergrading 
and antilinear in 1-supergrading), that commutes with ($=$ is natural 
with respect to) surface homeomorphisms. There is a canonical map $V\rightarrow V^{\ast}$,
which composed with the above antimorphism extends (from the particular 
case when $\T(\Sigma_1)=K$) to an antimorphism 
$\overline{\cdot}:Mor(\T(\Sigma_1),\T(\Sigma_2))\rightarrow
Mor(\T(-\Sigma_2),\T(-\Sigma_1))$, that commutes with 
homeomorphisms of 3-cobordisms, such that
$$\tau(-M)=\overline{\tau(M)}$$
\end{tqftax}

\rm
\noindent
We can not require multiplicativity or self-duality since in 
the category $\mathfrak Q$ all cobordisms are connected. 
Conditions (A1-A3) say that $\tau:{\mathfrak Q}\rightarrow\cal A$ is a
pseudo-functor. $\tau$ would is a true functor when there is no anomaly.
If the set of $\tau(M,S^2,\Sigma)$'s, spans (in the closure for infinite-dimensional
modules) $\T(\Sigma)$, the  TQFT is called {\it non-degenerate}.


\medskip
\subsec{$\tau^{\leq N}$ and $\tau$.}The 
definition of $\Omega_n(S^3_L)$ and the proof of its invariance under Kirby 
moves have been extended [\ref{MO}, proposition 2.4] to an invariant
of embeddings $L\sqcup G\hookrightarrow S^3$ and extended (generalized)
Kirby moves. Again, as long as $Z(L\cup G)$ has been shown well-defined, it does 
not matter which associator we use.

\medskip
\te{Proposition}{ 1) Let $n\in\nn$, and let $L\sqcup G\hookrightarrow S^3$ 
be an arbitrary embedding of a link and a (union of) chain graph(s)
in $S^3$, $\sigma_+$, $\sigma_-$ be the number of positive, 
respectively negative eigenvalues of $lk(L)$, the linking matrix 
of $L$, $g$ the total number of circle components of $G$. Define:

$$\Omega_n(L,G) := 
\left(\frac{\tilde{\iota}_n(\cZ(L\cup G))}{\iota_n(\cZ(\oo^{+1}))
^{\sigma_+} \cdot  \iota_n(\cZ(\oo^{-1}))^{\sigma_-}}\right)
^{[\leq n]}\in\aa^{\leq n}(\Gamma)\subset\aa(\Gamma)$$

\noindent 
where $\Gamma$ is $G$ as an abstract graph.
Then for every $m\leq n$ the internal degree part $\Omega_n(L,G)^{[m]}$
is invariant under extended (generalized) KI and KII moves, 
and under orientation change of components of $L$.

2) With the above notations, and denoting 
$d_L=|det(lk(L))|$, the following relation holds in $\aa(\Gamma)$ for any 
(not necessarily connected) chain graph $G$ and link $L$:

\begin{equation}\label{eqvogel.tilda}
\left(\tilde\iota_{n+1}\cZ(L,G)\right)^{[\leq n]} 
= 
(-1)^{|L|}det(lk(L))\left(\tilde\iota_n\cZ(L,G)\right)^{[\leq n]}
\end{equation}

\noindent
and therefore:

$$\Omega_{n+1}(L,G)^{[\leq n]}=d_L\cdot\Omega_n(L,G)$$

3) $\frac{1}{d_L}\Omega_n(L,G)$ is a group-like element of 
$\aa^{\leq n}(\Gamma)\subset\aa(\Gamma)$ of the form $1+$ higher order terms.
}

\dem{1) By the well defineness of the internal degree parts (see 2.8), it is 
enough to show that $\Omega_n(L,\Gamma)$ stays invariant under the moves, which is precisely 
the statement in Proposition 2.4 in [\ref{MO}] and the remarks after it. 
There the proof is similar to the case of links (see Section 3 in 
[\ref{LMO}], or Section 4 in [\ref{vogel}]). 

2) Follow the proof of Proposition 4.5 in [\ref{vogel}].

3) First note that $Z$ of any elementary pseudo-quasi-tangle is group-like of 
the desired form. Indeed, if one uses KZ associator, the elements $a, b$ in [\ref{MO}, page 503] used in 
the definition of $Z$ for the vicinity of trivalent vertices are clearly so. 
If one uses even associator, then (even simpler)
it follows from the fact that $\Delta\nu=\nu\otimes\nu$ and $\nu=1+h.o.t.$
Hence $\cZ(L\cup G)$ is group-like of the form $1+h.o.t.$ for any 
$L\cup G\hookrightarrow S^3$ (compare with [\ref{LMO}, subsection 1.4]). 
That $\Delta$ commutes with $\tilde\iota_n$ follows from the fact that $\Delta$ 
commutes with $\tilde\pfi$, and an explicit calculation of 
$\Delta\circ(\tilde q_n\circ\tilde\kappa_n)$ and 
$(\tilde q_n\circ\tilde\kappa_n)\otimes(\tilde q_n\circ\tilde\kappa_n)\circ\Delta$
 for any diagram with $2n$ legs of each colour $1,\dots,|L|$, just as in the case 
$G=\emptyset$ [\ref{vogel}, \ref{LMO}]. Similarly it follows that 
$\frac{1}{d_L}\Omega_n(L,G)$ has the form $1+h.o.t.$ 
(compare with [\ref{LMO}, Lemma 4.7]).
}

\medskip
Let $M$ be a morphism in ${\mathfrak Q}$ between $g_1$ and $g_2$. 
Let $(L,G_1,G_2)$ be such that $\kappa(L,G_1,G_2)=(M,f_1,f_2)$. By 
Proposition 2.1 in [\ref{MO}], the ambiguity in this choice 
is a finite sequence of extended KI and KII moves, and change of 
orientation of link components. In Theorem 1.4 of [\ref{MO}] 
using KZ associator, or alternatively in 2.6 here using even associator, 
the Kontsevich integral is extended to an isotopy invariant of chain graphs in $S^3$, 
and hence of embeddings $L\cup G_1\cup G_2$ in $S^3$.
Suppose that $G_1$, $G_2$ regarded as abstract graphs are $\Gamma^{g_1},\Gamma^{g_2}$. 
Then let us define:

\begin{eqnarray}\label{tau}
\tau(M,f_1,f_2) & = & \sum\limits_{n\geq 0}
\frac{1}{|det(lk(L))|^n}\cdot
\left(\frac{\tilde \iota_n(\cZ(L\cup G_1\cup G_2))}
{\iota_n(\cZ(\oo^{+1}))^{\sigma_+}\cdot\iota_n(\cZ(\oo^{-1}))
^{\sigma_-}}\right)^{[n]}\nonumber \\
& \stackrel{(\ref{eqvogel})}{=} & \sum\limits_{n\geq 0}
\frac{(-1)^{\sigma_+n}}{|det(lk(L))|^n}\cdot
\left(\frac{\tilde \iota_n(\cZ(L\cup G_1\cup G_2))}
{c_+^{\sigma_+}\cdot c_-^{\sigma_-}}
\right)^{[n]}\in\aa(\Gamma^{g_1},\Gamma^{g_2})
\end{eqnarray}

\noindent
where $[n]$ represents the internal degree part, 
$\sigma_+,\sigma_-$ are the number of positive and negative 
eigenvalues of $lk(L)$, and $\tilde  \iota_n$ refers to the circle 
components of chord diagrams, all coming here from the components 
of the link $L$. As before (see the proof of Proposition 2.3) we can 
assume that the vertices of  chord diagrams are off the horizontal lines. 
$c_+$, $c_-$ have been defined in 2.5.
We use the convention $det(lk(\emptyset))=1$. Also, let:

\begin{eqnarray}\label{tauN}
\tau^{\leq N}(M,f_1,f_2) & = & \frac{(-1)^{\sigma_+N}}{d(M)^{N}}\cdot
\left(\frac{\tilde \iota_{N}(\cZ(L\cup G_1\cup G_2))}
{c_+^{\sigma_+}\cdot c_-^{\sigma_-}}
\right)^{[\leq N]}\in\aa^{\leq N}(\Gamma^{g_1},\Gamma^{g_2})
\end{eqnarray}

\noindent
where $d(M)=|H_1(\widehat{M},\zz)|$.
Note that $\tau^{\leq N}(M,f_1,f_2)=\tau(M,f_1,f_2)^{[\leq N]}$.
Proposition 3.1.1) (also see Proposition 2.4 in [\ref{MO}]) implies 
that $\tau(M,f_1,f_2)$ and $\tau^{\leq N}(M,f_1,f_2)$ , 
being invariant under extended KI and KII moves and change of orientation 
of link  components, are independent of the choice of the triplet $(L,G_1,G_2)$.
Hence  $\tau(M,f_1,f_2)\in\aa(\Gamma^{g_1},\Gamma^{g_2})$ and
$\tau^{\leq N}(M,f_1,f_2)\in\aa^{\leq N}(\Gamma^{g_1},\Gamma^{g_2})$ 
are invariants of 3-cobordisms of the category ${\mathfrak Q}$. 
As we have seen in 2.9, 
$\tilde\ell^{\leq N}(\tau^{\leq N}(M,f_1,f_2))$ is then a 
$\aa^{\leq N}(\emptyset)$-homomorphisms 
from $\aa^{\leq N}(\arrowgOne)$ to $\aa^{\leq N}(\arrowgTwo)$.


\medskip
\subsec{Non-degeneracy of the TQFT(s).}We 
can define now similarly (\ref{tau}) and (\ref{tauN}) for cobordisms $(M,\emptyset,f)$ with
only one connected and parametrized boundary component, as long as $\widehat{M}$ is a $\qq$-homology
sphere. We think of the boundary as the top of the 3-cobordism. Then $\tau(M,\emptyset,f)\in\aa(\Gamma^g)$,
and hence via the isomorphism $\phi_{\ast}^{-1}$ we obtain associated to $(M,\emptyset,f)$ an
element in $\aa(\arrowg)$. (Similarly add ${}^{\leq N}$.)

\te{Theorem.}{Let $\aa_{\tau}^{\leq N}(\arrowg)$, respectively $\aa_{\tau}(\arrowg)$, be the
$\qq$-vector subspace of $\aa^{\leq N}(\arrowg)$,
respectively $\aa(\arrowg)$, generated by all $\phi_{\ast}^{-1}\tau^{\leq N}(M,\emptyset,f)$,
respectively all $\phi_{\ast}^{-1}\tau(M,\emptyset,f)$, such that $\widehat{M}$ is a $\zz$-homology
sphere. Then the completion of $\aa_{\tau}^{\leq N}(\arrowg)$ is $\aa^{\leq N}(\arrowg)$, and
the completion of $\aa_{\tau}(\arrowg)$ is $\aa(\arrowg)$.}

\medskip
It is known [\ref{HM}, Proposition 13.1] that for every chord diagram $\xi\in\aa(\arrowg)$ 
of degree $m$, with connected dashed graph, there exist string links $L^{\pm}$ such that
$Z(L^{\pm})=1\pm\xi+o(m+1)$. For a very intuitive geometric realization of $L^{\pm}$ see also [\ref{habiro},
\ref{LeGrenoble}].


\medskip
\subsec{Lemma.}{\it For every $n\geq 0$ and every chord diagram $\xi\in\aa(\arrowg)_{\leq n}$, 
there exist string links $L_1,\dots ,L_k$ and positive integers $a_1,\dots,a_k$ such that 
$\sum\limits_{i=1}^{k}a_iZ(L_i)=\xi+o(n+1)$.}

\dem{Induction on $n$. For $n=0$, $Z($trivial string link$)=1\in\aa(\arrowg)$. For $n=1$, $\xi$ must have
connected dashed graph, hence the claim follows from the mentioned result of Habegger and Masbaum, because
$Z($trivial string link$)=1$. 
For general $n$, suppose $\xi$ has degree $m$. We prove the statement first
for $m=n$, then for $m=n-1,\dots,1$ ($m=0$ is obvious). 
For arbitrary $m$, by the same argument as in the 
proof of Proposition 2.4.5) we can assume $\xi=\sum\pm\xi_1\bullet\dots\bullet\xi_k$, 
$\xi_i$ have connected dashed graph and degree $\geq 1$. 
If $k>1$, by the induction hypothesis there exist $\alpha^i_j\in\zz$ and string links $L^i_j$ such
that $\sum_ja^i_jZ(L^i_j)=\xi_i+o(n)$, $\forall i$. Therefore 
$\sum\limits_{i_1,\dots,i_k}a_1^{i_1}\cdots a_k^{i_k}\cdot Z(L_1^{i_1}\bullet\dots\bullet L_k^{i_k})=
\sum\limits_{i_1,\dots,i_k}a_1^{i_1}Z(L_1^{i_1})\bullet\dots\bullet a_k^{i_k}Z(L_k^{i_k})=
\left(\xi_1+o(n)\right)\bullet\dots\bullet\left(\xi_k+o(n)\right)=
\xi_1\bullet\dots\bullet\xi_k+o(n+1)$. 
If $k=1$, by Habegger-Masbaum result, there is a string link $L$ such that 
$Z(L)-Z($trivial string link$)=\xi_i+o(m+1)$. Therefore the statement for $m$ follows from the fact
that it holds for $m+1$ (express in the later formula the degree $m+1$ terms of $o(m+1)$). 
If $k=1$ and $m=n$, it is precisely the Habegger-Masbaum result.
Note that all coefficients $a_i$ appearing throughout the proof can be arranged positive or
negative as we wish [\ref{habiro},\ref{LeGrenoble}], hence the ones in the statement can be ensured positive.
}

\medskip
It is known [\ref{LeGrenoble}, Theorem 4.5; see also \ref{HO1}] that for any connected trivalent graph $D$ of 
degree $n$ there exist $\zz$-homology 3-spheres $M^{\pm}$ such that 
$Z^{LMO}(M^{\pm})=1\pm D+o(n+1)\in\aa(\emptyset)$. (This is proved there for $Z^{lmo}$, but it is obviously
then true for $Z^{LMO}$.) Since $Z^{LMO}(S^3)=1\in\aa(\emptyset)$, with a proof absolutely similar
to the one above, we have:


\medskip
\subsec{Lemma.}{\it For any $n\geq 0$ and any chord diagram $\xi\in\aa(\emptyset)_{\leq n}$, there exist 
$\zz$-homology spheres $M_1,\dots ,M_k$ and positive integers $b_1,\dots,b_k$ such that 
$\sum\limits_{i=1}^{k}b_iZ^{LMO}(M_i)=\xi+o(n+1)$. In particular the set
$\{\sum_ib_iZ^{LMO}|M_i\:\:\zz$-homology sphere, $b_i\in\nn^{\ast}\}$ is dense in $\aa(\emptyset)$.
}

\medskip
\demlung{Proof of Theorem 3.2.}{Let $\arrowtwog$ denote the graph with $2g$ edges oriented alternatively 
down- and upward. Lemma 3.3 is clearly true for $\xi\in\aa(\arrowtwog)$ as well. Therefore for every
$n\geq 0$ and every $\beta\in\aa(\arrowtwog)$ there exist string links $L_i$ and $a_i\in\qq$ such that
$\sum a_i\cdot Z(L_i)=\beta\bullet\left(\arrowdown\otimes\nu^{-1}\right)^{\otimes g}+o(n+1)$. Using
the operation $\ast$ defined in 2.8 attach $Z(\tangleNine)$ on top and $Z(\UP\dots\UP)$
below each side of this equality to obtain the existence of embedded framed graphs $G_i$ and 
$a_i\in\qq$ such that $\sum a_i\cdot Z(G_i)=Z(\tangleNine)
\ast\left(\beta\bullet\left(\arrowdown\otimes\nu^{-1}\right)^{\otimes g}\right)
\ast Z(\UP\dots\UP)+o(n+1)=\widehat{\beta}+o(n+1)$,
where $\left(\beta\mapsto\widehat{\beta}\right):\aa(\arrowtwog)\rightarrow\aa(\oGraph)$ 
is the map induced by inclusion. The later is well-defined, since at the level of
$\dd(\arrowtwog)$ AS, IHX and STU relations are sent to some of the same relations on on $\oGraph$.
It is clearly surjective by Proposition 2.3. Therefore for every $n\geq 0$ and $\alpha\in\aa(\arrowg)$
there exist $G_i$ and $a_i\in\qq$ such that $\sum a_i\cdot Z(G_i)=\phi_{\ast}(\alpha)+o(n+1)$.

Let $N=0$. Then by (\ref{tauN}) we have $\tau^{\leq 0}(M,\emptyset,f)=\tilde{\iota}_0\cZ(L\cup G)$
whenever $\kappa(L,G)=(M,\emptyset,f)$. Since $\tilde\iota_N$ refer only to link components, for every embedded
framed graph $G$ we obtain 
$\tau^{\leq N}(\kappa(\emptyset,G))=\tau^{\leq 0}(\kappa(\emptyset,G))=\tilde\iota_0\cZ(G)
=\cZ(G)=Z(G)$. Together with the conclusion of the previous paragraph this shows that for every
$n\geq 0$ and  every $\alpha\in\aa^{\leq 0}(\arrowg)=\aa_c(\arrowg)$ there exist cobordisms
$(M_i,\emptyset,f_i)$ having $\widehat{M_i}=S^3$, and $\a_i\in\qq$ such that 
$\sum_i a_i\cdot\tau^{\leq 0}(M_i,\emptyset ,f_i)=\phi_{\ast}(\alpha)+o(n+1)$. This proves the theorem
for $N=0$.

Let $N>1$. Recall that $\aa^{\leq N}(\arrowg)=\aa^{\leq N}(\emptyset)\otimes_{\qq}\aa_c(\arrowg)$.
Therefore the statement is enough to prove for $\xi\cdot\alpha$, $\xi\in\aa^{\leq N}(\emptyset)$
and $\alpha\in\aa_c(\arrowg)$. By Lemma 3.4 there exist $\zz$-homology spheres $M_i$ and $b_i\in\nn^{\ast}$
such that $\sum b_i Z^{LMO}(M_i)=\xi+o(N+1)$. Then, by the previous paragraph, for every $n\geq 0$
there exist cobordisms 
$(M_j,\emptyset,f_j)$ having $\widehat{M_j}=S^3$ and $a_i\in\qq$ such that
$\left(\sum_i b_i Z^{LMO}(M_i)\right)^{[\leq N]}\cdot$
$\cdot\left(\sum_j a_j\tau^{\leq 0}(M_j,\emptyset,f_j)\right)
=\xi\cdot\phi_{\ast}(\alpha)+o(n+1)$. Therefore:
$$\phi_{\ast}(\xi\cdot\alpha)=\xi\cdot\phi_{\ast}(\alpha)
=\sum\limits_{i,j}b_ia_j Z^{LMO}(M_i)^{[\leq N]}\tau^{\leq 0}(M_j,\emptyset,f_j)+o(n+1)
=\sum\limits_{i,j}b_ia_j \tau^{\leq N}(M_i)\tau^{\leq 0}(M_j,\emptyset,f_j)+$$
$$+o(n+1)
=\sum\limits_{i,j}b_ia_j\left(\tau^{\leq N}(S^3_{L_i})\tau^{\leq 0}(\kappa(\emptyset,G_j))\right)+o(n+1)
=\sum\limits_{i,j}b_ia_j \tau^{\leq N}\left(S^3_{L_i}\#\kappa(\emptyset,G_j)\right)+o(n+1),$$
where $b_ia_j\in\qq$, $M_i=S^3_{L_i}$, $(M_j,\emptyset,f_j)=\kappa(\emptyset,\G_j)$ and
$S^3_{L_i}\#\kappa(\emptyset,G_j)=\kappa(L_i\sqcup G_j)$, whose filling is $S^3_{L_i}$, 
a $\zz$-homology sphere. This proves the theorem for arbitrary $N$.

Since for any cobordism $(M,\emptyset,f)$ we have $\tau(M,\emptyset,f)^{[\leq N]}=\tau^{\leq N}(M,\emptyset,f)$,
by taking $N=n$, we can see that for every  $n\geq 0$ and any $\alpha\in\aa(\arrowg)$ there exist $c_i\in\qq$
and $(M_i,\emptyset,f_i)$ with $\widehat{M_i}$ $\zz$-homology spheres, such that
$\phi_{\ast}(\alpha)=\sum\limits_ic_i\tau(M_i,\emptyset,f_i)^{[\leq n]}+o(n+1)
=\sum\limits_ic_i\tau(M_i,\emptyset,f_i)+o(n+1)$.
}

\medskip\noindent
{\it Remark.} The statement of Theorem 3.2 remains true for $\arrowg$ replaced by any $\Gamma=$ union of chain graphs, 
since in the above proof we only used that the "capping map" $\beta\mapsto\widehat{\beta}$ is
linear and surjective, which is clear by the remark after Proposition 2.3.


\medskip
\subsec{Lemma.}{\it For any $\beta\in\aa(\Gamma^{g_2},\Gamma^{g_3})$ and
any sequence of elements $\alpha^n\in\aa(\Gamma^{g_1},\Gamma^{g_2})$, such that for every $n$, 
$(\alpha^n)_{\leq n}=(\alpha^{n+1})_{\leq n}=\dots$, both sides of the following equalities are
well-defined and the equalities holds:

\begin{equation}\label{passingToLimitN}  
\ell^{\leq N}(\lim\limits_n\alpha^n,\beta) = \lim\limits_n\ell^{\leq N}(\alpha^n,\beta)
\end{equation}

\noindent
A similar property holds for the r\^ ole of two arguments of $\ell^{\leq N}$ reversed.
}

\dem{ The existence of $\alpha:=\lim\limits_n\alpha^n\in\aa(\Gamma^{g_1},\Gamma^{g_2})$ 
follows directly from the fact that we defined the topology on $\aa(\Gamma^{g_1},\Gamma^{g_2})$
such that $dist(p,q)<\frac{1}{2^n}$ if and only if $p-q$ has degree $>n$. Then, since
$\alpha_{>n+gN}$ does not contribute to $\ell^{\leq N}(\alpha,\beta)_{\leq n}$, we have:

$$\lim\limits_n\ell^{\leq N}(\alpha,\beta)_{\leq n}
=\lim\limits_n\ell^{\leq N}(\alpha_{\leq n+gN},\beta)_{\leq n}
=\lim\limits_n\ell^{\leq N}(\alpha_{\leq n+gN},\beta)
=\lim\limits_m\ell^{\leq N}(\alpha_{\leq m},\beta)=$$
$$=\ell^{\leq N}(\lim\limits_m\alpha_{\leq m},\beta)
=\ell^{\leq N}(\alpha,\beta)$$

\noindent
The existence of the third limit and the second equality follow from a standard Cauchy-sequences argument.
The fourth equality is true since $\lim\limits_m$ commutes with $\ast$, $\widetilde{\iota_N}$ and $^{[\leq N]}$.
On the other hand $\ell^{\leq N}(\alpha^n,\beta)$ and $\ell^{\leq N}(\alpha,\beta)$ agree in degree $\leq n-2gN$. 
Hence $\lim\ell^{\leq N}(\alpha^n,\beta)=\lim\ell^{\leq N}(\alpha,\beta)_{\leq n}$.
Putting the two together we obtain (\ref{passingToLimitN}). 
}

\medskip\noindent
{\it Remark.} If in the statement of this Lemma we assume that $\lim\limits_n\alpha^n$ exists,
then we can relax the topology:
{\it distance}$(p,q)\leq\frac{1}{n}\Leftrightarrow p-q$ has no terms of degree $<n$.


\medskip
\subsec{The functors $\mathfrak Q\rightarrow{\cal A}^{\leq N}$ 
and $\mathfrak Z\rightarrow{\cal A}^{\leq N}$: gluing formula and normalization.}Set 
$\T^{\leq N}(g)=\aa^{\leq N}(\Gamma^g)$, if $g>0$,
$\T^{\leq N}(0)=\aa^{\leq N}(\emptyset)$.
In this case $K=\aa^{\leq N}(\emptyset)$.
Set $\T(g)=\aa(\Gamma^g)$, if $g>0$, and $\T(0)=\aa(\emptyset)$.
In this case $K=\aa(\emptyset)$.
Now, let us start verifying the axioms of TQFT. 
Set $\T(f|_{\Sigma})=id_{\aa(\Gamma^g)}$ for any homeomorphism $f$ of the 
parametrized surfaces. 
Then $\T$ is a covariant functor, and the naturality axiom (A1) is obvious. 
The same is true for $\T^{\leq N}$. We will derive now a gluing formula.

\te{Theorem.}{1) Let $(M_1,f_1,f_1^{\prime})$ and 
$(M_2,f_2,f_2^{\prime})$ be two  3-cobordisms. 
Suppose $(M_1, f_1, f_1^{\prime}) = \kappa(L_1,G_1,G_1^{\prime})$,
$(M_2,f_2,f_2^{\prime}) = \kappa(L_2,G_2,G_2^{\prime})$, 
and $(M_2\cup_{f_2\circ(f_1^{\prime})^{-1}} M_1, f_1,
f_2^{\prime}) = \kappa(L_1\cup L_0\cup L_2, G_1,G_2^{\prime})$,
the later triplet obtained from the previous two by the construction
described in Proposition 1.3.
Denote $\sigma^1_+=sign_+(lk(L_1))$, $\sigma^2_+=sign_+(lk(L_2))$,
$\sigma_+=sign_+(lk(L_1\cup L_0\cup L_2))$, and let $g$ be the
genus of the connected closed surface along which is this splitting.
Then the integer $s(M,M_1,M_2)=\sigma^1_++\sigma^2_++g-\sigma_+$ 
is an invariant of the decomposition $M=M_2\cup_{f_2\circ(f_1^{\prime})^{-1}}M_1)$, 
i.e. it does not depend on the choice of triplets 
representing the 3-cobordisms $M_1$ and $M_2$.

2) Let $(M_1,f_1,f_1^{\prime})$ and 
$(M_2,f_2,f_2^{\prime})$ be two $\qq$HH. Denote
$d=|H_1(\widehat{M_2\cup_{f_2\circ(f_1^{\prime})^{-1}}M_1},\zz)|$, 
$d_1=|H_1(\widehat{M_1},\zz)|$, $d_2=|H_1(\widehat{M_2},\zz)|$. Suppose that these cobordisms 
are glued along a surface of genus $g$. Then:
\begin{eqnarray}\label{anomaly1}
\tau^{\leq N}(M_2\bigcup_{f_2\circ(f_1^{\prime})^{-1}}M_1,f_1,f_2^{\prime}) & = &
\left((-1)^N\left(\frac{c_+}{c_-}\right)^{[\leq N]}\right)^{\sigma_+^1+\sigma_+^2+g-\sigma_+}
\left(\frac{d_1d_2}{d}\right)^N \cdot \nonumber\\  
& \cdot & \ell^{\leq N}(\tau^{\leq N}(M_1,f_1,f_1^{\prime}),
\tau^{\leq N}(M_2,f_2,f_2^{\prime}))
\end{eqnarray}

\smallskip
\noindent
where $(-1)^N\cdot\left(c_+/c_-\right)^{[\leq N]}\in\aa^{\leq N}(\emptyset)$,
the multiplication by scalars is thought
in the category $\aa^{\leq N}$,
and $\sigma_+^1+\sigma_+^2+g-\sigma_+$ is an integer.
}

\dem{1) Proposition 11 of [\ref{CL}].

2) Let $(L_1,G_1,G_1^{\prime})$, $(L_2,G_2,G_2^{\prime},)$, and
$(L_1\cup L_0\cup L_2, G_1,G_2^{\prime})$ be as above, 
and let $(\sigma_+,\sigma_-)$, $(\sigma_+^1,\sigma_-^1)$, 
respectively $(\sigma_+^2,\sigma_-^2)$ 
be the signatures of  $lk(L_1\cup L_0\cup L_2)$, $lk(L_1)$, resp. $lk(L_2)$.
Then, temporarily abbreviating $c_+^{[\leq N]}$ and $c_-^{[\leq N]}$
to $c_+$ and $c_-$: 

\begin{eqnarray*}
\tau^{\leq N}(M_2\bigcup_{f_2\circ(f_1^{\prime})^{-1}}M_1,f_1,f_2^{\prime}) 
& = & \frac{(-1)^{\sigma_+N}}{d^{N}}\cdot
\left(\frac{\tilde\iota_{N}
\cZ(L_1\cup L_0\cup L_2\cup G_1\cup G_2^{\prime})}
{c_+^{\sigma_+}\cdot c_-^{\sigma_-}}
\right)^{[\leq N]}\hspace{5cm}
\end{eqnarray*}

\begin{eqnarray*}
& = & \left((-1)^N\frac{c_+}{c_-}\right)^{\sigma_+^1+\sigma_+^2+g-\sigma_+}\cdot
\frac{d_1^{N}d_2^{N}}{d^{N}}
\cdot\left(\frac{\tilde\iota_{N}(\tilde\iota_{N}\cZ(L_1,G_1,G_1^{\prime})
\ast
(Z(T_{g})\otimes(\nu^{1/2})^{\otimes 2g})
\ast
\tilde\iota_{N}\cZ(L_2,G_2,G_2^{\prime}))}
{(-1)^{\sigma_+^1N}(-1)^{\sigma_+^2N}(-1)^{\sigma_+N} 
d_1^{N}d_2^{N}\cdot c_+^{\sigma_+^1} c_-^{\sigma_-^1}
\cdot c_+^{g} c_-^{g}
\cdot c_+^{\sigma_+^2} c_-^{\sigma_-^2}}\right)^{[\leq N]} \\
& = & \left((-1)^N\frac{c_+}{c_-}\right)^{\sigma_+^1+\sigma_+^2+g-\sigma_+}\cdot
\frac{d_1^{N}d_2^{N}}{d^{N}}
\cdot\left(\tilde\iota_{N}\left(
\frac{\tilde\iota_{N}\cZ(L_1,G_1,G_1^{\prime})}
{(-1)^{\sigma_+^1N}d_1^{N}c_+^{\sigma_+^1} c_-^{\sigma_-^1}}
\ast
\frac{Z(T_{g})\otimes(\nu^{1/2})^{\otimes 2g}}{(-1)^{gN}c_+^{g} c_-^{g}}
\ast
\frac{\tilde\iota_{N}\cZ(L_2,G_2,G_2^{\prime})}
{(-1)^{\sigma_+^2N}d_2^{N}c_+^{\sigma_+^2} c_-^{\sigma_-^2}}
\right)\right)^{[\leq N]} \\
& = & \left((-1)^N\left(\frac{c_+}{c_-}\right)^{[\leq N]}\right)
^{\sigma_+^1+\sigma_+^2+g-\sigma_+}\cdot
\left(\frac{d_1^Nd_2^N}{d^N}\right)\cdot
\ell^{\leq N}(\tau^{\leq N}(M_1),\tau^{\leq N}(M_2))
\end{eqnarray*}

\noindent
where we have used that $\sigma_+ + \sigma_- = 
\sigma_1^+ + \sigma_1^- + \sigma_2^+ + \sigma_2^- +2\cdot g$.
Observe that in the second equality, when "braking" $\cZ$ into three, on each component of $L_0$
a $\nu^{1/2}$ "goes" to $Z$ of $G_1^{\prime}$ or $G_2$, and another $\nu^{1/2}$ goes to
$z_g$. In fact, the two middle expressions are written for the even associator. For any
other associator we would insert between the $\ast$'s the element $A$ mentioned in the
remark at the end of 2.7.
}

\medskip
Let $(L,G,G^{\prime})$ be a triplet and $(M,f,f^{\prime})=\kappa(L,G,G^{\prime})$. 
We can talk about linking number between a link component $K$ 
and a circle $U$ of a chain graph, as well as between two circles $U$ and $V$ of chain graphs:
$lk(K,U)=lk(U,K)$ is defined to be the linking number between $K$ and the knot obtained 
from the graph by deleting all but the circle component $U$, and similarly for $lk(U,V)$.
{\it The linking matrix of a triplet} is then:
\begin{equation}\label{linkingMatrixTriplet}
lk(L,G,G^{\prime})=\left(\begin{array}{ccc}
lk(L) & lk(L,G) & lk(L,G^{\prime}) \\
lk(G,L) & lk(G,G) & lk(G,G^{\prime}) \\
lk(G^{\prime},L) & lk(G^{\prime},G) & lk(G^{\prime},G^{\prime})
\end{array}\right)=\left(\begin{array}{ccc}
A & B^T & C^T \\
B & D & E^T \\
C & E & F
\end{array}\right)
\end{equation}

\noindent
where $A, D, F$ are symmetric matrices. In [\ref{CL}] it has been shown that the
semi-Lagrangian condition can be expressed:
\begin{eqnarray}\label{lagrange2}
D & = & BA^{-1}B^T\nonumber\\
F & = & CA^{-1}C^T
\end{eqnarray}

\noindent
(for $\qq$-cobordisms this in particular means that the entries on the left-hand side, a priori in
$\zz\left[\frac{1}{\det A}\right]$, must be in $\zz$), and for the case $\ff=\qq$ additionally
$BA^{-1}C^T \in  {\cal M}_{g_1\times g_2}(\zz)$.
We will need the following elementary

\medskip\noindent
{\it Ramark.} The signature of a symmetric $2g\times 2g$-matrix 
$\left(\begin{array}{cc}A & -I \\ -I & {\bf 0}\end{array}\right)$ with
integer, respectively real entries is $(g,g)$. 
The determinant of such a matrix is $(-1)^g$.

\bigskip
\demlung{Proof of Proposition 2.11}{ Note that $w_g=
\tau(\Sigma_g\times[0,1],(\Sigma_g\times0,p_1),
(\Sigma_g\times1,p_2))$.
Using the gluing formula (\ref{anomaly1}), for any $\qq$HC $(M,f_1,f_2)$,
$\tau^{\leq N}((\Sigma_g\times[0,1])\cup_{p_1\circ(f_2)^{-1}}M,f_1,p_2)=
\left((-1)^N\left(\frac{c_+}{c_-}\right)^{[\leq N]}\right)^{\sigma_+^1+\sigma_+^2+g-\sigma_+}
\cdot\left(\frac{d_1d_2}{d}\right)^N\cdot
\ell^{\leq N}(\tau^{\leq N}(M,f_1,f_2),w_g)$.
If $\widehat{M}=S^3_L$, and the linking matrix of $L$ is $lk(L)$,
then the linking matrix of the link $L\cup L_0$ is 
$\left(\begin{array}{ccc}lk(L)&\ast &{\bf 0}\\ \ast &\ast &-I\\ {\bf 0}&-I&{\bf 0}\end{array}\right)
\sim\left(\begin{array}{ccc}lk(L)&{\bf 0}&{\bf 0}\\{\bf 0}&\ast&-I\\ {\bf 0}&-I&{\bf 0}\end{array}\right)$.
Using the above remark, $\sigma_+=\sigma_+^1+g$, $\sigma_-=\sigma_-^1+g$, 
$\sigma_+^2=\sigma_-^2=0$, $d_2=lk(\emptyset)=1$, $d_1=d$.
Observe that $((\Sigma_g\times[0,1])\cup_{p_1\circ(f_2)^{-1}}M,f_1,p_2)\cong(M,f_1,f_2)$.
Hence $\ell^{\leq N}(\tau^{\leq N}(M,f_1,f_2),w_g)=\tau^{\leq N}(M,f_1,f_2)$.
In particular, this holds if $(M,f_1,f_2)$ is a $\zz$HC with bottom $S^2$, and hence also
for any $(M,\emptyset,f)$ such that $\widehat{M}$ is a $\zz$-homology sphere. 
The statement now follows from Theorem 3.2 and Lemma 3.5.
}

\medskip 
Note that Proposition 2.11 verifies Axiom (A3) for the truncated TQFTs 
$\mathfrak Q\rightarrow{\cal A}^{\leq N}$ and 
$\mathfrak Z\rightarrow{\cal A}^{\leq N}$.


\medskip
\subsec{Absence of anomaly.}In [\ref{CL}] it has been shown that with the above notations
the linking matrix
\begin{equation}\label{lkmatrix}
lk(L_1\cup L_0\cup L_2)=\left(
\begin {array}{cccc}
A & B^T & {\bf 0} & {\bf 0} \\
B & BA^{-1}B^T & -I & {\bf 0} \\
{\bf 0} & -I & DC^{-1}D^T & D \\
{\bf 0} & {\bf 0} & D^T & C
\end{array}
\right)
\end{equation}

\noindent 
where $A=lk(L_1)\in{\cal M}_{|L_1|\times|L_1|}(\zz)$, 
$C=lk(L_2)\in{\cal M}_{|L_2|\times|L_2|}(\zz)$, 
$B=lk(G_1^{\prime},L_1)\in{\cal M}_{g\times|L_1|}(\zz)$, 
$D=lk(G_2,L_2)\in{\cal M}_{g\times|L_2|}(\zz)$, 
$BA^{-1}B^T, DC^{-1}D^T\in{\cal M}_{g\times g}(\zz)$. 
There it has been proven the following

\medskip
\te{Proposition {\rm [\ref{CL}, Proposition 13]}.}{The signature of the matrix (\ref{lkmatrix}) is 
$(\sigma_+^1+\sigma_+^2+g,\sigma_-^1+\sigma_-^2+g)$, where $(\sigma_+^1,\sigma_-^1)$,
respectively $(\sigma_+^2,\sigma_-^2)$ is the signature of $lk(L_1)$, respectively
$lk(L_2)$. Also the following holds:
\begin{equation}\label{detgluing}
\det(lk(L_1\cup L_0\cup L_2)) = (-1)^g\cdot\det(lk(L_1))\cdot\det(lk(L_2))
\end{equation}
}

\noindent
Therefore, we can re-write the {\it gluing formula} (\ref{anomaly1}):
\begin{equation}\label{gluingN}
\tau^{\leq N}(M_2\cup_{f_2\circ(f_1^{\prime})^{-1}}M_1,f_1,f_2^{\prime}) 
=\ell^{\leq N}(\tau^{\leq N}(M_1,f_1,f_1^{\prime}),
\tau^{\leq N}(M_2,f_2,f_2^{\prime}))
\end{equation}


\demlung{Proof of Theorem 2.12}{ 1) By construction, the inverse limits 
$\lim\limits_{\infty\leftarrow N}\aa^{\leq N}(\emptyset)=\aa(\emptyset)$ and 
$\lim\limits_{\infty\leftarrow N}\aa^{\leq N}(\Gamma^{g_1},\Gamma^{g_2})=\aa(\Gamma^{g_1},\Gamma^{g_2})$.
Let us show that the following diagram is commutative for every $N\in\nn$:

\begin{eqnarray}\nonumber\label{ellCommute}
\aa^{\leq N+1}(\Gamma^{g_1},\Gamma^{g_2})\otimes
\aa^{\leq N+1}(\Gamma^{g_2},\Gamma^{g_3}) & \rightarrow &
\aa^{\leq N}(\Gamma^{g_1},\Gamma^{g_2})\otimes
\aa^{\leq N}(\Gamma^{g_2},\Gamma^{g_3}) \\
\downarrow \ell^{\leq N+1} & & \downarrow \ell^{\leq N}\\
\aa^{\leq N+1}(\Gamma^{g_1},\Gamma^{g_3}) & \rightarrow &
\aa^{\leq N}(\Gamma^{g_1},\Gamma^{g_3})\nonumber
\end{eqnarray}

\noindent
where the horizontal arrows are the maps that forget the degrees $N+1$ 
parts. Let $\alpha=\tau^{\leq N+1}(M_1)$,
$\beta=\tau^{\leq N+1}(M_2)$ for some $\qq$HC $M_1$ and $M_2$. Then as previously observed
$\left(\tau^{\leq N+1}(M_i)\right)^{[\leq N]}=\tau^{\leq N}(M_i)$, $i=1,2$,
i.e. $\alpha^{\leq N}=\tau^{\leq N}(M_1)$, $\beta^{\leq N}=\tau^{\leq N}(M_2)$.
By the gluing formula (\ref{gluingN}) we then have
$\tau^{\leq N+1}(M_2\cup M_1) = \ell^{\leq N+1}(\alpha,\beta)$ and
$\tau^{\leq N}(M_2\cup M_1) = \ell^{\leq N}(\alpha^{\leq N},\beta^{\leq N})$.
Again, using now
$\left(\tau^{\leq N+1}(M_2\cup M_1)\right)^{[\leq N]}=\tau^{\leq N}(M_2\cup M_1)$,
we get 
$\left(\ell^{\leq N+1}(\alpha,\beta)\right)^{\leq N}=
\ell^{\leq N}(\alpha^{\leq N},\beta^{\leq N})$.
Hence the diagram (\ref{ellCommute}) is commutative for $\alpha, \beta$ as above.
By the remark after the proof of Theorem 3.2, and by Lemma 3.5, the diagram is then
commutative for arbitrary $\alpha, \beta$.

Therefore \hfill there \hfill exists \hfill a \hfill well-defined \hfill $\aa(\emptyset)$-bilinear \hfill map\hfill 
$\ell:\aa(\Gamma^{g_1},\Gamma^{g_2})\otimes\aa(\Gamma^{g_2},\Gamma^{g_3})\rightarrow$
 
\noindent
$\aa(\Gamma^{g_1},\Gamma^{g_2})\widetilde\otimes\aa(\Gamma^{g_2},\Gamma^{g_3})
\rightarrow\aa(\Gamma^{g_1},\Gamma^{g_3})$,
such that when restricting to the internal degree $\leq N$ parts one obtains the map $\ell^{\leq N}$.

2) By the proof of 1), 
$\tilde\ell(w_g)= \lim\limits_{\infty\leftarrow N}\tilde\ell^{\leq N}(w_g^{N})$.
By Proposition 2.11 the operators $\tilde\ell^{\leq N}(w_g^{N})$ are identities, 
hence so is the limit.
}


\medskip
\subsec{The functors $\mathfrak Q\rightarrow\cal A$ and $\mathfrak Z\rightarrow\cal A$.}Theorem 2.12.1) 
shows that $\ell^{\leq N}$ are the ${}^{[\leq N]}$-truncations of $\ell$. 
With (\ref{gluing}), this implies Axiom (A2) for the non-truncated TQFTs ${\mathfrak Z}\rightarrow\aa$ 
and ${\mathfrak Q}\rightarrow\aa$:
\begin{equation}\label{gluing}
\tau(M_2\cup_{f_2\circ(f_1^{\prime})^{-1}}M_1,f_1,f_2^{\prime}) 
=\ell(\tau(M_1,f_1,f_1^{\prime}),
\tau(M_2,f_2,f_2^{\prime}))
\end{equation}
\noindent
Theorem 2.12.2) in particular implies Axiom (A3) for the non-truncated TQFTs.


\medskip
\subsec{Lemma.}{\it For any $\beta\in\aa(\Gamma^{g_2},\Gamma^{g_3})$ and
any sequence of elements $\alpha^n\in\aa(\Gamma^{g_1},\Gamma^{g_2})$, such that for every $n$ 
$(\alpha^n)_{\leq n}=(\alpha^{n+1})_{\leq n}=\dots$, both sides of the following equalities are
well-defined and the equalities holds:

\begin{equation}\label{passingToLimit}
\ell(\lim\limits_n\alpha^n,\beta) = \lim\limits_n\ell(\alpha^n,\beta)
\end{equation}

\noindent
A similar property holds for the r\^ ole of two arguments of $\ell$ reversed.
}

\dem{$\ell^{\leq N}$ are the ${}^{[\leq N]}$-truncations of $\ell$. 
Apply (\ref{passingToLimitN}) and pass to the limit (keeping, for example $n=(2g+1)N$).
}

\medskip
Theorem 3.2 shows that our TQFTs are non-degenerate, and Lemmas 3.5 and 3.9 show that 
$\ell^{\leq N}$ and $\ell$ are continuous maps.


\medskip
\subsec{Conjugation.}The 
operation of {\em conjugation} in $\aa(\emptyset)$ can be extended as 
follows. Grade the modules $\aa(\emptyset)$, $\aa(\Gamma^g)$ and 
$\aa(\Gamma^{g_1},\Gamma^{g_2})$ by the internal degree, and define for an 
arbitrary chord diagram $D$, and an arbitrary natural number $k$
\begin{eqnarray*}
\overline{D^{[2k]}} & = & D^{[2k]}\\
\overline{D^{[2k+1]}} & = & -D^{[2k+1]}
\end{eqnarray*}
 
\noindent
Note that $\T(\Sigma_g)=\T(-\Sigma_g)=\aa(\Gamma^g)$, 
and $\T(f)=id$, $\forall f\in Homeo(\Sigma_g)$,
hence the naturality of this antimorphism is obvious.
$\overline{\cdot}:\aa(\Gamma^{g_1},\Gamma^{g_2})
\rightarrow\aa(\Gamma^{g_2},\Gamma^{g_1})$ satisfies the requirement of axiom (A4)
to commute with homeomorphisms of 3-cobordisms,
because already $\tau$ is defined for homeomorphism classes.
The same can be repeated with added $\leq N$.

Also $\tilde\iota_N\cZ\overline{(L,G)}=(-1)^{|L|N}\overline{\tilde\iota_N\cZ(L,G)}$
(compare with [\ref{LMO}, Proposition 5.2]).
This is true for the Murakami-Ohtsuki extension of $Z$ because 
$a$ and $b$ from [\ref{MO}], and hence $Z($vicinity of a trivalent vertex$)$ 
are ``mirrors'' of themselves, which is easy to check.
For the extension of $Z$ from 2.6 this property is obvious.
From the proof of Proposition 5.2 in [\ref{LMO}] it also follows 
that $c_-=\overline{c_+}$. Hence for any $N$:
\begin{equation}\label{conjugation}
\tau^{\leq N}(-M)
=
\frac{(-1)^{\sigma^{\prime}_+N}}{d^N}
\left(\frac{\tilde\iota_N(\cZ\overline{(L,G)})}
{(c_+^{[\leq N]})^{\sigma_+^{\prime}}\cdot
(c_-^{[\leq N]})^{\sigma_-^{\prime}}}\right)^{[\leq N]}
=
\frac{(-1)^{\sigma_+N}}{d^N}
\left(\frac{\overline{\tilde\iota_N(\cZ(L,G))}}
{\overline{(c_+^{[\leq N]})}^{\sigma_+}\cdot
\overline{(c_-^{[\leq N]})}^{\sigma_-}}\right)^{[\leq N]} 
= 
\overline{\tau^{\leq N}(M)}
\end{equation}

\noindent
Therefore also $\tau(-M)=\overline{\tau(M)}$.
Using this formula and (\ref{gluing}), it follows that for $\alpha=\tau(M_1)$, $\beta=\tau(M_2)$,
where $M_i$ are 3-cobordisms in our category, we have
$\overline{\ell(\alpha,\beta)}=\overline{\ell(\tau(M_1),\tau(M_2))}
=\overline{\tau(M_2\cup M_1)}=\tau(-(M_2\cup M_1))
=\tau((-M_1)\cup(-M_2))=\ell(\tau(-M_2),\tau(-M_1))
=\ell(\overline{\tau(M_2)},\overline{\tau(M_1)})=\ell(\overline{\beta},\overline{\alpha})$.
By Theorem 3.2 and (\ref{passingToLimit}), the same relation holds for arbitrary $\alpha, \beta$.
In particular, it remains true with $\leq N$ added.
Axiom (A4) is therefore verified for truncated and non-truncated TQFTs.

%


\medskip
\subsec{Conclusions and consequences.}The 
full and truncated TQFTs are now completely constructed. 
The full TQFT induces a linear representation 
${\cal L}_g\rightarrow GL_{\aa(\emptyset)}(\aa(\Gamma^g))$.
The truncated TQFTs induce linear representations
${\cal L}_g\rightarrow GL_{\aa^{\leq N}(\emptyset)}(\aa^{\leq N}(\Gamma^g))$.
It is known [\ref{FM}] that any $\zz$HS can be obtained as filling of a parametrized 
3-cobordism $(\Sigma_g\times I,w,id)$ for some $g\geq 0$ and some $w\in{\cal T}_g$, 
the Torelli group of genus $g$.
Furthermore [\ref{morita}] it even suffices to consider only $w\in{\cal K}_g$, 
the kernel of the Johnson homomorphism, or topologically the subgroup of 
${\cal T}_g$ generated by Dehn twists on bounding simple closed curves.
Our TQFTs, of cause, induce linear representations of both these subgroups of ${\cal L}_g$.
The group ${\cal L}_g$ has not been studied before, no explicit set of generators,
less so one of relations, is known.

Note, that theorem 3.2 and Lemmas 3.5 and 3.9 not only allow a well-defined non-truncated
TQFT (Theorem 2.12) and prove the non-degeneracy, but also solve the realization problem for links,
string links, three-dimensional manifolds and chain graphs, by showing (see Lemmas 3.3 and 3.4)
that $Z($links$)$, $\tau($closed 3-manifolds$)$, $Z($string links$)$ and $\tau($3-manifolds with
boundary$)$ in the closure generate the corresponding spaces of chord diagrams:
$\aa(\oOriented\dots\oOriented)$, $\aa(\emptyset)$ and $\aa(\arrowg)$.
(For links, the correspondent for Habegger-Masbaum result follows easily from Habiro's calculus of 
claspers [\ref{habiro}].) Without proving Theorem 3.2 even partial results of this sort were hard to obtain, 
as we can exemplify by the following

\te{Proposition.}{For every $N\geq 0$ and every $\qq$-homology 
handlebody $(M,f_1,f_2)$, $\tilde\ell^{\leq N}(\tau^{\leq N}(M,f_1,f_2))$
sends the $\aa^{\leq N}(\emptyset)$-submodule of $\aa^{\leq N}(\arrowg)$ 
generated by $\exp(\alpha)$, $\alpha\in\a(\arrowg)^{\leq N}$ to itself. 
}

\dem{By Proposition 3.1.3) $\tau^{\leq N}(M,f_1,f_2)$ is group-like. 
Observe that $\Delta$ commutes with $\ast$, and repeating the argument 
from the proof of 3.1.3) for $\tilde\iota_N$ in the definition of 
$\ell^{\leq N}$, we can see that $\tilde\ell^{\leq N}(\tau^{\leq N}(M,f_1,f_2))$ 
takes a group-like element of $\aa^{\leq N}(\Gamma^{g_1})$ of the form $1+h.o.t.$ 
to a group-like element of $\aa^{\leq N}(\Gamma^{g_2})$ of the form $1+h.o.t.$ 
Now apply Proposition 2.4.5) and 3) for the truncated case. 
Hence it sends the $\aa^{\leq N}(\emptyset)$-submodule 
of $\aa^{\leq N}(\arrowg)$ generated by $\exp(\alpha)$, 
$\alpha\in\a(\arrowg)^{\leq N}$ to itself. 
}

\medskip\noindent
{\it Remark.} This construction of TQFT can be done also in the language of the Aarhus integral [\ref{CM}].

\begin{figure}[thb]
\centerline{\psfig{file=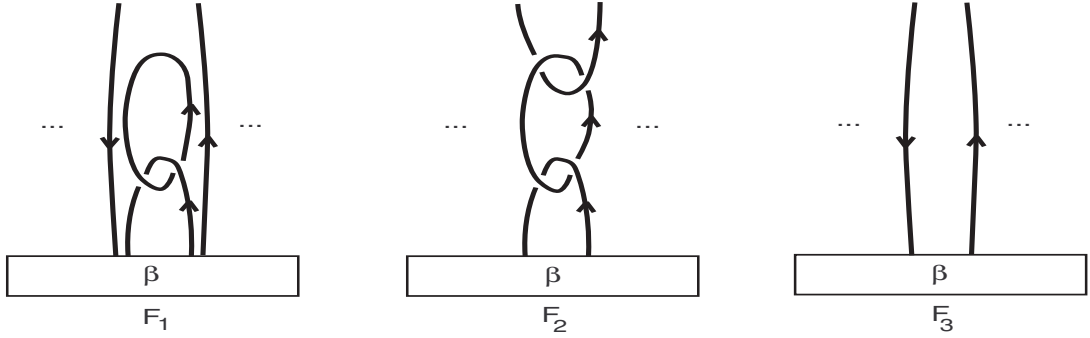,width=15cm,height=5cm,angle=0}}
\caption{\sl Chord-handle-canceling relation $F_1=F_3$}\label{figChordCansel}
\end{figure}

\medskip
\subsec{Chord-handle canceling.}In 
[\ref{LMO}] Le, Murakami and Ohtsuki have introduced the chord-KII move
to mirror the second Kirby move for links, which then allowed them to define
$Z^{LMO}$. However, it is well-known that handle canceling can not be obtained solely
by Kirby-2, and would require in addition Kirby-1. But no corresponding chord-KI move exists,
the invariance of $Z^{LMO}$ under Kirby-1 is achieved via normalization.
Therefore there is no a priori reason to suspect that a chord-canceling-handle relation
is true for arbitrary chord diagrams. But, the result obtained here above allow us to prove:

\te{Proposition.}{The chord-handle-canceling relation, schematically depicted in figure 
\ref{figChordCansel} holds for arbitrary $\beta\in\aa(\arrowg)$. (The upper part of each 
$F_i$ should be read as $Z($drawn tangle$)$.)}

\dem{For arbitrary $\beta$, $F_1$ differs from $F_2$ by a chord-KII move. 
(An argument similar to the one in [\ref{LMO}, Proposition 3.2] works.)
But now $F_2=\ell(\beta,w_g)=\beta=F_3$.
}

%
%

\section{A TQFT for the Casson-Walker-Lescop invariant}
\setcounter{equation}{0}
\setcounter{nsubsec}{1}



The term of degree one of $Z^{LMO}$ of a 3-manifold is
$(-1)^{b_1(M)}\frac{\lambda(M)}{2}\theta$, where $b_1(M)$ is the first Betti number,
$\lambda(M)$ is the Casson invariant (in Walker-Lescop extension) and $\theta$ is 
the (only) open chord diagram of degree 1, which looks like $\theta$ [\ref{LMO}].
Let us recall the definition and basic properties of Casson invariant. Let $K$ be
a knot in an oriented $\zz$-homology 3-sphere $M$, and 
$\Delta_K(t)=a_0+a_1(t+t^{-1})+a_2(t^2+t^{-2})+...$ be its Alexander polynomial
normalized such that $\Delta_K(1)=1$. Denote 
$\lambda^{\prime}(K)=\frac{1}{2}\Delta^{\prime\prime}(1)=\sum\limits_n n^2a_n$.


\medskip
\subsec{Theorem (Casson).}{\it There is an integer-valued invariant $\lambda$ for oriented
integer homology 3-spheres such that:

\rm(1)\it \quad $\lambda$ mod $2$ is the Rohlin invariant

\rm(2)\it \quad $\lambda(M)=0$ for any homotopy 3-sphere

\rm(3)\it \quad $\lambda(-M)=-\lambda(M)$

\rm(4)\it \quad $\lambda(M_1\# M_2)=\lambda(M_1)+\lambda(M_2)$

\rm(5)\it \quad If $K$ is a knot in an oriented integer homology 3-sphere $M$, and 
$M(K,\frac{1}{n})$ denotes the integer 

\quad\quad homology 3-sphere obtained from $M$ by a $\frac{1}{n}$-surgery
on $K$, then $\lambda(M(K,\frac{1}{n}))=\lambda(M)+n\lambda^{\prime}(K)$.}

\medskip
Property (5) from this theorem for $n=\pm 1$ with $\lambda(S^3)=0$ determine $\lambda$ uniquely,
since any integer homology 3-sphere can be obtained from $S^3$ by a succession of $\pm 1$-surgeries
on knots. $\lambda$ was extended to rational homology 3-spheres by Walker, and corresponding
properties (4) and (5) were given by Lescop[\ref{lescop}]:

\medskip
\rm (4$^{\prime}$)\it \quad $\lambda(M_1\# M_2)=|H_1(M_2,\zz)|\lambda(M_1)+|H_1(M_1,\zz)|\lambda(M_2)$

\rm (5$^{\prime}$)\it \quad $\lambda(M(L,\frac{p_1}{q_1},\dots,\frac{p_{|L|}}{q_{|L|}}))
=
\frac{|H_1(M(L,\frac{p_1}{q_1},\dots,\frac{p_{|L|}}{q_{|L|}}),\zz)|}{|H_1(M,\zz)|}\lambda(M)
+
{\cal F}_M(L,\frac{p_1}{q_1},\dots,\frac{p_{|L|}}{q_{|L|}})$,

\noindent\rm
where $M(L,\frac{p_1}{q_1},\dots,\frac{p_{|L|}}{q_{|L|}})$ is the manifold obtained from $M$
by performing rational surgery with indicated coefficients on the components of the link $L$,
and ${\cal F}_M(L,\frac{p_1}{q_1},\dots,\frac{p_{|L|}}{q_{|L|}})$ is a certain function on
the set of surgery presentations in $M$, which is essentially a function of the linking matrix,
homology and Alexander polynomial [\ref{lescop}].


\medskip
\subsec{The degree $^{N\leq 1}$ truncation} of our TQFT 
can be thought of as a TQFT for the Casson-Walker-Lescop invariant.
Note that the ring $\aa^{\leq 1}(\emptyset)=(\{r+s\theta|r,s\in\qq\}, [\leq 1]$-multiplication$)\cong
\qq[\theta]/(\theta^2)\cong\left\{\left(\begin{array}{cc}r&s\\ 0&r\end{array}\right)|r,s\in\qq\right\}$
Denote it by $R$. 
Observe [\ref{LMO}] that $c_+=1-\frac{\theta}{16}+h.o.t.$, $c_-=1+\frac{\theta}{16}+h.o.t.$, hence
$c_+c_-=1+$terms of degree $\geq 2$, and therefore 
$z_g^1=(-1)^{g\cdot 1}\left(\frac{Z(T_g)\otimes(\nu^{1/2})^{\otimes 2g}}{c_+c_-}\right)^{[\leq 1]}=(-1)^gZ(T_g)\otimes(\nu^{1/2})^{\otimes 2g}$.
If $\alpha\in\aa^{\leq 1}(\Gamma^{g_1},\Gamma^{g_2})$ and $\beta\in\aa^{\leq 1}(\Gamma^{g_2},\Gamma^{g_23})$,
then:
\begin{equation}\label{ellCasson}
\ell_{Casson}(\alpha,\beta)=\left(\widetilde{\iota_1}(\alpha\ast z_{g_1}^1\ast\beta)\right)^{[\leq 1]}
\end{equation}

\noindent
If $g=0$, $\ell_{Casson}$ is the disjoint union (multiplication). 
A formula for $Z(T_1)$ and $Z(W_1)$ is given for example in [\ref{gille}]. Using the even associator it is easy 
to write down $z_1^1$ and $w_1$ explicitly, at least the degree $\leq 3$ terms.

Let $\kappa(L,G_1,G_2)=(M,f_1,f_2)$. Then, keeping in mind that $c_+c_-=1+$terms of degree $\geq 2$, denote
\begin{equation}\label{tauCasson}
c(M,f_1,f_2):=\tau^{\leq 1}(M,f_1,f_2)=\frac{(-1)^{\sigma_+}}{d(M)}
\widetilde{\iota_1}(\cZ(L,G_1,G_2) )^{[\leq 1]}
\end{equation}

\noindent
where $d(M)=|H_1(\widehat{M},\zz)|$ and $\sigma_+=sign_+(L)$, is an invariant of 3-cobordisms 
of category $\mathfrak Q$. In particular, for cobordisms between $S^2$ and $S^2$, 
$c(M,id_{S^2},id_{S^2})=\frac{(-1)^{\sigma_+}}{d(M)}\widetilde{\iota_1}(\cZ(L))^{[\leq 1]}=Z^{LMO}(M)^{[\leq 1]}
=1+\frac{\lambda(M)}{2}\theta$, where we have identified 
$\aa^{\leq 1}(\Gamma^{0},\Gamma^{0})\equiv R=\aa^{\leq 1}(\emptyset)$. 
The filling of the composition of two 3-cobordisms between $S^2$ and $S^2$
is clearly the connected sum of the fillings. Hence $c(M_2\cup M_1,S^2,S^2)=c(M_1,S^2,S^2)c(M_2,S^2,S^2)$
implies, as it is easy to check, property (4) of the Casson invariant (the generalized version for $\qq$HS). 
As we have shown the following axioms of TQFT hold:
\begin{eqnarray}
c(M_2\cup_{f_2\circ (f_1^{\prime})^{-1}} M_1, f_1, f_2^{\prime}) 
& = & \ell_{Casson}(c(M_1,f_1,f_1^{\prime}),c(M_2,f_2,f_2^{\prime}))\label{ellCassonGlue}\\ 
c(\Sigma_g\times[0,1],p_1,p_2) & = & id_{\aa^{\leq 1}(\Gamma^g)}\\
c(-M,-f_2,-f_1) & = & \overline{c(M,f_1,f_2)}\label{conjugationCasson}
\end{eqnarray}

\noindent
where the notations are obvious. $R$, $\aa^{\leq 1}(\Gamma^0,\Gamma^g)$ and $\aa^{\leq 1}(\Gamma^{g_1},\Gamma^{g_2})$
are $\zz_2$-graded by the internal degree; the conjugation changes the sign of the internal degree 1 
part. In particular (\ref{conjugationCasson}) implies property (3) of the CWL invariant. 
It is natural to try now to obtain property (5) of the CWL invariant
as a consequence of the rational surgery formula of Bar-Natan and Lawrence.

Unfortunately, explicit calculations for $c(M,f_1,f_2)$, as expected, are rather hard to do. 
Now we would like to show that the induced representation 
${\cal L}_g\rightarrow GL_R(\aa^{\leq 1}(\Gamma^g))$
descends to Morita's homomorphism
$\lambda^{\ast}:{\cal K}_g\rightarrow\zz$.
($\lambda^{\ast}$ extends to ${\cal L}_g$, but fails to be a homomorphism there.)


\medskip
\subsec{Proposition.}{\it 1) Let $\mathfrak B$ be the completion of the $\qq$-vector subspace of 
$\aa(\Gamma^{g_1}\sqcup(\sqcup_mS^1)\sqcup\Gamma^{g_2})$ generated by finite sums of chord diagrams which intersect 
$\Gamma^{g_1}\sqcup\Gamma^{g_2}$. Then
$p:\aa(\Gamma^{g_1}\sqcup(\sqcup_mS^1)\sqcup\Gamma^{g_2})\rightarrow\aa(\sqcup_mS^1)$,
the natural map "erase $\Gamma^{g_1}$ and $\Gamma^{g_2}$ from a chord diagram",
if it does not intersect $\Gamma^{g_1}\sqcup\Gamma^{g_2}$, and set $=0$, otherwise,
is well-defined and the following sequence is short exact:
 
$$0\rightarrow\mathfrak B\rightarrow\aa(\Gamma^{g_1}\sqcup(\sqcup_mS^1)\sqcup\Gamma^{g_2})
\stackrel{p}{\rightarrow}\aa(\sqcup_mS^1)\rightarrow0$$

\noindent
We will denote also by $p$ the induced maps on minimal internal degree $^{\leq N}$ parts. They have the same property.

2) Denote by $r$ the maps similar to $p$ from 1) corresponding to the case $m=0$. Then the following diagram is commutative:

$$\aa(\Gamma^{g_1}\sqcup(\sqcup_mS^1)\sqcup\Gamma^{g_2})
\stackrel{\tilde{\iota}_N}{\longrightarrow}\aa(\Gamma^{g_1}\sqcup\Gamma^{g_2})
\longrightarrow\aa^{\leq N}(\Gamma^{g_1}\sqcup\Gamma^{g_2})$$
$$\downarrow p\quad\quad\quad\quad\quad\quad\quad\quad\downarrow r\quad\quad\quad\quad\quad\quad\downarrow r\quad\quad$$
$$\aa(\sqcup_mS^1)\quad\quad\stackrel{\iota_N}{\longrightarrow}\quad\quad
\aa(\emptyset)\quad\quad\longrightarrow\quad\aa^{\leq N}(\emptyset)\quad\quad\quad$$

3) For every embedding $L\cup G\hookrightarrow S^3$, such that $G$ as an abstract graph is
$\Gamma^{g_1}\sqcup\Gamma^{g_2}$, $p\cZ(L\cup G)=\cZ(L)$.

4) For every embedding $L\cup G\hookrightarrow S^3$, such that $G$ as an abstract graph is
$\Gamma^{g_1}\sqcup\Gamma^{g_2}$, and every $N\geq 1$, 
$p(\tau^{\leq N}(\kappa(L\cup G)))=Z^{LMO}(\widehat{\kappa(L\cup G)})^{[\leq N]}$.
In particular (if $N=1$), $p(c(M,f_1,f_2))=1+\frac{\lambda(\widehat{M})}{2}\theta$.

5) If $\pfi_1,\pfi_2\in{\cal K}_g$,
then $p(c(\Sigma_g\times I,\pfi_2\circ\pfi_1,id))
=p(c(\Sigma_g\times I,\pfi_1,id)) p(c(\Sigma_g\times I,\pfi_2,id))$. 
}

\dem{1) The following argument can be worked for every fixed degree, and since all relations are homogeneous,
we can use the universality property of the direct product as mentioned in 2.2 to obtain the
desired statement. Consider the corresponding diagram before introducing relations:
$$0\rightarrow\mathfrak B^{\prime}\rightarrow{\cal D}(\Gamma^{g_1}\sqcup(\sqcup_mS^1)\sqcup\Gamma^{g_2})
\stackrel{p^{\prime}}{\rightarrow}{\cal D}(\sqcup_mS^1)\rightarrow0$$

\noindent
The terms of any relation for diagrams on $\Gamma^{g_1}\sqcup(\sqcup_mS^1)\sqcup\Gamma^{g_2}$, either
all intersect $\Gamma^{g_1}\sqcup\Gamma^{g_2}$, or none does. Hence, if we denote by $R_1$ the $\qq$-vector space
generated by relations of the first type, by $R_2$ - the space generated by relation of the second type, 
and by $R$ - the one generated by all relations, then $R/R_1\cong R_2$.
All in all we get a diagram:

$$0\rightarrow R_1\quad\quad\longrightarrow\quad\quad\quad R\quad\quad
\longrightarrow\quad\quad R_2\quad\longrightarrow 0$$
$$\downarrow\quad\quad\quad\quad\quad\quad\quad\quad\downarrow\quad\quad\quad\quad\quad\quad\downarrow\quad\quad$$
$$0\rightarrow\mathfrak B^{\prime}\rightarrow{\cal D}(\Gamma^{g_1}\sqcup(\sqcup_mS^1)\sqcup\Gamma^{g_2})
\stackrel{p^{\prime}}{\rightarrow}{\cal D}(\sqcup_mS^1)\rightarrow0$$
$$\downarrow\quad\quad\quad\quad\quad\quad\quad\quad\downarrow\quad\quad\quad\quad\quad\quad\downarrow\quad\quad$$
$$0\rightarrow\mathfrak B\stackrel{i}{\rightarrow}\aa(\Gamma^{g_1}\sqcup(\sqcup_mS^1)\sqcup\Gamma^{g_2})
\stackrel{p}{\rightarrow}\aa(\sqcup_mS^1)\rightarrow0$$

\noindent
where all columns and the first two rows are short exact. The arrows $i$ and $p$ in the third row are then induced
and make the diagram commutative. They clearly are the maps described in the statement.
The exactness in the third row follows from the exactness in the second.

2) Let $\alpha\in\aa(\Gamma^{g_1}\sqcup(\sqcup_mS^1)\sqcup\Gamma^{g_2})$ and $\beta$ be such that
$\tilde{\pfi}(\beta)=\alpha$. (Recall that $\tilde{\iota_N}=\tilde{q_N}\circ\tilde{\kappa_N}\circ\tilde{\pfi}^{-1}$.)
A chord diagram $x$ from the expression of $\beta$ connects to $\Gamma^{g_1}\sqcup\Gamma^{g_2}$ if and only if
its image via $\pfi$ is a sum $y$ of chord diagrams expressing $\alpha$, all connected to
$\Gamma^{g_1}\sqcup\Gamma^{g_2}$. Again using the fact that the terms in any relation either all
connect or all do not, $p(y)=0$ implies that (in fact, if and only if) $\tilde{q_N}\circ\tilde{\kappa_N}(x)$ 
connects to $\Gamma^{g_1}\sqcup\Gamma^{g_2}$, i.e. 
$r(\tilde{q_N}\circ\tilde{\kappa_N}\circ\tilde{\pfi}^{-1}(y))=r(\tilde{q_N}\circ\tilde{\kappa_N}(x))=0$.

Now, if we decompose $\beta=\beta_1+\beta_2$ such that all terms in $\beta_1$ connect to
$\Gamma^{g_1}\sqcup\Gamma^{g_2}$ and all terms in $\beta_2$ do not, the result follows:
$(r\circ\tilde{\iota_N})(\alpha)=r(\pfi(\beta_1)+\pfi(\beta_2))
=r(\tilde{q_N}\circ\tilde{\kappa_N}(x_1)+\tilde{q_N}\circ\tilde{\kappa_N}(x_2))
=\tilde{q_N}\circ\tilde{\kappa_N}(x_2)=\tilde{\iota_N}(\pfi(\beta_2))
=\tilde{\iota_N}(\pfi(p(\beta)))=\tilde{\iota_N}(p(\pfi(\beta)))=(\tilde{\iota_N}\circ p)(\alpha)$.

3) Decompose $L\cup G$ into elementary pseudo-quasi-tangles. Observe that for everyone, 
except $\tangleThree$'s and $\tangleFour$'s (possibly with multiple strands),
$Z$ either returns diagrams either all in ${\mathfrak B}$, or all having no intersection between the dashed graph and
$\Gamma^{g_1}\sqcup\Gamma^{g_2}$.  Thus, suppressing $G$ for these elementary tangles corresponds precisely to
applying $p$.

The remaining cases. Observe, first, that one can "lift $L$ above $G$", leaving only some "fingers" from $L$
attached to $G$. To see this, from a generic plane projection of $L\cup G$ on $\rr^2\subset\rr^3$ obtain 
an isotopic embedding of $G\cup L$ in $\rr^3$, such that $G$ is in an $\varepsilon$-neighbourhood of the
plane $\{z=0\}\in\rr^3$, and $L$, except for some fingers that correspond to intersections between $G$ and $L$
in the original plane projection, lies in an $\varepsilon$-neighbourhood of the plane $\{z=1\}\in\rr^3$.
Hence, by "opening the two-page book", we can find such a tangle decomposition that all occurring associator-tangles 
are of one of the following three types:

(A) refer only to $G$ or only to $L$;

(B) a single middle strand, which comes from $L$, the left-most strand (with "big" multiplicity)
comes from $G$, the right-most strand (also with "big" multiplicity) comes from $L$.
Moreover, if such an associator-tangle occurs, its inverse (on the same strands) will occur "soon";

(C) one of the left-most two strands is a single strand coming from $L$, all other strands come from $G$.

\smallskip
We will assume that he associator $\Phi$ is horizontal, i.e it is
a formal series in two non-commuting variables $r_{12}, r_{23}$, which correspond to
a dashed line joining these indicated strands [\ref{LMO}].

(A) If all strands are from $G$ or none are from $G$, then all terms of
$\Phi^{\pm 1}=Z($tangle$)$, connect, respectively do not connect to $\Gamma^{g_1}\sqcup\Gamma^{g_2}$.

(C) If two of the three stands come from $G$, $\Phi^{\pm 1}=Z($tangle$)$ will
have all terms connected to $\Gamma^{g_1}\sqcup\Gamma^{g_2}$. Then, eliminating $G$ corresponds precisely 
to replacing this tangle-associator by the single strand from $L$, i.e. corresponds to applying $p$ in this case.

(B) If exactly one (multiple) strand comes from $G$, this corresponds to setting one of the two non-commutative 
variables $r_{12}$, $r_{23}$ zero. But, as we mentioned above, such tangles occur in pairs with their opposite.
Then, both $\Phi$ and $\Phi^{321}=\Phi^{-1}$ occur. Setting one of $r_{12}$, $r_{23}$ zero, still leaves
a series and its inverse (elementary exercise). Thus, eliminating $G$ corresponds again to applying $p$.

4) Recall the definitions of $\tau^{\leq N}$ (\ref{tauN}) and $Z^{LMO}$ (\ref{omega}). Apply $p$ and use the
result of part 3). Then, use the commutativity of the diagram from part 2) to obtain the desired relation.

5) Applying $p$ to (\ref{ellCassonGlue}),
$p(c(\Sigma_g\times I,\pfi_2\circ\pfi_1,id))
=p(\ell_{Casson}(c(\Sigma_g\times I,\pfi_1,id), c(\Sigma_g\times I,\pfi_2,id)))$.
Using part 4), $p(c(\Sigma_g\times I,\pfi_2\circ\pfi_1,id))=1+\frac{\lambda(W_{\pfi_2\circ\pfi_1})}{2}\theta$,
$p(c(\Sigma_g\times I,\pfi_1,id))=1+\frac{\lambda(W_{\pfi_1})}{2}\theta$, and
$p(c(\Sigma_g\times I,\pfi_2,id))=1+\frac{\lambda(W_{\pfi_2})}{2}\theta$,
where $W_{\pfi_i}=\widehat{(\Sigma_g\times I,\pfi_i,id)}$.
But 
$1+\frac{\lambda(W_{\pfi_2\circ\pfi_1})}{2}\theta
=\left(1+\frac{\lambda(W_{\pfi_1})}{2}\theta\right)
\left(1+\frac{\lambda(W_{\pfi_2})}{2}\theta\right)$
in $R$, because $\lambda^{\ast}:{\cal K}_g\rightarrow\zz$,
$\lambda^{\ast}(\pfi):=\lambda(W_{\pfi})$ satisfies
$\lambda^{\ast}(\pfi_2\circ\pfi_1)=\lambda(\pfi_1)+\lambda(\pfi_2)$
by [\ref{morita}].
}

\medskip\noindent
{\em Remark.}
Using expression (\ref{ellCasson}) for $\ell$  and observing that $p$ commutes with $\tilde{\iota_1}$ and 
with taking $^{[\leq 1]}$ by Proposition 4.3.2), we can re-write 
$p(\ell_{Casson}(c(\Sigma_g\times I,\pfi_1,id), c(\Sigma_g\times I,\pfi_2,id)))$ as
$\tilde{\iota_1}(p(c(\Sigma_g\times I,\pfi_1,id)\ast z_g^1\ast c(\Sigma_g\times I,\pfi_2,id)))^{[\leq 1]}$.
Expressing $c(\Sigma_g\times I,\pfi_i,id)$ by (\ref{tauCasson}), and keeping in mind 
the definition of $p$ and  properties of $\tilde{\iota_1}$, 
we can get to having to apply $p$ on $\left(\cZ(L_1,G_1,G_1^{\prime})\ast z_g^1\ast\cZ(L_2,G_2,G_2^{\prime})\right)$,
respectively to apply $p$ on $\left(\cZ(L_i,G_i,G_i^{\prime})\right)$, $i=1,2$.
On the other hand, it is possible to show directly that 
$\tilde{\iota_1}p\left(\cZ(L_1,G_1,G_1^{\prime})\ast z_g^1\ast\cZ(L_2,G_2,G_2^{\prime})\right)
=p\left(\cZ(L_1,G_1,G_1^{\prime})\right)\cdot p\left(\cZ(L_2,G_2,G_2^{\prime})\right)$, 
for suitably chosen $L_i$ in the triplets. This gives another proof of Proposition 4.3.5). 
We, thus, can obtain a proof of the fact that $\lambda^{\ast}:{\cal K}_g\rightarrow\zz$
is a homomorphism, using the Kontsevich integral.

%
%
%
%
%
%
%
%
%
%
%
%
%
%
%
%
%

\vspace{-4mm}
\section*{References}
\sloppy

\vspace{-2mm}
{\small
\begin{bib} 
\item\label{at88} 
\rm M.Atiyah, \it Topological Quantum Field Theories, \rm
Publications Math\'{e}matiques IHES {\bf 68} \rm , 175-186 (1988)

\item\label{BN} 
\rm D.Bar-Natan, \it On the Vassiliev knot invariants, \rm
Topology {\bf 34}, 423-472 (1995)


\item\label{BL}
\rm D.Bar-Natan, R.Lawrence, \it A rational surgery formula for the
LMO invariant, \rm arXiv math.GT/0007045 (2000)

\item\label{BLT}
\rm D.Bar-Natan, T.Le, D.Thurston, \it Two applications of elementary
knot theory to Lie algebras and Vassiliev invariants, \rm Geom. Topol.
{\bf 7}, 1-31 (2003)

\item\label{BHMV}
\rm C.Blanchet, H.Habegger, G.Masbaum, P.Vogel, \it Topological
quantum field theories derived from the Kauffman bracket, \rm
Topology {\bf 34}, no.4, 883-927 (1995)

\item\label{CL}
\rm D.Cheptea, T.Le, \it 3-cobordisms with their rational homology on the boundary, \rm preprint
 
\item\label{CM}
\rm D.Cheptea, G.Massuyeau, \it Tangles, cobordisms, and their LMO-type invariants, \rm in preparation

\item\label{C}
\rm D.Cheptea, \it Universal quantum invariants and the induced representation of the Torelli group, 
\rm in preparation

\item\label{FM}
\rm A.Fomenko, S.Matveev, \it Algorithmic and computer methods for 
three-manifolds, \rm Kluwer Academic Publishers (1997) 

\item\label{gille}
\rm C.Gille, \it On the Le-Murakami-Ohtsuki invariant in degree 2
for several classes of 3-manifolds, \rm J Knot Theory Ramifications
{\bf 12} (1), 17-45 (2003)

\item\label{GS}
\rm R.E.Gompf, A.I.Stipsicz, \it 4-manifolds and Kirby calculus, 
\rm Graduate Studies in Mathematics {\bf 20}, AMS (1999)

\item\label{HM}
\rm N.Habegger, G.Masbaum,\it The Kontsevich integral and Milnor's invariants, 
\rm Topology {\bf 39}, 1253-1289 (2000)

\item\label{HO1}
\rm N.Habegger, K.Orr,\it Finite type three manifold invariants -realization 
and vanishing, \rm  J Knot Theory Ramifications {\bf 8} (8), 1001-1007 
(1999)

\item\label{HO2}
\rm N.Habegger, K.Orr, \it Milnor link invariants and quantum 3-manifold 
invariants, \rm  Comment. Math. Helv. {\bf 74}, no.2, 322-344 (1999)

\item\label{habiro}
\rm K.Habiro, \it Claspers and finite-type invariants of links, \rm  Geometry
and Topology {\bf 4}, 1-83 (2000) 

\item\label{le}
\rm T.T.Q.Le, \it An invariant of integral homology 3-spheres which is
universal for all finite type invariants,\rm AMS Translation series
{\bf 2}, 179, 75-100 (1997)

\item\label{LeGrenoble}
\rm T.T.Q.Le, \it The LMO invariant, \rm ``Invariants de noeuds at de
vari\'{e}t\'{e}s de dimension 3'', \'{E}cole d'\'{e}t\'{e} de Math\'{e}matiques, 
Institut Fourier, Grenoble (1999)

\item\label{LMO}
\rm T.T.Q.Le, J.Murakami, T.Ohtsuki, \it On a universal perturbative
invariant of 3-manifolds, \rm Topology {\bf 37}, no.3, 539-574 (1998)

\item\label{lescop}
\rm C.Lescop, \it Global surgery formula for the Casson-Walker invariant, 
\rm Princeton University Press (1996)

\item\label{matveev}
\rm S.V.Matveev, \it Generalized surgery of three-dimensional 
manifolds and representations of homology spheres, \rm Matematicheskie
Zametki {\bf 42}, no.2, 268-278 (1986)

\item\label{morita}
\rm S.Morita, \it Casson's invariant for homology 3-spheres and 
characteristic classes of surface bundles I, \rm Topology {\bf 28}, no.3, 
305-323 (1989)

\item\label{MO}
\rm J.Murakami, T.Ohtsuki, \it Topological Quantum Field Theory
for the Universal Quantum Invariant, \rm Commun. Math. Phys. {\bf
188}, 501-520 (1997)

\item\label{serre}
\rm J-P.Serre, \it Lie Algebras and Lie Groups, \rm 2nd ed., Lecture
Notes in Mathematics {\bf 1500}, Springer, New York (1992)

\item\label{turaev}
\rm V.Turaev, \it Quantum Invariants of Knots and 3-Manifolds, \rm
Walter de Gruyter (1994)

\item\label{vogel}
\rm P.Vogel, \it Invariants de type fini, \rm en ``Nouveaux
Invariants en G{\' e}om{\' e}trie et en Topologie'', publi{\' e}
par D. Bennequin, M. Audin, J. Morgan, P. Vogel, Panoramas et
Synth{\` e}ses {\bf 11}, Soci{\' e}t{\' e} Math{\' e}matique de
France, 99-128 (2001)
\end{bib}

\begin{flushleft}
Dorin Cheptea\newline
\sc UFR de Math\'ematique et d'Informatique, Universit\'e Louis Pasteur,\newline
7, rue Ren\'e Descartes, 67084, Strasbourg, France\newline
and\newline
\sc Institute of Mathematics, P.O.Box 1-764, Bucharest, 70700, Romania \newline
e-mail: \tt cheptea@math.u-strasbg.fr

\medskip\rm
Thang T Q Le\newline
\sc School of Mathematics, Georgia Institute of Technology,\newline
Atlanta, GA 30332-0160, USA \newline
e-mail: \tt letu@math.gatech.edu
\end{flushleft}
}

\end{document}